\documentclass[11pt,leqno]
{amsart}
\usepackage{amsmath,epsfig,graphicx,color, float}
\usepackage{subcaption}
\usepackage{relsize}
\usepackage{comment}
\usepackage{natbib}

\numberwithin{table}{section}
\numberwithin{figure}{section}
\numberwithin{equation}{section}

\definecolor{darkblue}{rgb}{.2, 0.2,.8}
\definecolor{darkgreen}{rgb}{0,0.5,0.3}
\definecolor{darkred}{rgb}{.8, .1,.1}

\newcommand{\markcol}{\color{black}}
\newcommand{\bfY}{\vect{Y}}
\newcommand{\bfy}{\vect{y}}
\newcommand{\bfx}{\vect{x}}
\newcommand{\bfX}{\vect{X}}
\newcommand{\bfr}{\vect{r}}

\newcommand{\bfalp}{\vect{\alpha}}

\newcommand{\bfpi}{\vect{\pi}}
\newcommand{\bfT}{\mat{T}}
\newcommand{\bft}{\vect{t}}
\newcommand{\bfe}{\vect{e}}
\newcommand{\bfR}{\mat{R}}

\newcommand{\0}{\mat{0}}

\newcommand{\R}{\mathbb{R}}
\newcommand{\E}{\mathbb{E}}
\renewcommand{\P }{{\mathbb P}}

\newcommand{\ov}{\overline}

\newcommand{\ex}{{\rm e}}

\newcommand{\eqd}{\stackrel{d}{=}}

\newtheorem{example}{Example}[subsection]

\newtheorem{remark}{Remark}[subsection]
\newtheorem{algorithm}{Algorithm}[subsection]

\usepackage{bm}
\newcommand{\vect}[1]{\pmb{#1}}
\newcommand{\mat}[1]{\boldsymbol{\bm #1}}

\DeclareMathOperator*{\argmax}{arg\,max}

\allowdisplaybreaks

\begin{document}
\title[Fitting IPH distributions to data]{Fitting inhomogeneous phase-type distributions to data: the univariate and the multivariate case}
\thanks{ The research of Jorge Yslas is supported by Danmarks Frie Forskningsfond Grant No 9040-00086B.}

\author[H.Albrecher ]{Hansj\"org Albrecher}
\address{Department of Actuarial Science, Faculty of Business and Economics, University of Lausanne and Swiss Finance Institute, UNIL-Dorigny,
1015 Lausanne}
\email{hansjoerg.albrecher@unil.ch}

\author[M. Bladt]{Mogens Bladt}
\address{Department  of Mathematics,
University of Copenhagen,
Universitetsparken 5,
DK-2100 Copenhagen,
Denmark}
\email{Bladt@math.ku.dk}
\author[J. Yslas]{Jorge Yslas}
\address{Department  of Mathematics,
University of Copenhagen,
Universitetsparken 5,
DK-2100 Copenhagen,
Denmark}
\email{jorge@math.ku.dk}

\begin{abstract}
The class of inhomogeneous phase-type distributions (IPH) was recently introduced in \cite{albrecher2019inhomogeneous} as an extension of the classical phase-type (PH) distributions. {\markcol Like PH distributions, the class of IPH is dense in the class of distributions on the positive halfline, but}
 leads to more parsimonious models in the presence of heavy tails. In this paper we propose a fitting procedure for this class to given data. We furthermore consider an analogous extension of Kulkarni's multivariate phase-type class \citep{kulkarni1989new} to the inhomogeneous framework and study parameter estimation for the resulting new and flexible class of multivariate distributions. As a by-product, we amend a previously suggested fitting procedure for the homogeneous multivariate phase-type case and provide appropriate adaptations for censored data. The performance of the algorithms is illustrated in several numerical examples, both for simulated and real-life insurance data.
\end{abstract}
\maketitle

\section{Introduction}\setcounter{equation}{0}
The development, study and fitting of flexible distributions for random phenomena is an important branch of applied probability and statistics. Some respective approaches are based on a nice blend of theory and practice, among which the class of phase--type (PH) distributions is a prominent example. Originally initiated by \cite{neuts75}, the realization of a (univariate) phase--type distributed random variable is interpreted as the time until absorption of a  time--homogeneous, finite state--space Markov jump process with one absorbing state and the rest being transient. The explicit description through matrix exponentials makes the resulting class of distributions at the same time versatile and analytically tractable (see e.g. \cite{Bladt2017} for a recent survey). The class of phase--type distributions is known to be dense (in the sense of weak convergence) among all distributions on the positive halfline, but for distributions whose shape is very different from combinations of exponential components (which are the building blocks of the probabilistic Markov jump process construction), a suitable phase--type approximation will need a large dimension of the involved matrix (representing the number of phases of the underlying Markov process) and -- in addition to computational challenges -- may then be seen unnatural. This is particularly the case for heavy--tailed distributions, where the focus in modelling often lies on the tail of the distribution, and the latter is not well captured by the combination of exponential components of the PH construction. After some first amendment involving infinite--dimensional matrices was suggested in \cite{bladt2015calculation} and \cite{bladt2017fitting}, recently a new way to circumvent this problem was proposed in \cite{albrecher2019inhomogeneous}. Concretely, when the Markov jump process is allowed to be  time--inhomogeneous, one gains a lot of flexibility in terms of the structure of the individual components entering the matrix framework, which can reduce the complexity of appropriate fitting distributions drastically, in particular for distributions with heavy tails. The intensity matrices of the Markov jump process are then a function of time. In the general case, they may not commute at different time epochs, which complicates their statistical estimation due to a lack of appropriate sufficient statistics. However, there is an important sub--class for which the intensity matrices can be written as a constant matrix scaled by some real function. In this class all matrices commute, and it was shown in \cite{albrecher2019inhomogeneous} that along this way one in fact obtains, for instance, Pareto, Weibull and Generalized Extreme Value (GEV) distributions with matrix-valued parameters. These distribution classes are all dense in the class of distributions on the positive halfline and inherit the computational advantages of the PH-type class, but also provide excellent fits for heavy--tailed data already for small dimensions, something that the original PH class could not achieve. In particular, if by some preliminary exploratory analysis one has a good guess for an appropriate scaling function (typically suggested by the empirical tail behavior), the resulting matrix distributions can be very parsimonious yet effective model improvements of the respective base distributions with a genuinely heavy tail. However, while parameter estimation for univariate PH distributions by a standard maximum likelihood procedure based on an EM algorithm has been studied in the seminal paper of \cite{asmussen1996fitting} (see also the later extension of \cite{olsson1996estimation} dealing with censored observations and \cite{bladt2003estimation} for an MCMC approach), parameter estimation for the time--inhomogeneous case has not yet been addressed. {\markcol PH distributions have played a crucial role in various application areas beyond finance and insurance, such as queueing theory, biology and operational research. The time--inhomogeneous extension is relatively new and we are confident that it  may prove equally useful in those areas, in particular due to the flexible tail behaviour of this new class.   } \\

Motivated by the flexibility of the approach, in this paper we will also consider an inhomogeneous extension of the  multivariate version of the PH distribution. The multivariate phase--type   distribution (of MPH$^*$ type) was originally introduced by \cite{kulkarni1989new} and is constructed as the joint distribution of certain state-dependent accumulated rewards earned on the same underlying Markov jump process. It has PH-distributed marginals and also enjoys a denseness property in the class of all distributions on the respective positive orthant. Multivariate phase--type distributions have found applications in diverse  areas. For instance,  \cite{cai2005conditional} consider them for determining conditional tail expectations in risk management, 
\cite{cai2005multivariate} studied several types of ruin probabilities for a multivariate compound Poisson risk model when the claim size vector follows an MPH$^*$ distribution, and \cite{herbertsson2011modelling}  used this class to model default contagion in credit risk. 
More recently, \cite{bladt2019parisian} applied MPH$^*$ distributions for the calculation of Parisian type ruin probabilities. In terms of fitting of the (time--homogeneous) multivariate MPH$^*$ distribution, \cite{ahlstrom1999parametric} introduced an algorithm for a bivariate subclass of MPH$^*$, and an EM algorithm for parameter estimation in the general case was proposed in \cite{breuer2016semi}. However, the latter was not actually implemented and contains an inconsistency in the maximum likelihood estimator (which we amend in this paper). 

The inhomogeneous extension of the MPH$^*$ to be proposed in this paper will then again serve the purpose of keeping the dimension of the involved matrices low when one faces a non-exponential behavior in the marginals and the joint  multivariate behaviour. We would like to point out that an alternative analytically tractable deviation from exponential behavior utilizing Mittag-Leffler distributions in both the univariate and multivariate case can be found in \cite{albrecher2019matrix,abb2020multivariate,multiml}. A number of commonly used heavy-tailed multivariate distributions are in fact transformed multivariate exponential distributions. For instance, \cite{mardia1962multivariate} was the first to systematically study multivariate Pareto distributions, which he introduced by transforming a Wicksell--Kibble--type multivariate exponential distribution (see \cite{kibble1941two}). He also noticed that estimation methods for the multivariate exponential can then be translated directly towards the estimation of the multivariate Pareto distribution. \cite{arnold2015pareto} presents some approaches to extend Mardia's analysis to obtain other multivariate distributions with Pareto marginals. Likewise, multivariate versions of the Weibull distribution have been obtained as power transforms of multivariate exponential distributions, see e.g.\ \cite{lee1979multivariate}. The inhomogeneous MPH$^*$ extension that we propose in this paper can to some extent be seen as a generalization and unification of these above models. \\

The main purpose of this paper is to provide algorithms for the statistical fitting of all these flexible classes of distributions and illustrate and discuss their implementation. We will present a unified maximum--likelihood based approach to fitting phase--type distributions (PH), inhomogeneous phase--type distributions (IPH), multivariate phase--type distributions (of MPH$^*$ type) and its newly introduced inhomogeneous extension. These classes contain a large number of mathematically tractable distributions that are sufficiently general to fit any non--negative data set, in the body and for both light or heavy tails. We will also consider extensions of the procedures to adapt for censored data and to the fitting of theoretically known joint distributions. \\

The structure of the paper is as follows. In Section~\ref{sec:pre} we provide an overview of the class of IPH distributions and present a new fitting procedure, which we then exemplify on two particular cases, one on a simulated data set and the other on actual data for lifetimes of the Danish population. In Section~\ref{sec:mph} we shortly recollect some facts about the MPH* class, review existing methods for parameter estimation and provide a substantiation and correction of an algorithm that was previously proposed in the literature. We then extend the algorithm to the case of censored observations, and give more details on an important particular bivariate subclass with explicit density. The section finishes with illustrations of the algorithms for a simulated bivariate sample as well as a phase--type approximation to a known bivariate exponential distribution. In Section~\ref{sec:miph} we introduce some multivariate extensions to distributions in the IPH class, derive basic properties, provide an EM algorithm for its parameter estimation and again illustrate its use in several examples, including multivariate matrix--Pareto models, multivariate matrix--Weibull models as well as a real data application to a bivariate Danish fire insurance data set. Section~\ref{sec:conclusions} concludes.

\section{Inhomogeneous phase--type distributions}\label{sec:pre}
\subsection{Preliminaries}
Let $ \{ J_t \}_{t \geq 0}$ denote a time--inhomogeneous Mar\-kov jump process on a state space $\{1, \dots, p, p+1\}$, where states $1,\dots,p$ are transient and state $p+1$ is absorbing. Then $ \{ J_t \}_{t \geq 0}$ has an intensity matrix of the form
\begin{align*}
	\mat{\Lambda}(t)= \left( \begin{array}{cc}
		\bfT(t) &  \bft(t) \\
		\0 & 0
	\end{array} \right)\,, \quad t\geq0\,,
\end{align*}
where $\bfT(t) $ is a $p \times p$ matrix and $\bft(t)$ is a $p$--dimensional column vector. Here, for any time $t\geq0$, $\bft (t)=- \bfT(t) \, \bfe$, where $\bfe $ is the $p$--dimensional column vector of ones. Let $ \pi_{k} = \P(J_0 = k)$, $k = 1,\dots, p$, $\bfpi = (\pi_1 ,\dots,\pi_p )$ and assume that $\P(J_0 = p + 1) = 0$. Then we say that the time until absorption 
\begin{align*}
	\tau = \inf \{ t \geq  0 \mid J_t = p+1 \}
\end{align*}
has an inhomogeneous phase--type distribution with representation $(\bfpi,\bfT(t) )$ and we write $\tau \sim \mbox{IPH}(\bfpi,\bfT(t) )$.  If $\bfT(t) = \lambda(t)\,\bfT$, where $\lambda(t)$ is some known non--negative real function and $\bfT$ is a sub--intensity matrix, then we write $\tau \sim  \mbox{IPH}(\bfpi , \bfT , \lambda )$. {\markcol Note that for  $\lambda(t) \equiv 1$ one returns to the time-homogeneous case, which corresponds to the conventional phase--type distribution with notation $\mbox{PH}(\bfpi , \bfT )$ (a comprehensive account of phase--type distributions can be found in \cite{Bladt2017})}.
If $X \sim  \mbox{IPH}(\bfpi , \bfT , \lambda )$, then there exists a function $g$ such that \begin{equation}\label{gtrans}
X \sim g(Y) \,,
\end{equation}where $Y \sim \mbox{PH}(\bfpi , \bfT )$.
 Specifically, $g$ is defined by
\begin{equation*}
g^{-1}(x) = \int_0^x \lambda (t)dt   \label{eq:transformation-g}
\end{equation*}
or, equivalently,
\begin{equation*}
\lambda (t) = \frac{d}{dt}g^{-1}(t) \,.  \label{eq:transformation-lambda} 
\end{equation*}
The density $f_X$ and distribution function $F_X$ for $X \sim  \mbox{IPH}(\bfpi , \bfT , \lambda )$ are given by
\begin{eqnarray*}
 f_X(x) &=& \lambda (x)\, \vect{\pi}\exp \left( \int_0^x \lambda (t)dt\ \mat{T} \right)\vect{t} \,, \label{eq:dens-IPH} \\
 F_X(x)&=&1- \vect{\pi}\exp \left( \int_0^x \lambda (t)dt\ \mat{T} \right)\vect{e} \,.  \label{eq:cdf-IPH}
\end{eqnarray*}
For further reading on inhomogeneous phase--type distributions and motivations for their use in modelling we refer to \cite{albrecher2019inhomogeneous}. {\markcol For the representation of some IPH distributions, we make use of functional calculus. 
 If $h$ is an analytic function and $\mat{A}$ is a matrix, we define
\begin{align*}
	h( \mat{A})=\dfrac{1}{2 \pi i} \oint_{\gamma}h(z) (z \mat{I} -\mat{A} )^{-1}dz \,,
\end{align*}
where $\gamma$ is a simple path enclosing the eigenvalues of $\mat{A}$ (cf.\ \cite[Sec. 3.4]{Bladt2017} for details). 
Another standard way to define an analytic matrix-valued function is using the corresponding series expansion, while for non-analytic but sufficiently smooth functions $h(\mat{A})$ can be defined using the Jordan decomposition of $\mat{A}$. We refer to \cite{higham2008functions} for these and further equivalent ways to define $h(\mat{A})$.}


As illustrated in \cite{albrecher2019inhomogeneous}, a number of IPH distributions can be expressed as classical distributions with matrix-valued parameter. 
Important examples include the transformation $g(y)=\beta \left( e^y-1\right) $ for $\beta>0$  in \eqref{gtrans} leading to a  matrix--Pareto distribution with density function and survival function 
\begin{gather}\label{denspar}
f_{X}(x) = \bfpi \left( \dfrac{x}{\beta}+1\right)^{\bfT - \mat{I}} \bft \;\frac{1}{\beta}\,, \quad \bar{F}_{X}(x)=1- {F}_{X}(x) = \bfpi \left(\dfrac{x}{\beta}+1\right)^{\bfT} \bfe \,,
\end{gather} 
respectively, as well as the matrix--Weibull distribution with density and survival function 
\begin{gather*}\label{denswei}
f_{X}(x) = \bfpi e^{\bfT x^{\beta}} \bft \beta x^{\beta-1} \,,\quad \bar{F}_{X}(x) = \bfpi e^{\bfT x^{\beta}} \bfe \,,
\end{gather*} 
obtained from $g(y)=y^{1/\beta}$ ($\beta>0$), see \cite{albrecher2019inhomogeneous} for further details.

\subsection{Parameter estimation}\label{sec:fit-to-IPH}
For the matrix--Pareto distribution \eqref{denspar} and $\beta = 1 $, the transform is parameter-independent, so that the distribution can be fitted to i.i.d.\ data $x_1,\dots,x_N$ by fitting a phase--type distribution $\mbox{PH}(\vect{\pi},\mat{T})$ to the transformed data $\log (1+x_1),\dots,\log (1+x_N)$ using an EM algorithm \citep{asmussen1996fitting}. This was the procedure employed in \cite{albrecher2019inhomogeneous} for the numerical illustration there. The general case -- where the transform does depend on parameters -- is more subtle and shall be dealt with here. The key will be to apply a parameter-dependent transformation in each step of the EM algorithm. \\

Let $x_1, \dots, x_N$ be an i.i.d.\ sample of an inhomogeneous phase--type distribution with representation $X\sim \mbox{IPH}(\bfpi , \bfT, $ $ \lambda (\, \cdot\, ;\vect{\beta} ) )$, where $ \lambda ( \,\cdot\, ; \vect{\beta} )$ is a parametric non--negative function depending on the vector $\vect{\beta}$. We then know that $X \eqd g(Y; \vect{\beta})$ with $Y \sim \mbox{PH}(\bfpi , \bfT )$ and $g$ is defined in terms of its inverse function $g^{-1}(x; \vect{\beta})=\int_{0}^{x}\lambda(t; \vect{\beta})dt$. In particular $g^{-1}( X; \vect{\beta}) \eqd Y \sim \mbox{PH}(\bfpi , \bfT )$. The EM algorithm for fitting  $\mbox{IPH}(\bfpi , \bfT,  \lambda (\, \cdot\, ;\vect{\beta} ) )$ then works as follows.

\begin{algorithm} [EM algorithm for transformed phase--type distributions]\label{alg:transPH} \
	
	0. Initialize with some ``arbitrary'' $( \bfpi,\bfT , \vect{\beta} )$.
	
	1. Transform the data into $y_i=g^{-1}(x_i; \vect{\beta})$, $i=1,\dots,N$, and 
	apply the E-- and M--steps of the conventional EM algorithm of \cite{asmussen1996fitting} by which we obtain the estimators $( \hat{\bfpi},\hat{\bfT})$.
	
	2. Compute  
	\begin{align*}
	\hat{\vect{\beta}}  & = \argmax_{\vect{\beta}} \sum_{i=1}^{N} \log (f_{X}(x_i; \hat{\bfpi}, \hat{\bfT},\vect{\beta} )) 
	= \argmax_{\vect{\beta}} \sum_{i=1}^{N} \log \left( \lambda(x_i ; \vect{\beta}) \hat{\bfpi} \exp\left({ \int_{0}^{x_i} \lambda(t ; \vect{\beta})  dt \ \hat{\bfT} }\right) \hat{\bft} \right) \,.
	\end{align*}
	
	3. Assign $(\vect{\pi},\mat{T},\vect{\beta}) =(\hat{\vect{\pi}},\hat{\mat{T}}, \hat{\vect{\beta}})$ and GOTO 1. 
\end{algorithm}
Then the likelihood function increases for each iteration, and hence converges to a (possibly local) maximum.
\begin{proof}
	Since the data points $x_i$ are assumed to be i.i.d.\ realisations from the unknown distribution $\mbox{IPH}(\vect{\pi},\mat{T},\lambda)$, there exists a function $g$ such that $y_i=g^{-1}(x_i;\vect{\beta})$ are i.i.d.\ realisations of phase--type distributed random variables $\mbox{PH}(\vect{\pi},\mat{T})$. That function $g$ is assumed to be known up to the value of $\vect{\beta}$. In turn, $x_i = g(y_i;\vect{\beta})$, so a data point $x_i$ can be interpreted as the absorption time of the Markov jump process corresponding to $\mbox{PH}(\vect{\pi},\mat{T})$, which is $y_i$,  but with the scale of the time axis for the $y_i$--data converted (stretched) into $g(\,\cdot\, ;\vect{\beta})$--coordinates instead. The full data likelihood is then given by 
	\begin{eqnarray*} 
		L (\vect{\pi},\mat{T},\vect{\beta};\vect{y})&=& \prod_{k=1}^p \pi_k^{B_k} \prod_{k=1}^p\prod_{l\neq k} t_{kl}^{N_{kl}}e^{-t_{kl}Z_k(\vect{\beta})} \prod_{k=1}^p t_k^{N_k}e^{-t_k Z_k(\vect{\beta})}\,,
	\end{eqnarray*}
	where $B_k$ is the number of times the Markov process underlying the phase--type distribution initiates in state $k$, $N_{kl}$ denotes the total number of transitions from state $k$ to $l$, $N_k$ denotes the number of times an exit to the absorbing state was caused by a jump from state $k$, and $Z_k(\vect{\beta})$ is the total time the Markov process has spent in state $k$. We notice that $Z_k(\vect{\beta})$ is the only sufficient statistic which depends on the transformation of the time axis for the $y$--data and hence on $\vect{\beta}$. Consequently, for any given $\vect{\beta}$, the $E$--step is simply the one as in \cite{asmussen1996fitting}, and so is the $M$--step for $(\vect{\pi},\mat{T})$.

	The $\vect{\beta}$ update in 2. requires a general, usually numerical, maximization of the incomplete data likelihood. 
    Each iteration of the algorithm increases the likelihood. Indeed, let $L^I$ denote the incomplete data likelihood{\markcol , i.e.,
    \begin{align*}
	L^{I}(\vect{\pi},\mat{T},\vect{\beta} ; \bfx)  & = \prod_{i=1}^{N}  f_{X}(x_i; {\bfpi}, {\bfT},\vect{\beta} )
	= \prod_{i=1}^{N} \lambda(x_i ; \vect{\beta}) {\bfpi} \exp\left({ g^{-1}(x_i;\vect{\beta} )  \ {\bfT} }\right) {\bft}  \,,
	\end{align*}
    and consider  parameter values $(\vect{\pi}_n,\mat{T}_n,\vect{\beta}_n)$ after the $n$-th iteration. In the $(n+1)$-th iteration, we first obtain $(\vect{\pi}_{n+1},\mat{T}_{n+1})$ in 1. so that
    \begin{align*}
	\prod_{i=1}^{N}  \bfpi_{n} \exp\left({ g^{-1}(x_i;\vect{\beta}_n )  \ \bfT_{n} }\right) \bft_{n} \leq \prod_{i=1}^{N}  \bfpi_{n+1} \exp\left({ g^{-1}(x_i;\vect{\beta}_n )  \ \bfT_{n+1} }\right) \bft_{n+1}  \,,
	\end{align*}
   }
    By monotonicity of $g$ and the transformation theorem, 
    \[  L^I(\vect{\pi}_n,\mat{T}_n,\vect{\beta}_n;\vect{x})\leq   L^I(\vect{\pi}_{n+1},\mat{T}_{n+1},\vect{\beta}_n;\vect{x})   \]
    and hence, by 2., 
    \[  L^I(\vect{\pi}_n,\mat{T}_n,\vect{\beta}_n;\vect{x})\leq   L^I(\vect{\pi}_{n+1},\mat{T}_{n+1},\vect{\beta}_n;\vect{x})\leq  \sup_{\vect{\beta}}L^I(\vect{\pi}_{n+1},\mat{T}_{n+1},\vect{\beta};\vect{x}) =  L^I(\vect{\pi}_{n+1},\mat{T}_{n+1},\vect{\beta}_{n+1};\vect{x}) \,. \]
\end{proof}

\begin{example}\normalfont{(Matrix--Gompertz)}
Let $X= \log( \beta Y  + 1 ) /\beta$, where $Y \sim \mbox{PH}(\bfpi , \bfT )$ and $\beta > 0$. Then
	\begin{equation}\label{gomp}
		\bar{F}_{X}(x) = \bfpi e^{\bfT (e^{\beta x} -1) / \beta} \bfe  \quad \text{and}\quad 
		f_{X}(x) = \bfpi e^{\bfT (e^{\beta x} -1) / \beta} \bft  e^{\beta x} \,.
	\end{equation}
We refer to the distribution of $X$ as a \textit{matrix--Gompertz distribution}, since the scale parameter of the usual Gompertz distribution is now replaced by a matrix. Note that the resulting distribution has a lighter  tail than a conventional phase--type distribution. The Gompertz distribution is used in a number of applications, most notably it is historically used for the modelling of human lifetimes \citep{gompertz}. Its matrix version \eqref{gomp} provides a natural flexible extension.
As an illustration, we fitted a matrix--Gompertz distribution with 3 phases using Algorithm~\ref{alg:transPH} with $2\,500$ iterations to the lifetime of the Danish population that died in the year $2\,000$ at ages $50$ to $100$ (data obtained from the Human Mortality Database (HMD) and available in the R-package \textit{MortalitySmooth} \citep{mortaitys}).
Here and in later examples, the number of iterations in the algorithm is chosen in such a way that the changes in the successive log--likelihoods become negligible. Concerning running times, our implementation makes use of the gradient ascent method for the maximization part of the algorithm, in which the running times highly depend on the step--length and the actually chosen stopping criterion. In the present example we employed a step--length of $10^{-8}$ and run gradient ascent until the absolute value of the derivative is less than $0.001$ leading to a running time of about $35$ seconds on a usual PC (with $2.9$ GHz Dual--Core Intel Core i$5$ processor $5287$U) for the $2\,500$ iterations of the EM algorithm. Note that this  choice of stopping criterion is to prioritize precision over speed, and an improvement on running times can be attained by using a different maximization procedure.
The obtained parameters are as follows: 
\begin{gather*} 
 \hat{\bfpi}=\left(
0.0450, \,0.1303,\, 0.8246 \right)\,, \\ 
\hat{\bfT}=\left( \begin{array}{ccc}
-0.1357  & 0.1214 & 0  \\ 
0.0130  & -0.0421 & 0.0288 \\
0.1415 & 0.0184 & -0.1620  
\end{array} \right) \,, \\ 
\hat{\beta}= 0.1019 \,.
\end{gather*}

Figure~\ref{fig:exGompertz} shows that the fitted density recovers the structure of the data quite well. Note that conventional phase--type distributions have been used to model the distribution of lifespans (see for instance \cite{asmussen2019phase}). However, the number of phases required to capture the tail behavior of the data with the latter is rather large, due to the lighter than exponential tail. In contrast, the matrix--Gompertz distribution provides an excellent fit with comparably fewer parameters (phases). 

\begin{figure}[H]
	\centering
	\begin{subfigure}{0.22\textwidth}
		\includegraphics[width=\textwidth]{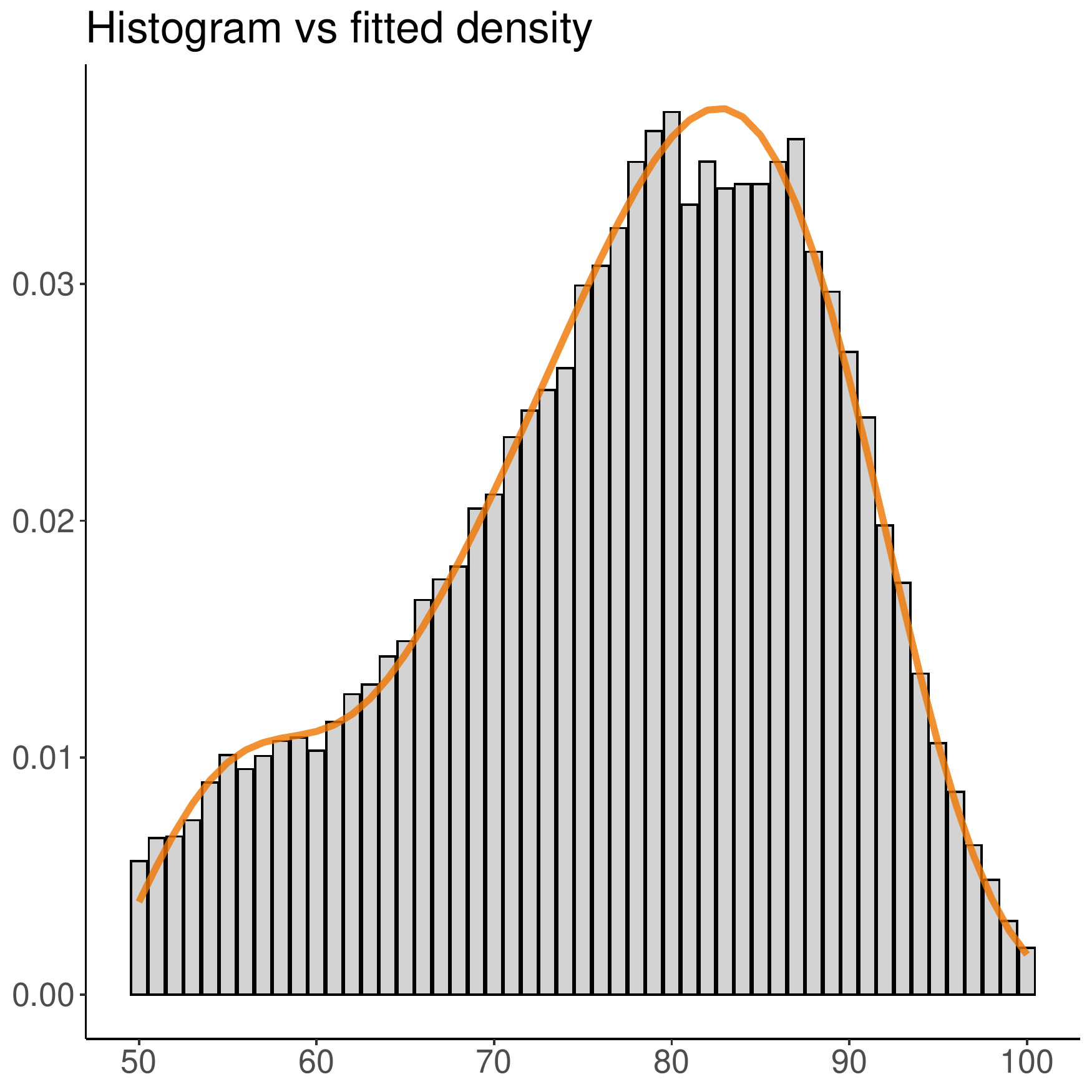}
	\end{subfigure}
	\caption{Histogram of lifetimes of the Danish population that died in the year $2\,000$ at ages $50$ to $100$ versus the density of the fitted matrix--Gompertz distribution. } \label{fig:exGompertz}
\end{figure}
\end{example}

\begin{example}\normalfont{(Matrix--GEV)}
	Algorithm~\ref{alg:transPH} can also be applied to estimate distributions that are not IPH in a strict sense, but that are defined as a transformation of a PH distribution. This is for instance the case for $g(y) = \mu - \sigma(y^{-\xi} - 1)/\xi$ with $\mu \in \R$, $\sigma >0$, and $\xi \in \R$. Recall from \cite{albrecher2019inhomogeneous} that
	\begin{gather*}
		{F}_{X}(x) =\left\{\begin{array}{ll}
\bfpi \exp \left(\bfT \left(1+\xi \dfrac{x-\mu}{\sigma}\right)^{-1 / \xi}\right) \bfe\,, & \xi \neq 0\,, \\
\bfpi \exp \left(\bfT \exp \left(-\dfrac{x-\mu}{\sigma}\right)\right)\bfe\,, & \xi=0\,,
\end{array}\right.  \\
		f_{X}(x) =\left\{\begin{array}{ll}
\dfrac{1}{\sigma} \bfpi \exp \left(\bfT \left(1+\xi \dfrac{x-\mu}{\sigma}\right)^{-1 / \xi}\right) \bft \left(1+\xi \dfrac{x-\mu}{\sigma}\right)^{-(1 + \xi) / \xi} \,, & \xi \neq 0\,, \\
\dfrac{1}{\sigma} \bfpi \exp \left(\bfT \exp \left(-\dfrac{x-\mu}{\sigma}\right)\right)\bft \exp \left(-\dfrac{x-\mu}{\sigma}\right) \,, & \xi=0\,,
\end{array}\right.  
	\end{gather*} 
	from which it becomes clear that this distribution can be interpreted as a matrix version of the generalized extreme value (GEV) distribution, see e.g.\ \cite{beirlant2004}. As an illustration, we generated an i.i.d.\ sample of size $5\,000$ from such a distribution of $3$ phases with parameters 
	\begin{gather*} 
	{\bfpi}=\left(
	1, \,0,\, 0\right)\,, \\ 
	{\bfT}=\left( \begin{array}{ccc}
	-1 & 0.5 & 0  \\
	0.2 & -2 & 0.8  \\
	1 & 1 & -5 \\
	\end{array} \right) \,, \\ 
	\mu= 2 \,,\quad \sigma = 0.5 \,,\quad \xi =0.4 \,,
	\end{gather*}
	which has theoretical moments $\E(X)=2.2524$ and $SD(X)=1.4423$. The generated sample has moments  $\hat{\E}(X)=2.2607$ and $\hat{SD}(X)=1.3307$. We then fitted such a matrix--GEV distribution with the same number of phases using Algorithm~\ref{alg:transPH} with $1\,500$ steps, obtaining the following parameters: 
	\begin{gather*} 
	\hat{\bfpi}=\left(
	0.0772 , \,0.1268,\, 0.7960 \right)\,, \\ 
	\hat{\bfT}=\left( \begin{array}{ccc}
	-8.9772 & 0.0964 & 0.0001  \\
	0.2891 & -2.8439 & 0.3542  \\
	3.2353 & 0.0137 & -5.7731 \\
	\end{array} \right) \,, \\ 
	\hat{\mu}= 1.3852 \,,\quad \hat{\sigma} = 0.2285 \,,\quad \hat{\xi} = 0.4251 \,.
	\end{gather*}
	We observe that the algorithm estimates pretty  well the shape parameter $\xi$, which determines the heaviness of the tail. Moreover, the fitted distribution has moments $\E(X)=2.2640$ and $SD(X)=1.6587$, which resemble the ones of the sample, and Figure~\ref{fig:exGEVD} shows that the algorithm recovers both body and tail of the data. Note also that the log--likelihood of the fitted matrix--GEV is $-4\,104.541$, while the log--likelihood using the original matrix--GEV distribution is $-4\,107.005$. 
	Such a comparison of the log--likelihoods works as additional evidence for the performance of the algorithm. 
    One can observe that the parameters estimated for $\vect{\pi}$, $\mat{T}$, $\mu$ and $\sigma$ do not resemble the original parameter values, but this is linked with the well-known identifiability issue for phase--type distributions (namely that other parameter combinations may lead to a very similar density shape). In fact, the algorithm finds the parameters that maximize the likelihood for the given sample, and as the concrete numbers above show, the present parameters even outperform the original model underlying the sample(!), see also the convincing QQ--plot in Figure~\ref{fig:exGEVD}. Here, the step--length is $10^{-5}$ and the gradient ascent is run until the norm of the derivative is less than $0.1$ leading to a running time of $2\,609$ seconds for the $1\,500$ iterations of Algorithm~\ref{alg:transPH}.
	
	\begin{figure}[H]
		\centering
		\begin{subfigure}{0.22\textwidth}
			\includegraphics[width=\textwidth]{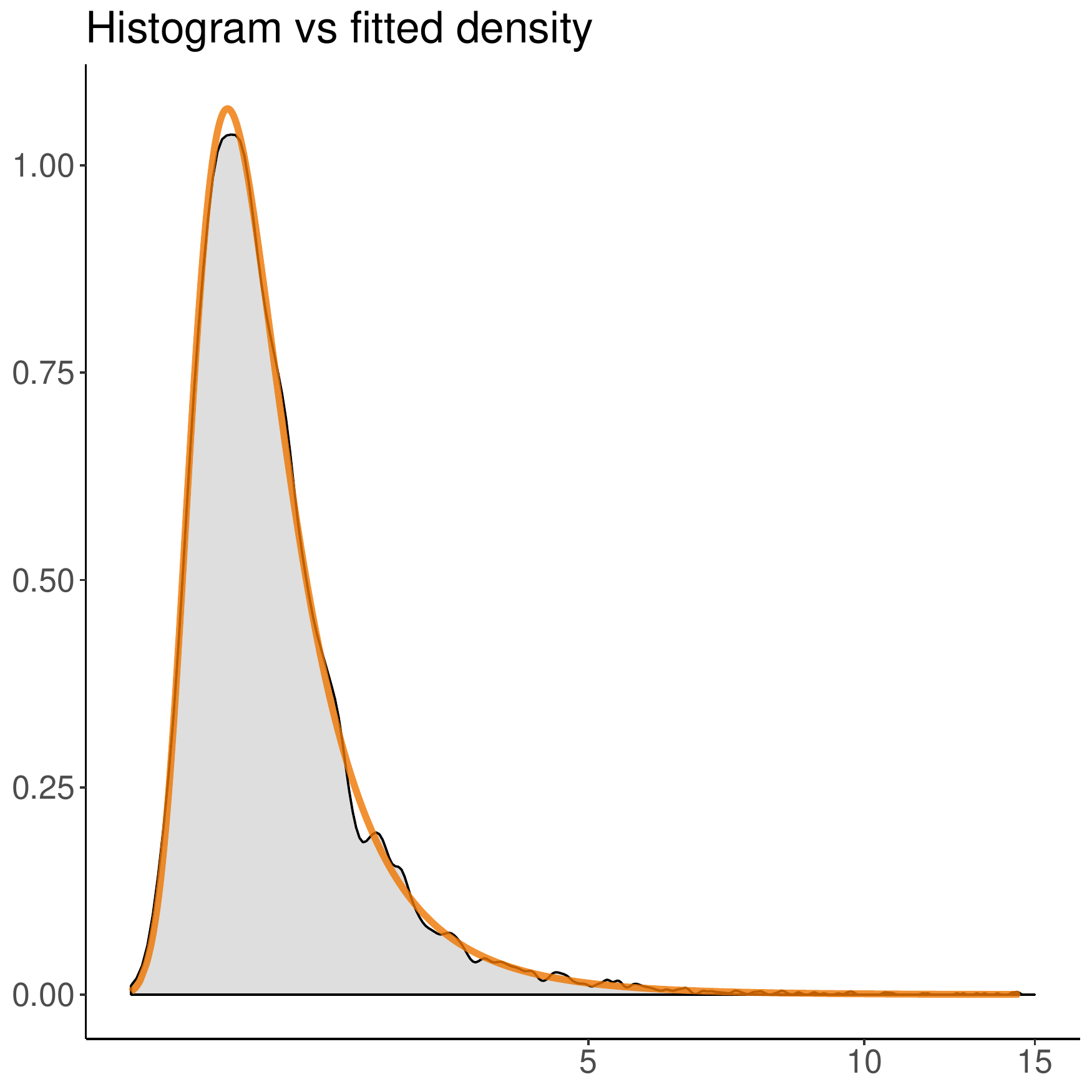}
		\end{subfigure}
		\begin{subfigure}{0.22\textwidth}
			\includegraphics[width=\textwidth]{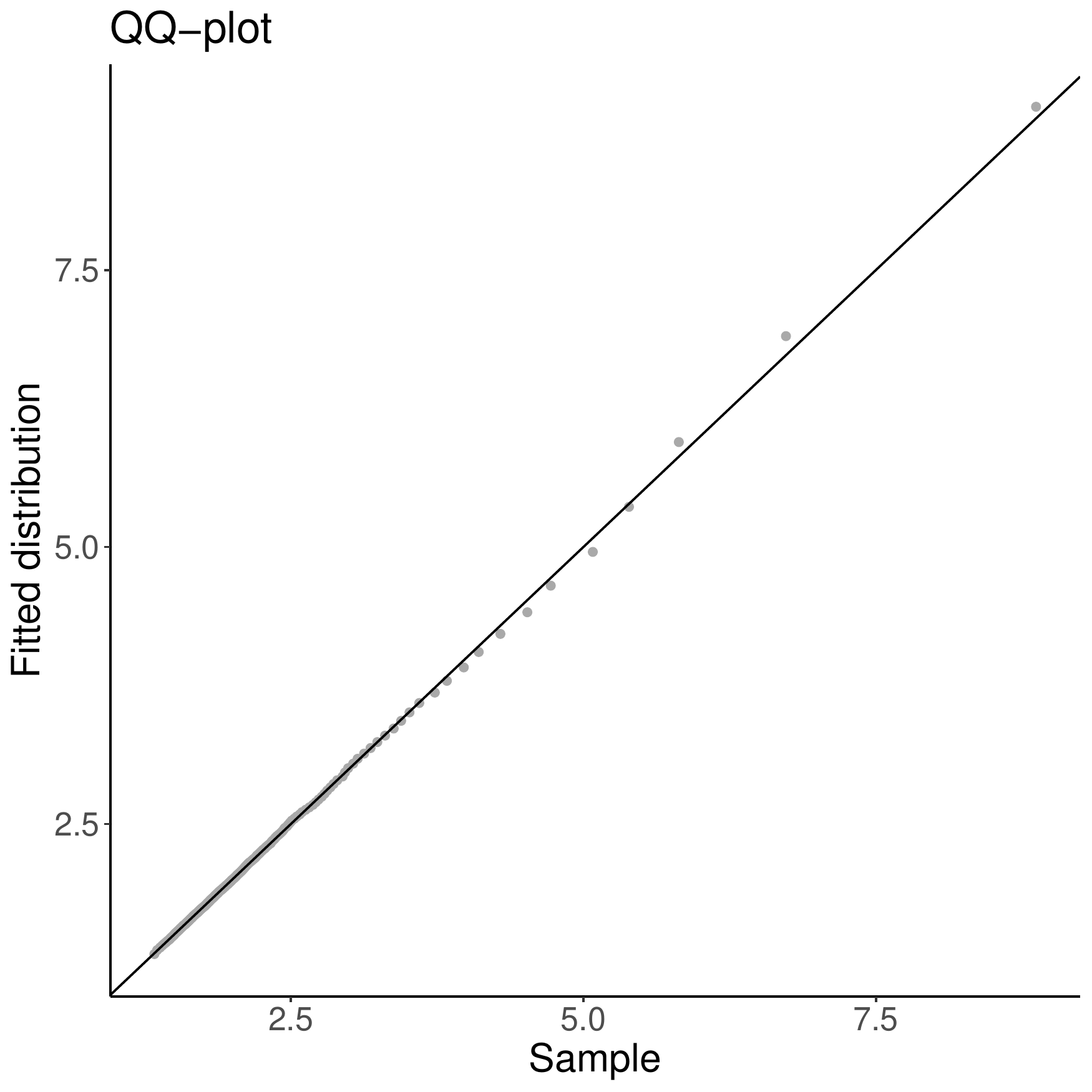}
		\end{subfigure}
		\caption{ Histogram of simulated sample versus density of the fitted matrix--GEV distribution in log--scale (left) as well as QQ--plot of simulated sample versus fit (right). } \label{fig:exGEVD}
	\end{figure}
\end{example}
\section{Multivariate phase-type distributions}\label{sec:mph}
\subsection{Preliminaries}Let $\tau \sim \mbox{PH} \left( \bfpi, \bfT \right)$ be a (conventional) $p$--dimensional phase--type distributed random variable with underlying time--homogeneous Markov jump process $\{ J_{t} \}_{t\geq 0} $. Let $\bfr_{j}=\left( r_{j}(1),\dots, r_{j}(p) \right)^{\prime}$ be non--negative $p$--dimensional column vectors, $j=1,\dots,d$, and let
\begin{align*}
\bfR=\left( \bfr_{1}, \bfr_{2}, \dots, \bfr_{d}  \right)
\end{align*}
be a $p \times d$--dimensional \textit{reward matrix}. 
 Now define 
\begin{align*}
Y^{(j)}=\int_{0}^{\tau}r_{j}\left(J_{t} \right)dt
\end{align*} 
for all $j=1,\dots,d$. If we interpret 
$r_{j}(k)$ as the rate at which a reward is obtained while $J_t$ is in state $k$, then
 $Y^{(j)}$ is the total reward for component $j$ obtained prior to absorption. We then say that the random vector  $\bfY = \left(Y^{(1)},\dots, Y^{(d)} \right)^{\prime}$ has a multivariate phase--type distribution of the MPH$^*$ type (as defined in  \cite{kulkarni1989new}, see also \cite{Bladt2017}) and we write $\bfY \sim 	\mbox{MPH}^{*}\left( \bfpi, \bfT, \bfR \right)$.  

While each member of the MPH$^*$ class has an explicit expression for the (joint) Laplace transform and the joint moments of any order (see Section 8.1.1 of \cite{Bladt2017}), there are no general explicit expressions 
for the density and distribution functions. However, for certain structures and sub--classes explicit expressions for the latter do exist (like Example 8.1.13 of \cite{Bladt2017}). 

If $\bfY=\left(Y^{(1)},\dots, Y^{(d)} \right)^{\prime} \sim \mbox{MPH}^{*}\left( \bfpi, \bfT, \bfR \right)$, then each marginal $Y^{(j)}$ has a phase--type distribution, $\mbox{PH}(\vect{\pi}_j,\mat{T}_j)$ say. First we decompose
\[  \vect{r}_j= \begin{pmatrix}\vect{r}_j^+ \\ \vect{r}_j^0 \end{pmatrix} 	\,, \quad  \vect{\pi}= \left( \vect{\pi}^+ \ \vect{\pi}^0 \right) \quad \mbox{and} \quad \mat{T} = 
\begin{pmatrix}
 \mat{T}^{++} & \mat{T}^{+0} \\
 \mat{T}^{0+} & \mat{T}^{00}
 \end{pmatrix}, \]
 where we have reordered the state space such that the $+$ terms correspond to the states $k$ for which the rewards $r_k^{(j)
 }$ are strictly positive, and the $0$ terms to the states with zero rewards. E.g., $\mat{T}_{+0}$ corresponds to the intensities by which the underlying Markov jump process $\{ J_t \}_{t\geq 0}$ jumps from a state with positive reward to a state with zero reward. Then the phase--type distribution of $Y^{(j)}$ is given by an atom 
 at zero of size $\vect{\pi}_0 (\mat{I}-(-\mat{T}_{00})^{-1}\mat{S}_{0+})\,\vect{e}$, where $\vect{e}$ is the column vector of ones of appropriate dimension, and
 \begin{equation}
  \vect{\pi}_j = \vect{\pi}^+ + \vect{\pi}^0(-\mat{T}^{00})^{-1}\mat{T}^{0+} \ \ \mbox{together with} \ \ \mat{T}_j =\mat{\Delta}(\vect{r}_j^+)^{-1}\left(  \mat{T}^{++} + \mat{T}^{+0}(-\mat{T}_{00})^{-1}\mat{T}^{0+} \right) \,,  \label{eq:T_j}\end{equation}
 where $\Delta(\vect{a})$ denotes the $d' \times d'$ diagonal matrix with entries $a^{(m)}$, $m=1,\dots,d'$, from a $d'$--dimensional vector $\vect{a}$. 
 The atom appears in case there is a positive probability of starting in a non--reward--earning state (0) and the underlying Markov process gets absorbed before visiting a reward earning state (+). The Markov jump process generating $Y^{(j)}$ starts 
 in the same state as $J_t$ if the reward is positive (hence $\vect{\pi}_+$) or it starts in the first state with positive rewards that $J_t$ enters after starting in a zero reward state (hence the term $\vect{\pi}^0(-\mat{T}^{00})^{-1}\mat{T}^{0+} $). Similar arguments apply to the generator $\mat{T}_j$, where only reward-earning terms will form part of the state space for $Y^{(j)}$. We refer to \cite{Bladt2017} for further details.

 Summarizing, each marginal $Y^{(j)}$ has a phase--type distribution, which is based on the original Markov process $\{ J_t\}_{t\geq 0}$, but with a possibly smaller state space and with rescaled parameters.

\subsection{Parameter estimation}\label{subsec:mphest}
We next provide an algorithm for estimating MPH$^*$ distributed data. The data consist of a $d$--dimensional multivariate sample of $N$ i.i.d.\ observations 
\begin{align*}
	 \bfy_i =(y_i^{(1)},\dots, y_i^{(d)})^{\prime}\,, \quad i=1, \dots, N \,.
\end{align*}
That is, we only observe the times to absorption, $y_i^{(j)}$, of each phase--type distributed marginal. Hence we are clearly in an incomplete data set--up and we shall employ the EM algorithm for fitting $(\vect{\pi},\mat{T},\mat{R})$.  

The EM algorithm works by replacing unavailable sufficient statistics by their conditional expectations given data under given parameters, and thereby updating the parameters by using known formulas for the maximum likelihood estimator in the complete data domain. Iteration of the procedure then produces a sequence of parameter values which increases the likelihood in each step. 

For the present situation, we define the complete data as both the trajectories of the underlying Markov process which generates the phase--type distribution from which the marginals of the multivariate vector are constructed, and the Markov jump processes representing the rewards in all marginal distributions. It is not sufficient with complete knowledge of the marginal trajectories only. Indeed, one can easily construct examples where the underlying processes cannot be reconstructed from the marginals only. The complete knowledge of both marginals and the underlying Markov process which generates the marginals creates another problem in relation to the incomplete data since we do not have observations for the absorption times of the underlying Markov process. We can get around this problem by assuming that the rows of the reward matrix $\mat{R}$ sum to one, i.e. $\mat{R}\vect{e}=\vect{e}$. This assumption is not restrictive and can be imposed without losing generality due to the great ambiguity of (multivariate) phase--type representations.  
Hence our data consists of marginals $ \bfy^{(j)} =(y_1^{(j)},\dots, y_N^{(j)})^{\prime}$, $j=1,\dots,d$, and their sums $\bfy^{(S)} =  \sum_{l=1}^{d} \bfy^{(l)} $.

In the complete data domain, the estimation is straightforward and works as follows.
{\markcol Using the notation introduced in the proof of Algorithm \ref{alg:transPH},} the complete data MLE for $(\vect{\pi},\mat{T})$ is given by
\[  \hat{\pi}_k =\frac{B_k}{N}\,, \quad \hat{t}_{kl}=\frac{N_{kl}}{Z_k}\,, \quad  \hat{t}_k = \frac{N_k}{Z_k}\,, \quad \hat{t}_{kk}=-\sum_{l\neq k}\hat{t}_{kl}-\hat{t}_k \,.   \]
The rewards of the marginals are then given by
\[  \hat{r}_j(k) = \frac{Z_k^{(j)}}{Z_k} = \frac{Z_k^{(j)}}{\sum_{l=1}^{d} Z_k^{(l)}}, \]
where $\markcol Z_k^{(j)} = r_{j}(k) Z_k $ is the over-all amount of time the $j$'th component has spent in state $k$. 

In the EM algorithm, we now must replace all aforementioned sufficient statistics by their conditional expectations given data. Concerning $Z_k$, $N_{kl}$, $N_k$ and $B_k$, these only depend on the underlying Markov jump process and are computed conditionally on $\vect{y}^{(S)}$ only. 
Their formulas are then as stated in the algorithm below (see \cite{asmussen1996fitting}). 

Concerning the conditional expectation of $Z_k^{(j)}$, we must calculate the expected reward (under $(\vect{\pi},\mat{T})$) given all data of marginal $j$, which amounts to calculating the conditional expected time given data for the corresponding phase--type representation of the $j$-th marginal, $(\vect{\pi}_j,\mat{T}_j)$. These are readily given by (again using \cite{asmussen1996fitting}) 
\begin{align*}
	\E\left(Z_{k}^{(j)} \mid \bfY^{(j)}=\bfy^{(j)} \right)=\sum_{i=1}^{N} \dfrac{\mathlarger{\int_{0}^{y^{(j)}_{i}} \bfe_{k}^{\prime}  \ex^{\bfT_{j}(y^{(j)}_{i}-u)  }\bft_{j} \bfpi_{j} \ex^{  \bfT_{j}u  } \bfe_{k} du }}{\bfpi_{j} \ex^{ \bfT_{j}y^{(j)}_{i}  }\bft_{j}}\,.
\end{align*}
Then
\begin{align*}
	\hat{r}_j(k)=\dfrac{\E\left(Z_{k}^{(j)} \mid \bfY^{(j)}=\bfy^{(j)} \right)}{\E\left(Z_{k} \mid \bfY^{(S)}=\bfy^{(S)} \right)} .
\end{align*}
Iterating the above finally provides a (single) full EM algorithm for the estimation of $(\vect{\pi},\mat{T},\mat{R})$. 
We summarize the results in the following. 

\begin{algorithm}[EM algorithm for MPH* distributions]\label{alg:MEM}\

	0. Initialize with some ``arbitrary'' $(\bfpi,\bfT , \mat{R})$ with $\mat{R}\vect{e}=\vect{e}$, and compute $\vect{\pi}_j$ and $\mat{T}_j$, $j=1,\dots,d$, using \eqref{eq:T_j}.
	
	1. (E--step) Calculate
	\begin{align*}
		& \E\left( B_{k} \mid \bfY^{(S)}=\bfy^{(S)} \right)=\sum_{i=1}^{N} \dfrac{ \pi_{k} \bfe_{k}^{\prime} \ex^{ \bfT y_{i}^{(S)} }\bft}{\bfpi \ex ^ { \bfT y_{i}^{(S)} }\bft }
		\\
		& \E\left( Z_{k} \mid \bfY^{(S)}  = \bfy^{(S)} \right) = \sum_{i=1}^{ N }  \dfrac {\mathlarger{\int_{0}^{y_{i}^{(S)}} \bfe_{k}^{\prime} \ex^{\bfT (y_{i}^{(S)}-u)} \bft \bfpi \ex^{\bfT u} \bfe_{k} du}} { \bfpi \ex^{\bfT y_{i}^{(S)} } \bft } 
		\\
		& \E \left( N_{kl} \mid \bfY^{(S)}  = \bfy^{(S)} \right) = \sum_{i=1}^{N}t_{kl}  \dfrac{\mathlarger{\int_{0}^{y_{i}^{(S)}} \bfe_{l}^{\prime} \ex^{\bfT (y_{i}^{(S)}-u)} \bft \bfpi \ex^{ \bfT u }} \bfe_{k} du} { \bfpi \ex^{\bfT y_{i}^{(S)}} \bft }
		\\
		&\E \left( N_{k} \mid \bfY^{(S)}  = \bfy^{(S)} \right) = \sum_{i=1}^{N} t_{k} \dfrac{ \bfpi \ex^{ \bfT y_{i}^{(S)}} \bfe_{k}} { \bfpi \ex^{ \bfT y_{i}^{(S)} } \bft }  \\
		&\E\left(Z_{k}^{(j)} \mid \bfY^{(j)}=\bfy^{(j)} \right)=\sum_{i=1}^{N} \dfrac{\mathlarger{\int_{0}^{y^{(j)}_{i}} \bfe_{k}^{\prime}  \ex^{\bfT_{j}(y^{(j)}_{i}-u)  }\bft_{j} \bfpi_{j} \ex^{  \bfT_{j}u  } \bfe_{k} du }}{\bfpi_{j} \ex^{ \bfT_{j}y^{(j)}_{i}  }\bft_{j}}\,.
	\end{align*}
	
	2. (M--step) Let
	\begin{align*}
		&\hat{\pi}_{k} = \dfrac{1}{N} \E\left( B_{k} \mid \bfY^{(S)} = \bfy^{(S)} \right)\,,\quad
		\hat{t}_{kl} = \dfrac{ \E\left( N_{kl} \mid \bfY^{(S)}  = \bfy^{(S)} \right) } { \E \left( Z_{k} \mid \bfY^{(S)}  = \bfy^{(S)} \right) }\,,\quad
		\hat{t}_{k} = \dfrac{ \E\left( N_{k} \mid \bfY^{(S)}  = \bfy^{(S)} \right) } { \E \left( Z_{k} \mid \bfY^{(S)}  = \bfy^{(S)} \right) } \,,\\
		&\hat{t}_{kk} = -\sum_{l\neq k} \hat{t}_{kl} - \hat{t}_ {k}\,, \quad \hat{\bfpi} = ( \hat{\pi}_{1} , \ldots , \hat{\pi}_{p} )\,,\quad \hat{\bfT} = \{ \hat{t}_{kl} \}_{ k,l = 1,\ldots,p}\ \ \mbox{and} \ \ \hat{\bft} = ( \hat{t}_{1},\ldots,\hat{t}_{p} )^{\prime}\,.
	\end{align*}
	and
	\begin{align*}  
	\hat{r}_j(k):= \dfrac{\E\left(Z_{k}^{(j)} \mid \bfY^{(j)}=\bfy^{(j)} \right)}{\mathlarger{ \E\left(Z_{k} \mid \bfY^{(S)}=\bfy^{(S)} \right)}} = \dfrac{\E\left(Z_{k}^{(j)} \mid \bfY^{(j)}=\bfy^{(j)} \right)}{\mathlarger{\sum_{l=1}^{d} \E\left(Z_{k}^{(l)} \mid \bfY^{(l)}=\bfy^{(l)} \right)}} \quad \mbox{and} \quad  \hat{\bfR}=\{ \hat{r}_j(k) \}_{ k=1,\dots, p, j=1,\dots, d}\,.
\end{align*}
	
	3. Assign $\bfpi:=\hat { \bfpi }$, $ \bfT :=\hat { \bfT }$, $ \bft :=\hat{\bft}$, $\bfR:=\hat{\bfR}$  and compute $\vect{\pi}_j$, $\mat{T}_j$, $j=1,\dots,d$, using \eqref{eq:T_j}.  GOTO 1.\\
\end{algorithm}

\begin{remark}\rm  
Algorithm \ref{alg:MEM} was originally proposed in \cite{breuer2016semi} as two consecutive EM algorithms and its original statement contained a minor error in the M--step update for the reward matrix. 
 To see why Algorithm~\ref{alg:MEM} can be decomposed into  the two consecutive EM algorithms, we argue as follows. Running the EM Algorithm~\ref{alg:MEM}, $(\hat{\vect{\pi}},\hat{\mat{T}})$ will eventually converge (without input from the part involving the reward components). For constant  $(\hat{\vect{\pi}},\hat{\mat{T}})$, Algorithm~\ref{alg:MEM} is indeed equivalent to the second EM algorithm in \cite{breuer2016semi}. More specifically, the algorithm takes the following form.

 {\normalfont\textbf{First EM.}}

 	0. Initialize with some ``arbitrary'' $(\bfpi,\bfT)$.
	
 	1. (E--step) Calculate
 	\begin{align*}
 		& \E\left( B_{k} \mid \bfY^{(S)}=\bfy^{(S)} \right)=\sum_{i=1}^{N} \dfrac{ \pi_{k} \bfe_{k}^{\prime} \ex^{ \bfT y_{i}^{(S)} }\bft}{\bfpi \ex ^ { \bfT y_{i}^{(S)} }\bft }
 		\\
 		& \E\left( Z_{k} \mid \bfY^{(S)}  = \bfy^{(S)} \right) = \sum_{i=1}^{ N }  \dfrac {\mathlarger{\int_{0}^{y_{i}^{(S)}} \bfe_{k}^{\prime} \ex^{\bfT (y_{i}^{(S)}-u)} \bft \bfpi \ex^{\bfT u} \bfe_{k} du}} { \bfpi \ex^{\bfT y_{i}^{(S)} } \bft } 
 		\\
 		& \E \left( N_{kl} \mid \bfY^{(S)}  = \bfy^{(S)} \right) = \sum_{i=1}^{N}t_{kl}  \dfrac{\mathlarger{\int_{0}^{y_{i}^{(S)}} \bfe_{l}^{\prime} \ex^{\bfT (y_{i}^{(S)}-u)} \bft \bfpi \ex^{ \bfT u }} \bfe_{k} du} { \bfpi \ex^{\bfT y_{i}^{(S)}} \bft }
 		\\
 		&\E \left( N_{k} \mid \bfY^{(S)}  = \bfy^{(S)} \right) = \sum_{i=1}^{N} t_{k} \dfrac{ \bfpi \ex^{ \bfT y_{i}^{(S)}} \bfe_{k}} { \bfpi \ex^{ \bfT y_{i}^{(S)} } \bft } \, .
 	\end{align*}
	
 	2. (M--step) Let
 	\begin{align*}
 		&\hat{\pi}_{k} = \dfrac{1}{N} \E\left( B_{k} \mid \bfY^{(S)} = \bfy^{(S)} \right)\,,\quad
 		\hat{t}_{kl} = \dfrac{ \E\left( N_{kl} \mid \bfY^{(S)}  = \bfy^{(S)} \right) } { \E \left( Z_{k} \mid \bfY^{(S)}  = \bfy^{(S)} \right) }\,,\quad
 		\hat{t}_{k} = \dfrac{ \E\left( N_{k} \mid \bfY^{(S)}  = \bfy^{(S)} \right) } { \E \left( Z_{k} \mid \bfY^{(S)}  = \bfy^{(S)} \right) }\,,\\
 		&\hat{t}_{kk} = -\sum_{l\neq k} \hat{t}_{kl} - \hat{t}_ {k}\,, \quad \hat{\bfpi} = ( \hat{\pi}_{1} , \ldots , \hat{\pi}_{p} )\,,\quad \hat{\bfT} = \{ \hat{t}_{kl} \}_{ k,l = 1,\ldots,p}\quad \mbox{and} \quad \hat{\bft} = ( \hat{t}_{1},\ldots,\hat{t}_{p} )^{\prime}\,.
 	\end{align*}
	
 	3. Assign $\bfpi:=\hat { \bfpi }$, $ \bfT :=\hat { \bfT }$, $ \bft :=\hat{\bft}$ and GOTO 1.\\
	
 \noindent {\normalfont\textbf{Second EM.}} Use the estimated $(\bfpi,\bfT )$ of the first EM.

 0. Initialize with some ``arbitrary'' $\bfR$ with $\mat{R}\vect{e}=\vect{e}$, and compute $\vect{\pi}_j$ and $\mat{T}_j$, $j=1,\dots,d$, using \eqref{eq:T_j}.

 1. (E--step) Calculate
 \begin{align*}
 	\E\left(Z_{k}^{(j)} \mid \bfY^{(j)}=\bfy^{(j)} \right)=\sum_{i=1}^{N} \dfrac{\mathlarger{\int_{0}^{y^{(j)}_{i}} \bfe_{k}^{\prime}  \ex^{\bfT_{j}(y^{(j)}_{i}-u)  }\bft_{j} \bfpi_{j} \ex^{  \bfT_{j}u  } \bfe_{k} du }}{\bfpi_{j} \ex^{ \bfT_{j}y^{(j)}_{i}  }\bft_{j}}\,.
 \end{align*}	
 
 	2. (M--step) Let
 \begin{align*}
 	\hat{r}_j(k):= \dfrac{\E\left(Z_{k}^{(j)} \mid \bfY^{(j)}=\bfy^{(j)} \right)}{\mathlarger{\sum_{l=1}^{d} \E\left(Z_{k}^{(l)} \mid \bfY^{(l)}=\bfy^{(l)} \right)}} \quad \mbox{and} \quad  \hat{\bfR}=\{ \hat{r}_j(k) \}_{ k=1,\dots, p, j=1,\dots, d}\,.
 \end{align*}
	
 3. Assign  $\bfR:=\hat{\bfR}$ and compute $\vect{\pi}_j$ and $\mat{T}_j$, $j=1,\dots,d$, using \eqref{eq:T_j}. GOTO 1.

\end{remark}

 \begin{remark} \rm 
The main computational burden lies in the E--steps, where matrix exponentials and integrals thereof must be evaluated. 
 In \cite{asmussen1996fitting}  this is done by converting the problem into a system of ODEs, which are then solved via a Runge--Kutta method of fourth order (a C implementation, called EMpht, is available online \citep{olsson1998empht}).  While this approach is adequate for fitting univariate phase--type distributions, the Runge--Kutta method fails to work in some cases in the multivariate setting, in particular for the second EM, when an element in the reward matrix approaches zero. The reason is that the sub-intensity matrix of (at least) one of the marginals will adjust to this change by increasing some of the entries of the matrix in each iteration, and thus requiring an increasingly smaller step--size in the Runge--Kutta method to accurately approximate the solution to the system. 
 Our implementation includes an approach for the computation of matrix exponentials based on uniformization, and it is a slight variation of the method in \cite[p.232]{neuts1995algorithmic}. We explain briefly the method. By taking $\phi= \max (-t_{kk})_{k=1,\dots,p}$ and defining $\mat{P} := {\phi}^{-1}\left(  \phi \mat{I} + \bfT \right)$, which is in fact a transition matrix, we have that
 \begin{align*}
 	\exp( \bfT y)= \sum_{n=0}^{\infty} \frac{(\phi y)^{n}}{n !} e^{-\phi y }  	\mat{P} ^{n} \,. 
 \end{align*}
 Then
 \begin{align*}
 	\left| e^{\bfT y} - \sum_{n=0}^{ M } \frac{(\phi y)^{n}}{n !} e^{-\phi y }  	\mat{P} ^{n}  \right| \leq  \sum_{n=M+1}^{\infty} \frac{(\phi y)^{n}}{n !} e^{-\phi y }  	  \left|  \mat{P} ^{n}  \right| \leq \sum_{n=M+1}^{\infty} \frac{(\phi y)^{n}}{n !} e^{-\phi y } = \mathbb{P}(N_{\phi y}>M) \,,
 \end{align*}
 where $N_{\phi y}$ is Poisson distributed with mean ${\phi y}$. Hence, we can find $M$ such that the difference of the matrix exponential with a finite sum is less than or equal to a given error $\epsilon>0$. Of course, larger values of $\phi y$ give bigger values of $M$, dismissing any computational improvement for large observations. A way to circumvent this problem is to observe that $e^{\bfT y}=(e^{\bfT y/2^{m}})^{2^m}$, thus we can find $m$ such that $\phi y /2^{m} <1 $, compute $e^{\bfT y/2^{m}}$ by a finite sum and then retrieve $e^{\bfT y}$ by squaring.
 
 To compute the integrals involving matrix exponentials, we observe that by defining 
	\begin{align*}
		\mat{G} (y; \bfpi, \bfT ):= \int_{0}^{y} \ex^{ \bfT (y-u)} \bft \bfpi \ex^{ \bfT u}du \, ,  
	\end{align*}
	we have that (see \cite{van1978computing})
	\begin{align*}
		\exp \left( \left(  \begin{array}{cc}
		\bfT & \bft \, \bfpi \\
		\0 & \bfT
	\end{array}  \right) y \right)=  \left(  \begin{array}{cc}
		\ex^{\bfT y } & \mat{G} (y; \bfpi, \bfT ) \\
		\0 & \ex^{\bfT y } 
	\end{array}  \right)\,.
	\end{align*}
	Correspondingly, a simple (and efficient) way to compute $\mat{G} (y; \bfpi, \bfT )$ is by calculating the matrix exponential of the left hand side.
	
 Approaches to improve the speed of the EM algorithm in the univariate case exist in the literature; for instance, \cite{okamura2011refined} proposed a method also based on uniformization. \hfill $\Box$
 \end{remark}

\subsection{Parameter estimation for censored data}\label{subsub:censoring}
In certain applications, some or all of the data may be censored. We call a data point \textit{right--censored} at $v$ if it takes an unknown value above $v$, \textit{left--censored} at $w$ if it takes an unknown value below $w$, and  \textit{interval--censored} if it is contained in the interval $(v , w]$, but its exact value is unknown. Left--censoring is a special case of interval--censoring with $v=0$, while right--censoring can be obtained by fixing $v$ and letting $w \to \infty$.

The EM Algorithm \ref{alg:MEM} works much in the same way as for uncensored data, with the only difference that we are no longer observing exact data points $Y^{(j)} = y^{(j)}$,   but only $Y^{(j)} \in (v^{(j)},w^{(j)}]$. This will only change the E--steps, where the conditional expectations can be calculated using the formulas in \cite{olsson1996estimation}. We now explain in detail how to adapt Algorithm \ref{alg:MEM} to censored data. \\

\noindent\textbf{ First EM}

 It is possible that a data point consists of a combination of marginals with both censored (not necessarily in the same intervals) and uncensored data (this is relevant in the  first EM algorithm when considering data of the sum of the marginals). 
Table~\ref{tb:censored} contains all possible combinations one might have in the data and the way of treating them. Note that for $d>2$, one simply repeats the same rules iteratively.

\begin{table}[H]
\begin{center}
\begin{tabular}{|l|l|l|}
\hline
  $Y^{(1)}$& $Y^{(2)}$ & $Y^{(S)}=Y^{(1)}+Y^{(2)}$ \\
  \hline
  Uncensored  with value $y^{(1)}$ & Uncensored  with value  $y^{(2)}$ & Uncensored   with value $y^{(1)}+y^{(2)}$ \\[2mm]
  Right--censored at $v^{(1)}$ & Uncensored with value $y^{(1)}$ & Right--censored at $v^{(1)}+y^{(2)}$ \\[2mm]
  Right--censored at $v^{(1)}$ & Right--censored at $v^{(2)}$ & Right--censored at $v^{(1)}+v^{(2)}$ \\[2mm]
  Right--censored at $v^{(1)}$ & Interval--censored $(v^{(2)},w^{(2)}]$ & Right--censored at $v^{(1)}+v^{(2)}$ \\[2mm]
  Interval--censored $(v^{(1)},w^{(1)}]$ & Uncensored with value $y^{(2)}$ & Interval--censored $(v^{(1)}+y^{(2)},w^{(1)}+y^{(2)}]$ \\[2mm]
  Interval--censored $(v^{(1)},w^{(1)}]$  & Interval--censored $(v^{(2)},w^{(2)}]$  & Interval--censored $(v^{(1)}+v^{(2)},w^{(1)}+w^{(2)}]$\\
  \hline
\end{tabular}
\vspace{0.5cm}
\caption{Rules for censored data.}
\label{tb:censored}
\end{center}
\end{table}
For completeness, we include here the conditional expectations needed (see also \cite{olsson1996estimation}). 

\begin{align*}
		& \E\left( B_{k} \mid Y^{(S)} \in (v,w] \right) = \frac{ \pi_{k} \bfe_{k}^{\prime} \ex^{\bfT v}\bfe - \pi_{k} \bfe_{k}^{\prime} \ex^{\bfT w}\bfe} { \bfpi \ex^{\bfT v}\bfe - \bfpi \ex^{\bfT w}\bfe}\,,
		\\
		& \E\left( Z_{k} \mid Y^{(S)} \in (v,w] \right) =  \dfrac{\mathlarger{ \int_{v}^{w} \bfpi \ex^{\bfT u} \bfe_{k}du - \left( \int_{0}^{w} \bfe_{k}^{\prime} \ex^{\bfT(w-u)} \bfe \bfpi \ex^{\bfT u} \bfe_{k}du - \int_{0}^{v} \bfe_{k}^{\prime} \ex^{\bfT(v-u)} \bfe \bfpi \ex^{\bfT u} \bfe_{k}du \right)}} {\bfpi \ex^{\bfT v}\bfe - \bfpi \ex^{\bfT w}\bfe }\,,
		\\
		& \E \left( N_{kl} \mid Y^{(S)} \in (v,w] \right) =  t_{kl}  \dfrac {\mathlarger{\int_{v}^{w}  \bfpi \ex^{\bfT u} \bfe_{k}du - \left( \int_{0}^{w} \bfe_{l}^{\prime} \ex^{\bfT(w- u)} \bfe \bfpi \ex^{\bfT u} \bfe_{k} du  - \int_{0}^{v} \bfe_{l}^{\prime} \ex^{\bfT(v- u)} \bfe \bfpi \ex^{\bfT u} \bfe_{k} du \right) }} { \bfpi \ex^{\bfT v}\bfe - \bfpi \ex^{\bfT w}\bfe}\,,
		\\
		&\E\left( N_{k} \mid Y^{(S)} \in (v,w] \right) = t_{k} \dfrac{\mathlarger{ \int_{v}^{w} \bfpi \ex^{\bfT u} \bfe_{k} du}} { \bfpi \ex^{ \bfT v }\bfe - \bfpi \ex^{\bfT w}\bfe}\,. 
\end{align*}

\noindent\textbf{ Second EM} 

The second EM algorithm works as above, with the only difference that for marginals with censored data the corresponding conditional expectation is calculated as 
\begin{eqnarray*}
	\lefteqn{\E\left( Z_{k}^{(j)} \mid Y^{(j)} \in (v^{(j)},w^{(j)}] \right) } \\
	& &= \frac{ \mathlarger{\int_{v^{(j)}}^{w^{(j)}}  \bfpi_j \ex^{\bfT_j u} \bfe_{k} du - \left( \int_{0}^{w^{(j)}} \bfe_{k}^{\prime} \ex^{\bfT_j (w^{(j)}-u)} \bfe \bfpi_j \ex^{\bfT_j u} \bfe_{k}du  - \int_{0}^{v^{(j)}} \bfe_{k}^{\prime} \ex^{\bfT_j (v^{(j)}-u)} \bfe \bfpi_j \ex^{\bfT_j u} \bfe_{k}du \right)}} {\bfpi_j \ex^{\bfT_j v^{(j)}}\bfe - \bfpi_j \ex^{\bfT_j w^{(j)}}\bfe } \,.
\end{eqnarray*}

\subsection{A bivariate phase--type distribution with explicit density}\label{subsec:bivmph}
 For a general $\mbox{MPH}^{*}$ distribution an explicit density is not available. \cite{kulkarni1989new} characterized the density by a system of partial differential equations, and in \cite{breuer2016semi} a semi--explicit form is deduced. The following type of bivariate phase--type distributions does lead to an explicit density:\\ Let $\bfY=(Y^{(1)},Y^{(2)})^{\prime} \sim \mbox{MPH}^{*}( \bfpi, \bfT, \bfR) $ with
\begin{align} \label{par:biv}
\bfT= \left( \begin{array}{cc}
\bfT_{11} & \bfT_{12} \\ 
\mat{0} & \bfT_{22}
\end{array} \right)\,, \quad \bfpi= \left(
\bfalp ,\, \vect{0} \right) \quad \text{and }\quad \bfR= \left( \begin{array}{cc}
\bfe & \mat{0} \\ 
\mat{0} & \bfe
\end{array} \right) \,,
\end{align}
where $\bfT_{11}$ and $\bfT_{22}$ are sub--intensity matrices of dimensions $p_1$ and $p_2$ ($p=p_1+p_2$), respectively, and   $\bfT_{11}\,\bfe + \bfT_{12} \,\bfe = \vect{0}$. Then the joint density of $\bfY$ is given by 
\begin{align}\label{eq:denbiv}
	f_{\bfY} \left( y^{(1)}, y^{(2)} \right)=\bfalp \ex^{\bfT_{11} y^{(1)}} \bfT_{12}\ex^{\bfT_{22} y^{(2)}}(-\bfT_{22})\,\bfe \,,
\end{align}
with marginals $Y^{(1)} \sim \mbox{PH}( \bfalp, \bfT_{11})$ and $Y^{(2)} \sim \mbox{PH}( \bfalp(-\bfT_{11})^{-1}\bfT_{12}, \bfT_{22})$.  Note that the Baker--type bivariate distributions introduced in \cite{bladt2019parisian} are a particular case. The latter have some remarkable properties: one can construct a distribution of this type with specific given marginals and a given Pearson correlation coefficient; this class is also dense within the set of bivariate distributions with support in $\mathbb{R}_{+}^{2}$  (this follows from the fact that the class of Bernstein copulas can be used to approximate arbitrarily well any copula (see \cite{sancetta2004bernstein}) and that the class of phase--type distributions can approximate arbitrarily well any distribution with support on $\mathbb{R}_{+}$), making the bigger class of bivariate distributions also dense.

\subsubsection{Tail independence}\label{subsubsec:bivtailmph}
 The existence of an explicit form of the density allows us to compute the upper tail dependence coefficient $\lambda_{U}$. Recall that the latter is defined as
\begin{align*}
	\lambda_{U}&= \lim_{q \to 1^{-}}\P\left( Y^{(1)}>F_{Y^{(1)}}^{-1}(q) \mid Y^{(2)}>F_{Y^{(2)}}^{-1}(q) \right). 
\end{align*}
It is a classical measure of dependence in the tail and of considerable interest in applications in insurance and finance, where the modelling of tail events is crucial. From \eqref{eq:denbiv} we have 
\begin{align*}
	\bar{F}_{\bfY} (y^{(1)}, y^{(2)}) = \P\left( Y^{(1)}>y^{(1)}, Y^{(2)}>y^{(2)} \right)=\bfalp \left(- \bfT_{11} \right)^{-1} \ex^{\bfT_{11} y^{(1)}} \bfT_{12}\ex^{\bfT_{22} y^{(2)}}\bfe \,. 
\end{align*}
Then, if $-\lambda_{j}$ is the real part of the eigenvalue of $\bfT_{jj}$ with largest real part and $k_{j}$ is the dimension of the Jordan block of $\lambda_{j}$ for $j=1,2$, it is easy to see that
\begin{align*} 
	\bar{F}_{\bfY} (y^{(1)}, y^{(2)})  \sim b (y^{(1)})^{k_{1}-1} \ex^{-\lambda_{1}y^{(1)}} (y^{(2)})^{k_{2}-1} \ex^{-\lambda_{2}y^{(2)}} \,,\quad \text{as} \quad y^{(1)},y^{(2)} \to \infty \,,
\end{align*}
where $b$ is a positive constant. Hence
\begin{align*}
	\lambda_{U}
	&= \lim_{q \to 1^{-}}\frac{b (F_{Y^{(1)}}^{-1}(q))^{k_{1}-1} \ex^{-\lambda_{1}(F_{Y^{(1)}}^{-1}(q))} (F_{Y^{(2)}}^{-1}(q))^{k_{2}-1} \ex^{-\lambda_{2}(F_{Y^{(2)}}^{-1}(q))}}{c  (F_{Y^{(2)}}^{-1}(q))^{k_{2}-1} \ex^{-\lambda_{2}(F_{Y^{(2)}}^{-1}(q))}}
	=0 \,,
\end{align*}
with $c$ positive constant. In other words, $\bfY$ is upper-tail-independent.  


\subsubsection{Estimation}
The density \eqref{eq:denbiv} allows for a special form of EM algorithm. Such an algorithm  was introduced in \cite{ahlstrom1999parametric} and we include it  for completeness, subsequent use and comparison purposes. 

\begin{algorithm}\label{alg:bivariate}\

	0. Initialize with some ``arbitrary'' $(\bfalp,\bfT )$.
	
	1. (E--step) Calculate
	\begin{align*}
	& \E\left( B_{k} \mid \bfY = \bfy \right) = \sum_{i=1}^{N} \dfrac{\alpha_{k}\bfe_{k}^{\prime} \ex^{\bfT_{11}y^{(1)}_{i}}\bfT_{12} \ex^{\bfT_{22}y^{(2)}_{i}} \left( -\bfT_{22}  \right) \bfe}{f_{\bfY}(y^{(1)}_{i},y^{(2)}_{i}; \bfalp ,\boldsymbol { T } )}\,, \quad k=1,\dots, p_1
		\\[1em]
		& \E\left( Z_{k} \mid \bfY = \bfy \right) = \left \lbrace \begin{array}{ll}
\mathlarger{\sum_{i=1}^{N}} \dfrac{ \mathlarger{\int_{0}^{y^{(1)}_{i}} \bfalp \ex^{\bfT_{11}u}\bfe_{k} \bfe_{k}^{\prime} \ex^{\bfT_{11}(y^{(1)}_{i} - u)} \bfT_{12} \ex^{\bfT_{22}y^{(2)}_{i}} \left( - \bfT_{22} \right) \bfe du}}{f_{\bfY}(y^{(1)}_{i},y^{(2)}_{i}; \bfalp ,\bfT)}\,,  & k=1,\dots, p_{1} \\[2em]
 \mathlarger{\sum_{i=1}^{N}} \dfrac{\mathlarger{ \int_{0}^{y^{(2)}_{i}} \bfalp \ex^{\bfT_{11}y^{(1)}_{i}}\bfT_{12}\ex^{\bfT_{22}u}\bfe_{(k-p_{1})} \bfe_{(k-p_{1})}^{\prime} \ex^{\bfT_{22}(y^{(2)}_{i}-u)} \left( -\bfT_{22} \right) \bfe du}}{f_{\bfY}(y^{(1)}_{i},y^{(2)}_{i}; \bfalp ,\bfT )}\,,  & k=p_{1}+1,\dots, p
\end{array}
  \right.
		\\
		& \E\left( N_{kl} \mid \bfY=\bfy \right) = \\[1em]
		& \quad \left\lbrace\begin{array}{ll}
 \mathlarger{\sum_{i=1}^{N}} t_{kl}  \dfrac {\mathlarger{ \int_{0}^{y^{(1)}_{i}}\bfalp \ex^{\bfT_{11}u }\bfe_{k} \bfe_{l}^{\prime} \ex^{\bfT_{11}(y^{(1)}_{i}-u )} \bfT_{12} \ex^{\bfT_{22}y^{(2)}_{i}} \left( - \bfT_{22}  \right) \bfe du} } { f_{\bfY}(y^{(1)}_{i},y^{(2)}_{i};\bfalp ,\bfT) }\,,  & k,l=1,\dots,p_{1}, \, k \neq l \\[2em]
  \mathlarger{\sum_{i=1}^{N}} t_{kl} \dfrac{\bfalp \ex^{\bfT_{11}y^{(1)}_{i}} \bfe_{k}  \bfe_{l-p_{1}}^{\prime}
	 \ex^{\bfT_{22}y^{(2)}_{i}} \left( -\bfT_{22} \right) \bfe }{f_{\bfY}(y^{(1)}_{i},y^{(2)}_{i}; \bfalp ,\bfT)}\,, & k=1,\dots,p_{1},\, l=p_{1}+1,\dots,p \\[2em]
\mathlarger{ \sum_{i=1}^{N}} t_{kl}  \dfrac { \mathlarger{ \int_{0}^{y^{(2)}_{i}} \bfalp \ex^{\bfT_{11}y^{(1)}_{i}}\bfT_{12} \ex^{\bfT_{22}u } \bfe_{k-p_{1}} \bfe_{l-p_{1}}^{\prime} \ex^{\bfT_{22}(y^{(2)}_{i}- u  )} \left( -\bfT_{22} \right) \bfe  du }} { f_{\bfY}(y^{(1)}_{i},y^{(2)}_{i}; \bfalp ,\bfT) }\,, &  k, l=p_{1}+1,\dots,p, \,k \neq l
\end{array}
 \right.
		\\[1em]
		&\E\left( N_{k} \mid \bfY= \bfy \right) = \sum_{i=1}^{N}t_{k} \dfrac{\bfalp \ex^{\bfT_{11}y^{(1)}_{i}} \bfT_{12} \ex^{\bfT_{22}y^{(2)}_{i} } \bfe_{k-p_{1}}  }{f_{\bfY}(y^{(1)}_{i},y^{(2)}_{i} ; \bfalp ,\bfT)}\,, \quad k=p_1+1,\dots, p \,.
	\end{align*}
	
	2. (M--step) Let
	\begin{align*}
		&\hat{\alpha}_{k} = \dfrac{1}{N} \E\left(B_{k} \mid \bfY = \bfy \right)\,,\quad
		\hat{t}_{kl} = \dfrac{\E\left( N_{kl} \mid \bfY = \bfy \right)}{\E\left( Z_{k} \mid \bfY = \bfy \right) }\,,\quad
		\hat{t}_{k} = \dfrac {\E\left( N_{k} \mid \bfY = \bfy \right) } { \E\left( Z_{k} \mid \bfY = \bfy \right)}\,,\quad
		\hat{t}_{kk} =-\sum_{l\neq k} \hat{t}_{kl} -\hat{t}_{k} \,,\\
		& \hat{\bfalp} = ( \hat{\alpha}_{1}, \ldots, \hat{\alpha}_{p_{1}} )\,,\quad \hat{\bfT} = \{ \hat{t}_{kl} \}_{k, l=1, \ldots, p} \quad \mbox{and} \quad \hat{\bft} = (\hat{t}_{p_1+1}, \ldots, \hat{t}_{p})^{\prime} \,.
	\end{align*}
	
	3. Assign $\bfalp:=\hat{\bfalp}$, $\bfT:=\hat{\bfT}$, $\bft:=\hat{\bft}$ and GOTO 1.

\end{algorithm}

We now provide two detailed illustrations. 
When Algorithm~\ref{alg:MEM} is employed, given that an explicit form of the joint density is not available, we  choose the number of iterations in such a way that the changes in the successive log--likelihoods in the first EM become negligible and the changes in the successive parameter estimates become negligible in the second EM. For Algorithm~\ref{alg:bivariate} we used a criterion similar to the univariate case. {\markcol For a visual assessment of the quality of the fits in the multivariate setting, we also added some contour plots the interpretation of which, however, is more subjective than the one of the corresponding one-dimensional graphs. Nevertheless, we believe that they do provide some helpful insight concerning the shape of the distributions.}
   
\begin{example}[Simulation study]\label{ex:mph4}\normalfont 
The objective of the present example is to compare the performance of Algorithm~\ref{alg:MEM} and Algorithm~\ref{alg:bivariate}. 
We will illustrate that the more general Algorithm~\ref{alg:MEM} also provides reasonable results when dealing with a sample from a bivariate distribution with density \eqref{eq:denbiv}, for which the more specific Algorithm~\ref{alg:bivariate} is particularly well-suited. We generated an i.i.d.\ sample of size $ 10\,000$ from a $\mbox{MPH}^{*}$ distribution with parameters 
\begin{align*} 
\bfpi=\left( 
 0.15,\, 0.85,\, 0,\, 0 \right)\,, \quad 
\bfT=\left( \begin{array}{cccc}
-2 & 0 & 2 & 0 \\ 
9 & -11 & 0 & 2 \\
0 & 0 & -1 & 0.5 \\
0 & 0 & 0 & -5 
\end{array} \right) \,, \quad \text{and} \quad 
\bfR= \left( \begin{array}{c c}
1 & 0 \\ 
1 & 0  \\
0 & 1 \\
0 & 1
\end{array} \right)\,,
\end{align*}
which has theoretical mean  
	$\E(\bfY)= \left( \E(Y^{(1)}) , \E(Y^{(2)}) \right)^{\prime}= \left(  0.5 , 0.9609 \right)^{\prime} $, 
	 and correlation coefficient $\rho(Y^{(1)},Y^{(2)})=0.1148 $.
Moreover, we know that $\lambda_U=0$. The simulated sample has numerical values $\hat{\E}(\bfY)= \left(  0.5046 , 0.9650 \right)^{\prime} $, 
	 $\hat{\rho}=0.1235 $ and Kendall's tau $\hat{\rho}_\tau=0.1613$. We now use Algorithms \ref{alg:MEM} and \ref{alg:bivariate} to recover the underlying structure of the data, and we assess the quality of the estimation by comparing densities, QQ plots, numerical properties of the distributions and contour plots.

Using Algorithm \ref{alg:MEM} with the same number of phases $p=4$, random initial values and $3\,500$ steps in each EM algorithm, we obtain the following parameters:
\begin{gather*} 
\hat{\bfpi}=\left(
 0.8592 ,\, 0.0002 ,\, 0.0005 ,\, 0.1402 \right)\,, \\ 
\hat{\bfT}=\left( \begin{array}{cccc}
-9.2654 & 0.3805 & 7.2657 & 1.6192 \\ 
0.0039 & -1.1038 & 0.0002 & 0.0118 \\
0.7429 & 0.1640 & -7.1356 & 4.9516 \\
0.3295 & 1.6162 & 0.5278 & -2.4740 
\end{array} \right) \,, \\ 
\hat{\bfR}= \left( \begin{array}{c c}
0.5950 & 0.4050 \\ 
0 & 1  \\
0.4254 & 0.5746 \\
0.8347 & 0.1653
\end{array} \right)\,.
\end{gather*} 
The fitted distribution has mean 
$	\E(\bfY)= \left(
0.5205, 0.9491 \right)^{\prime}$ and
	$\rho=0.2559 $, 
which approximates reasonably well the mean of the original distribution and to a lesser degree well the correlation coefficient. We also approximated $\lambda_U$ and $\rho_\tau$ via simulation  obtaining $\hat{\lambda}_U=0.0058$ and $\hat{\rho}_\tau=0.2843$. Figures \ref{examph:den} and \ref{examph:qq} show that the algorithm is able to recover the structure of the marginals and the sum of the marginals. 
Moreover, Figure \ref{examph:contour} shows that the contour plot of the fitted distribution is similar to the one of the sample. Here, the running time of Algorithm  \ref{alg:MEM} with its $3\,500$ iterations was $924$ seconds.

\begin{figure}[hbt]
	\centering
  	\begin{subfigure}{0.22\textwidth}
  		\includegraphics[width=\textwidth]{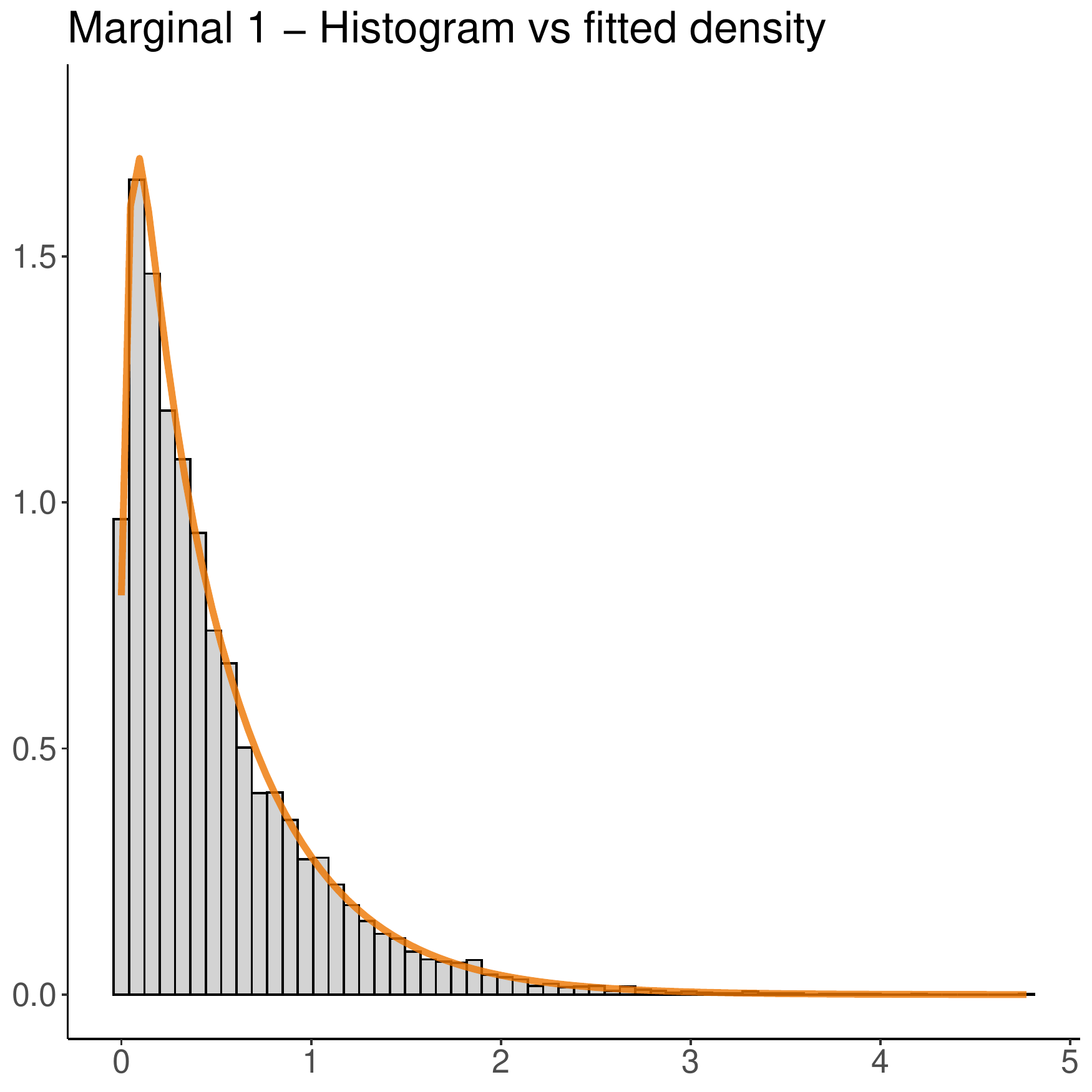}
  	\end{subfigure}
  	\begin{subfigure}{0.22\textwidth}
  		\includegraphics[width=\textwidth]{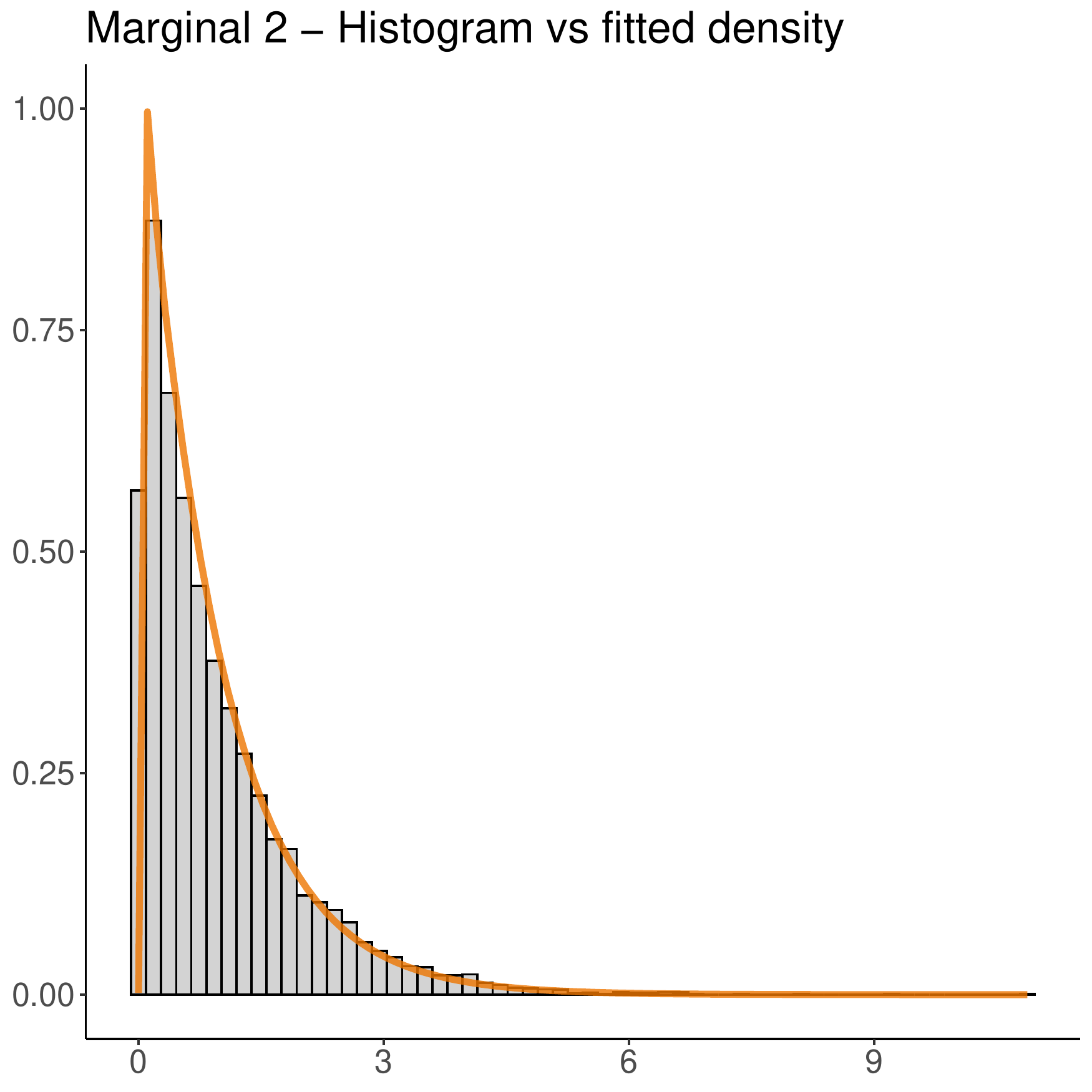}
  	\end{subfigure}
  		\begin{subfigure}{0.22\textwidth}
  		\includegraphics[width=\textwidth]{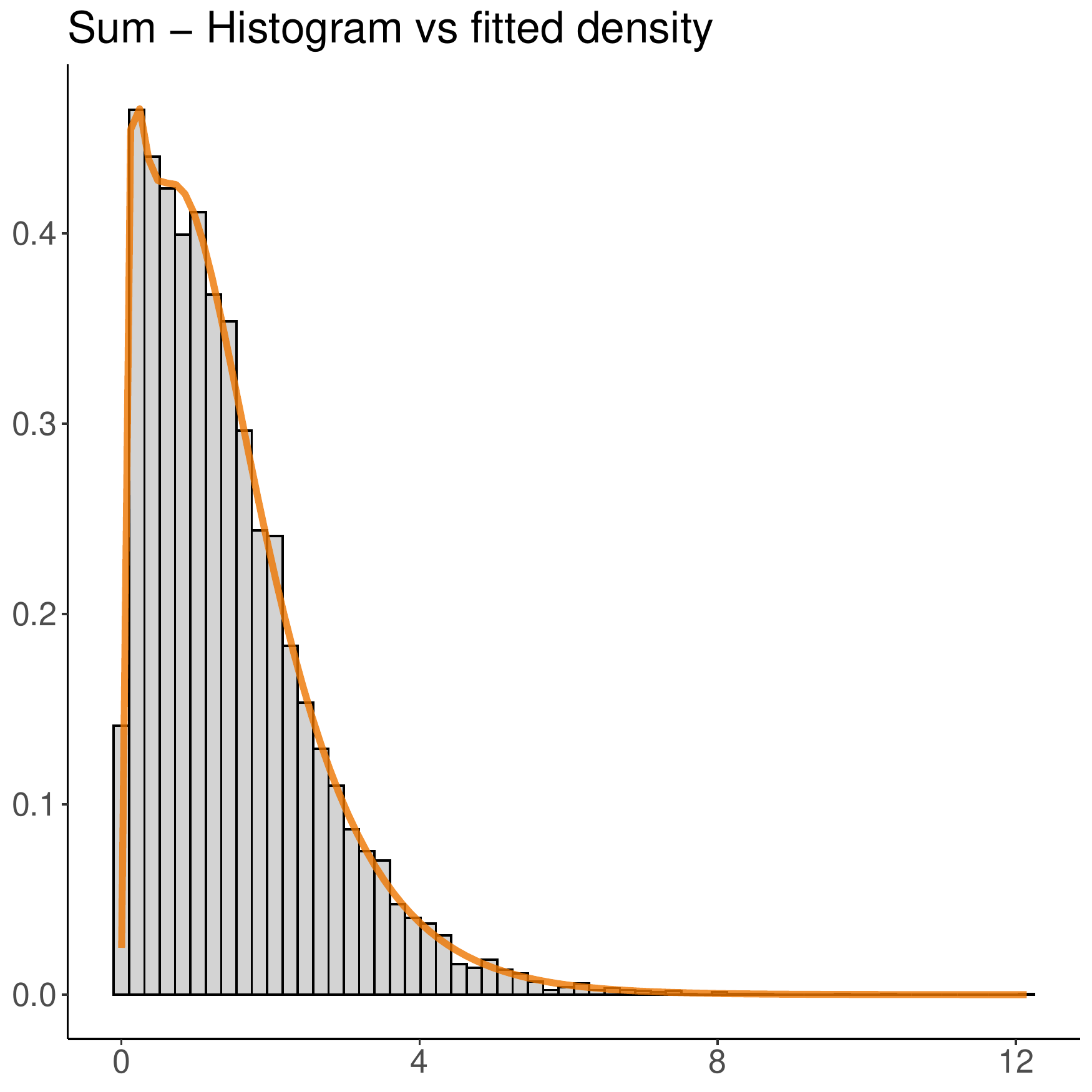}
  	\end{subfigure}
  \caption{Histograms of simulated sample versus densities of the MPH* distribution fitted using Algorithm \ref{alg:MEM}. } \label{examph:den}
\end{figure}


\begin{figure}[hbt]
	\centering
  	\begin{subfigure}{0.22\textwidth}
  		\includegraphics[width=\textwidth]{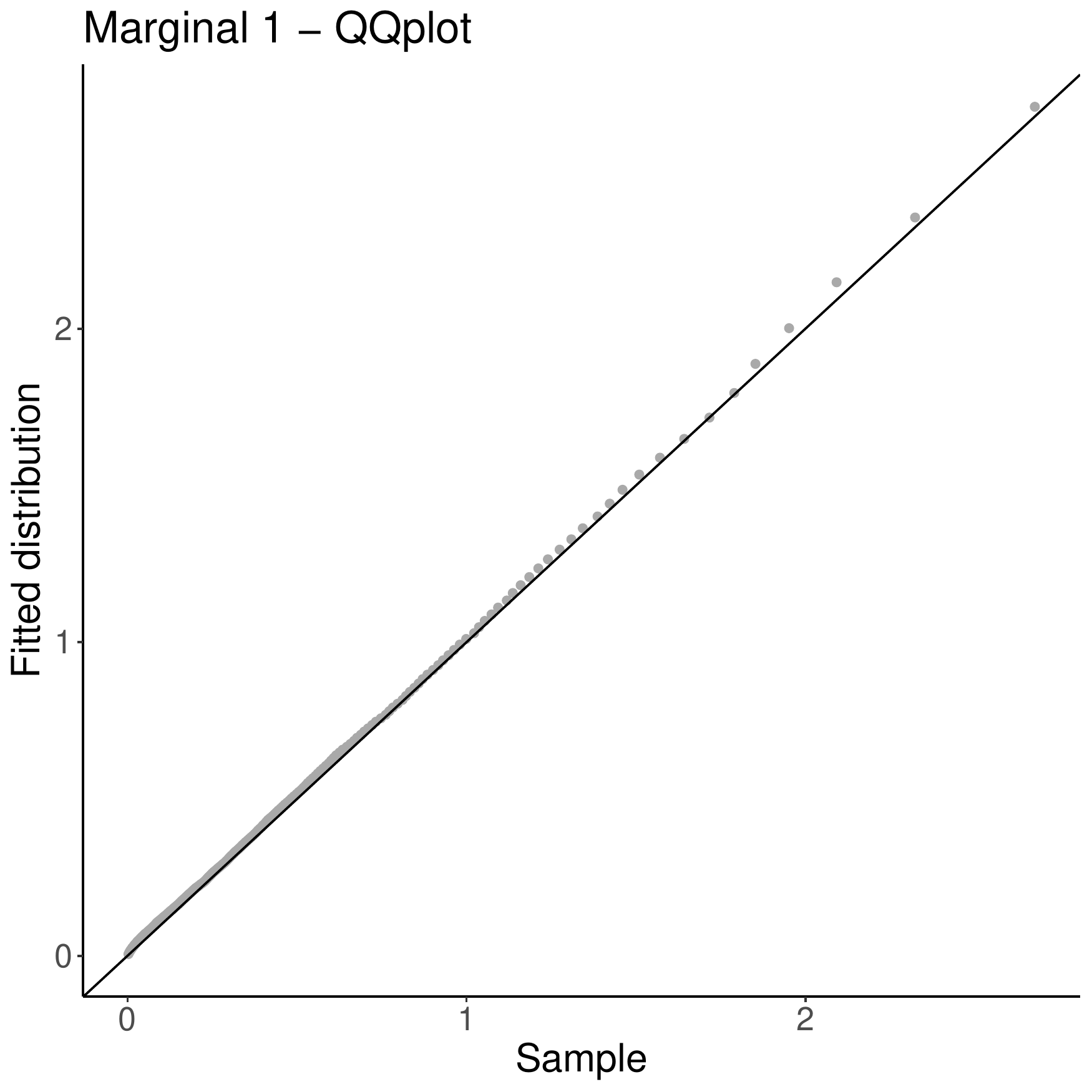}
  	\end{subfigure}
  	\begin{subfigure}{0.22\textwidth}
  		\includegraphics[width=\textwidth]{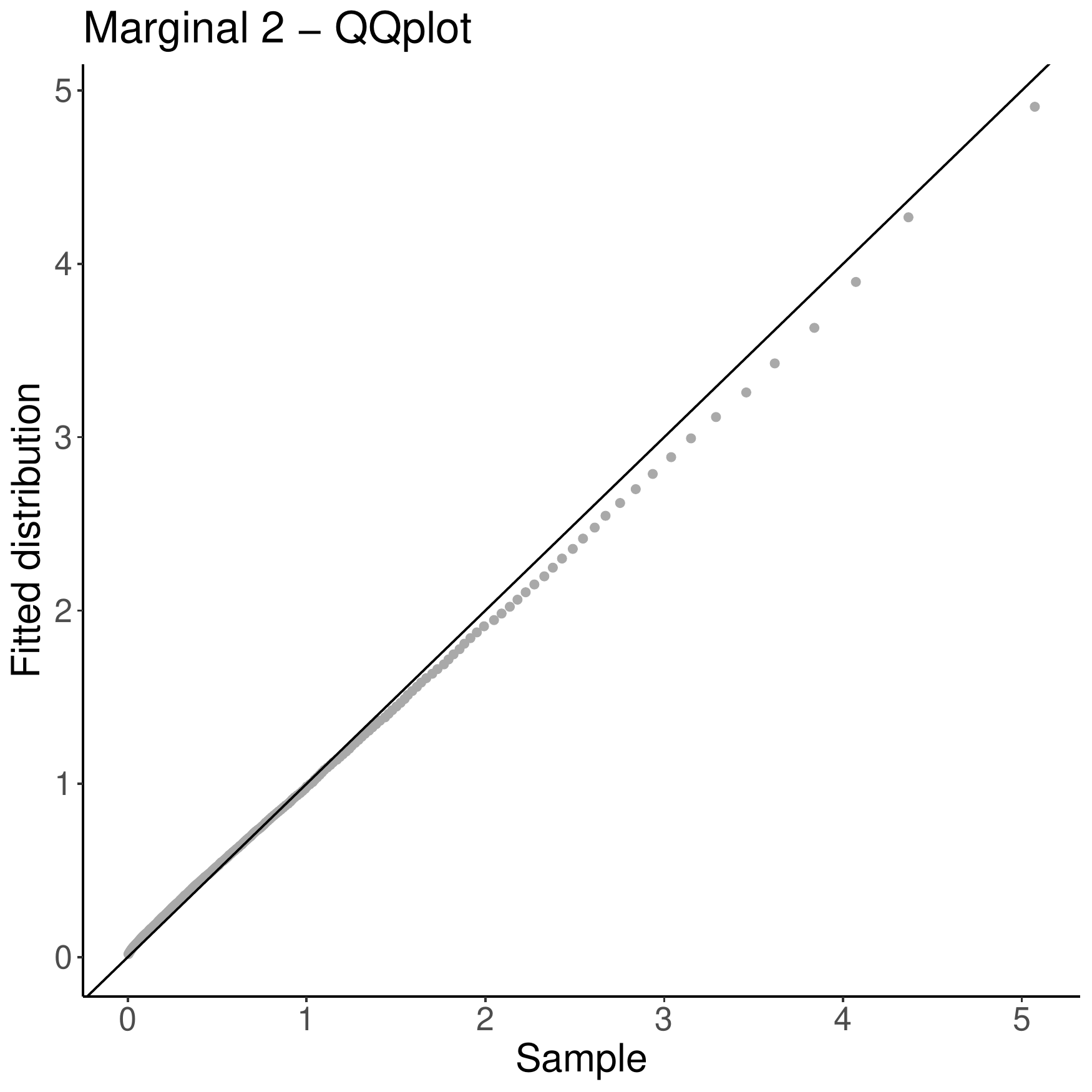}
  	\end{subfigure}
  	  \begin{subfigure}{0.22\textwidth}
  		\includegraphics[width=\textwidth]{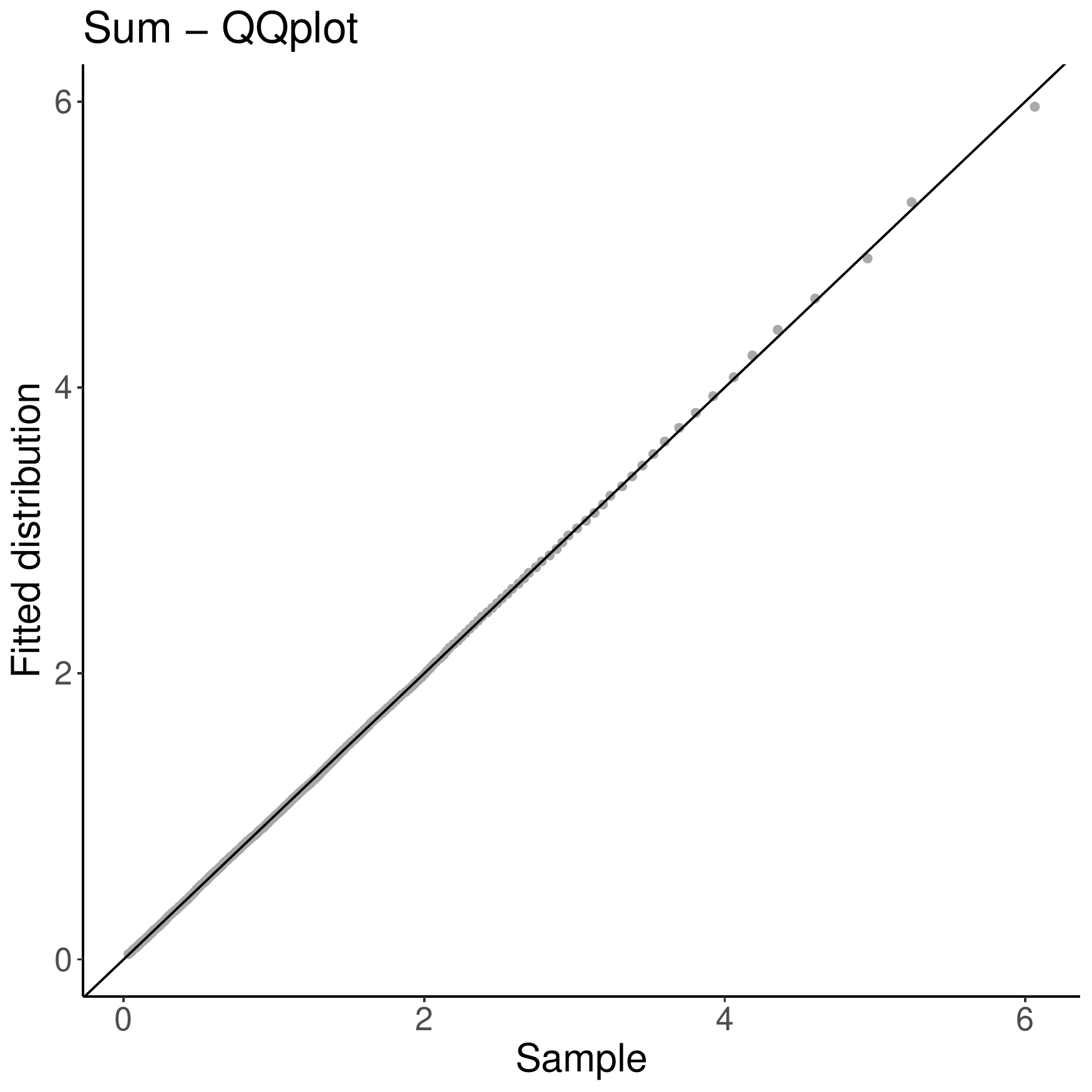}
  	\end{subfigure}
  \caption{ QQ plots of simulated sample versus  fitted MPH* distribution using  Algorithm \ref{alg:MEM}. }\label{examph:qq}
\end{figure}

\begin{figure}[hbt]
	\centering
  	\begin{subfigure}{0.22\textwidth}
  		\includegraphics[width=\textwidth]{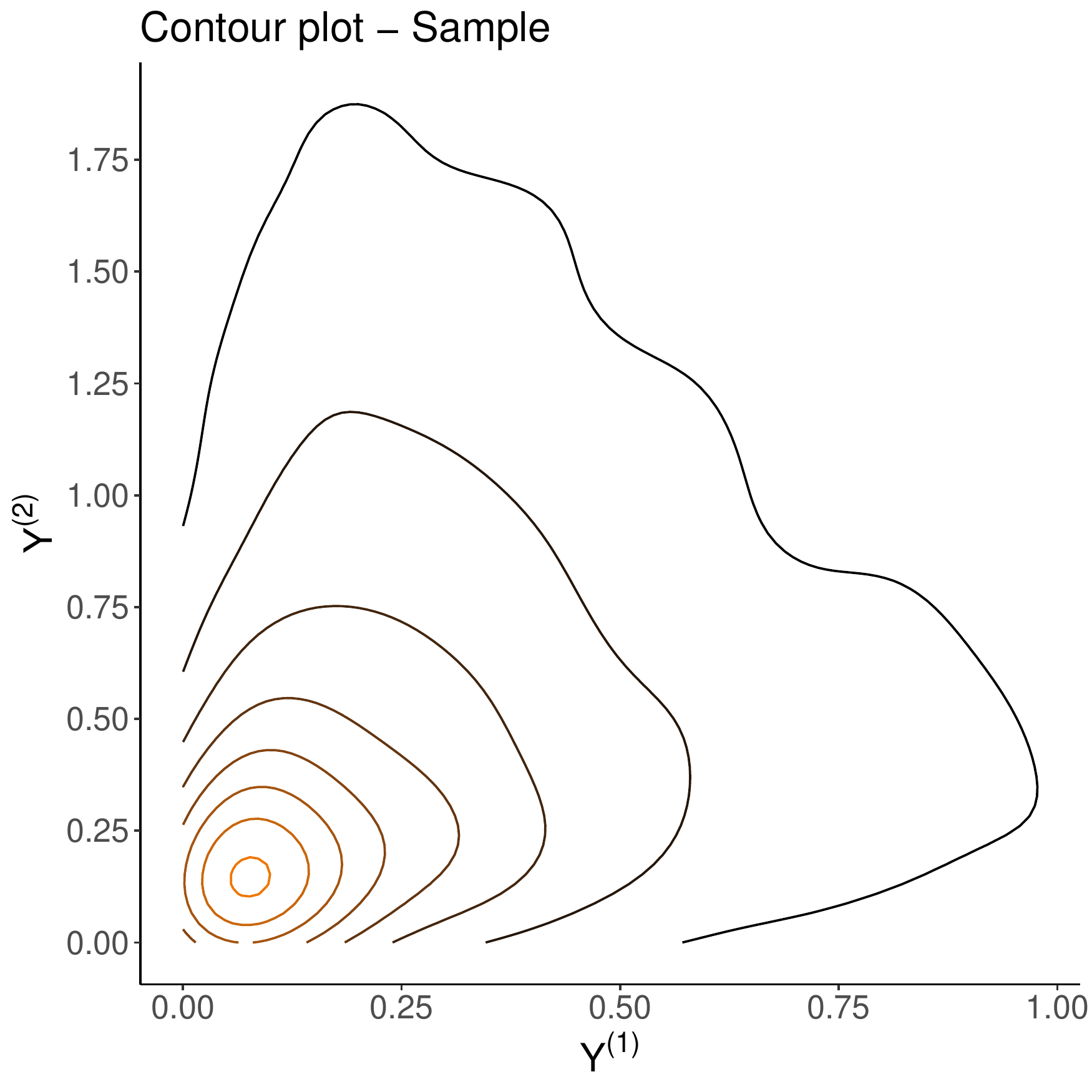}
  	\end{subfigure}
  	\begin{subfigure}{0.22\textwidth}
  		\includegraphics[width=\textwidth]{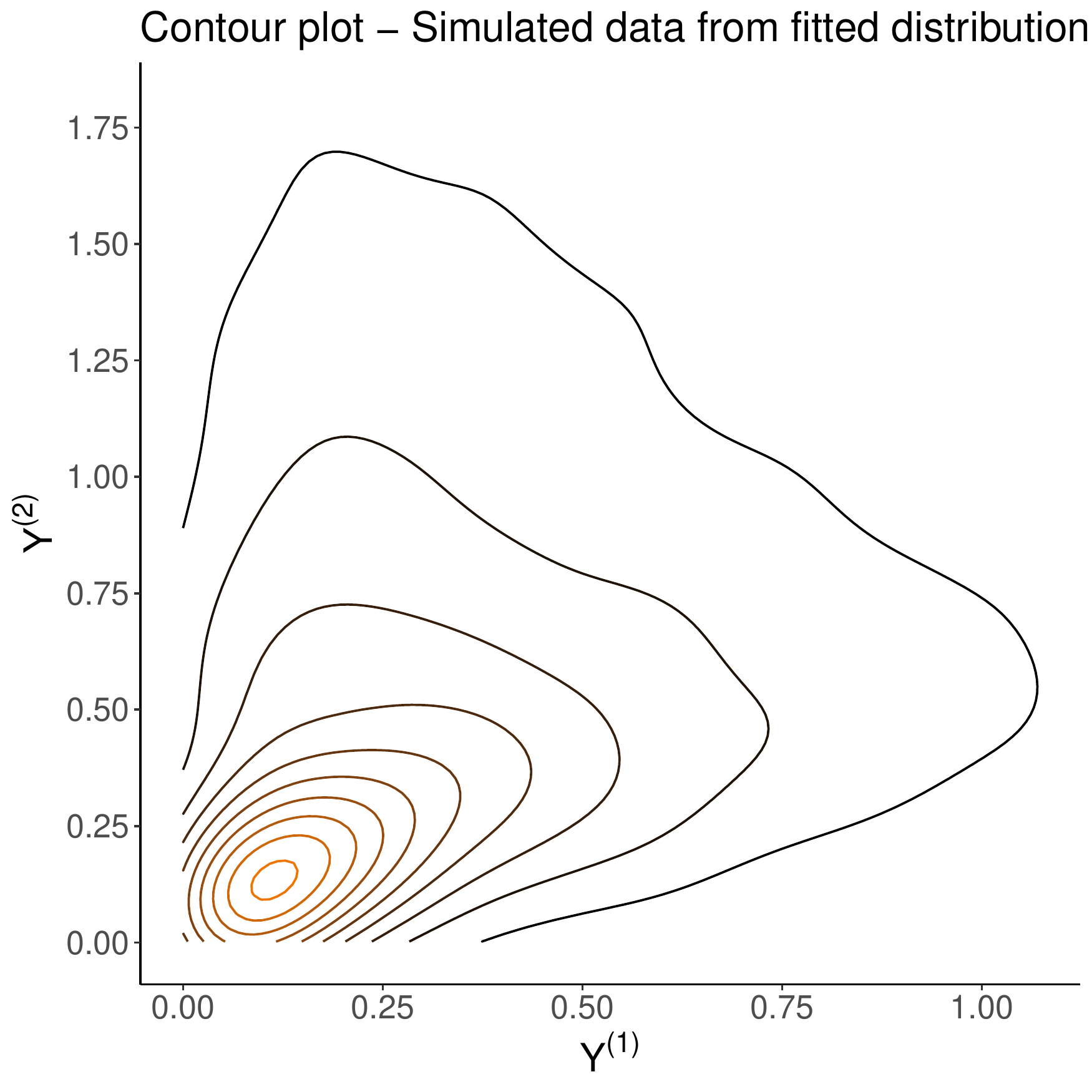}
  	\end{subfigure}
  	\begin{subfigure}{0.22\textwidth}
  		\includegraphics[width=\textwidth]{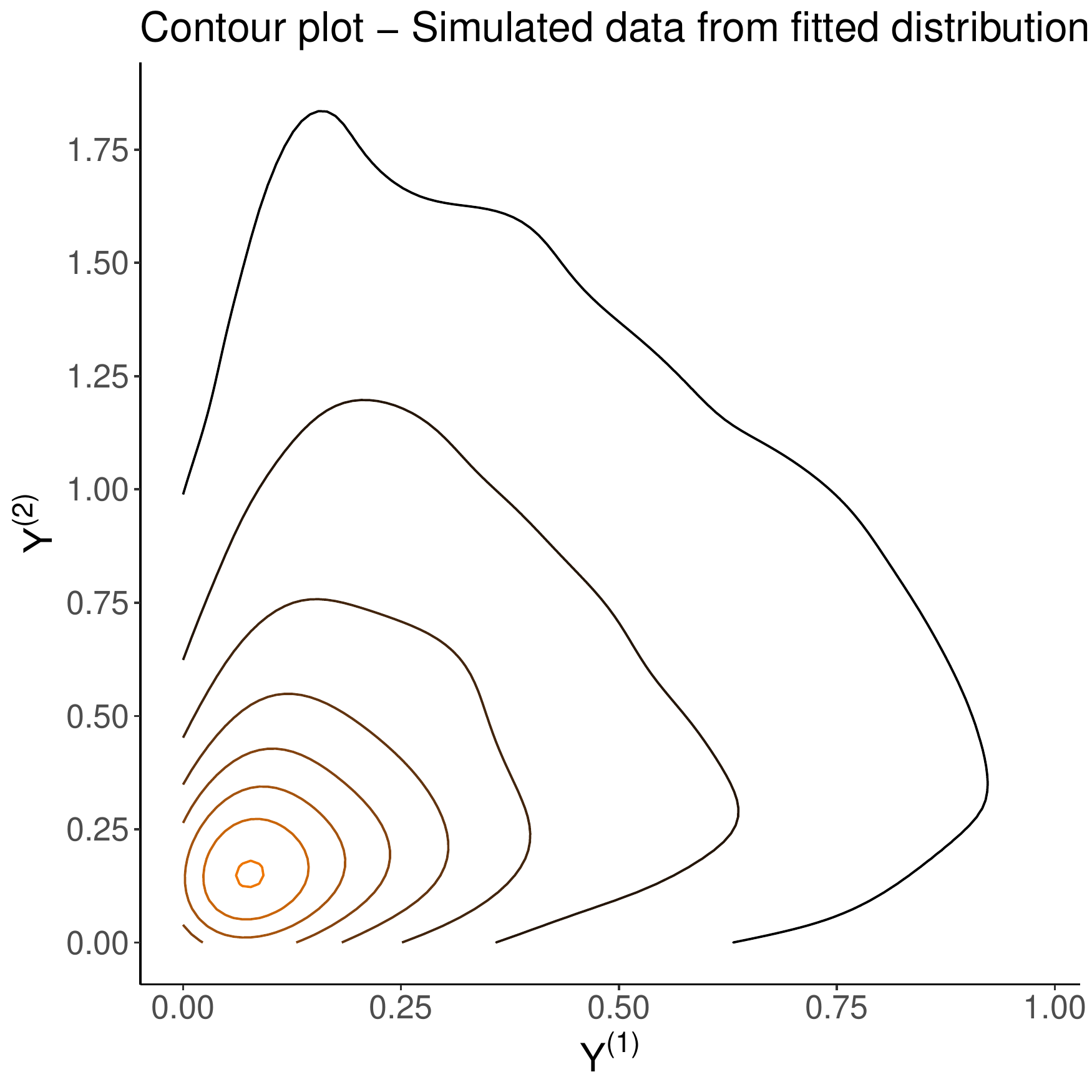}
  	\end{subfigure}
  \caption{Contour plot of sample (left), contour plot of a simulated sample from the MPH* distribution fitted with Algorithm  \ref{alg:MEM} (center) and  contour plot of of a simulated sample from the distribution fitted with Algorithm \ref{alg:bivariate} (right). }\label{examph:contour}
\end{figure}

Next we use Algorithm \ref{alg:bivariate} with 4 phases ($p_1=2$ and $p_2=2$), random initial values and $2\,500$ steps in the  EM algorithm (leading to a running time of $623$ seconds). The estimated parameters are:
\begin{gather*} 
\hat{\bfpi}=\left(
 0.1218,\, 0.8782,\, 0,\, 0 \right)\,, \\ 
\hat{\bfT}=\left( \begin{array}{cccc}
-2.2656 & 0.2832 & 1.9823 & 0.0002 \\ 
8.2508 & -10.2295 & 0.0257 & 1.9530 \\
0 & 0 & -1.0115  & 0.6419 \\
0 & 0 & 0.0992 & -4.7781 
\end{array} \right) \,, \\ 
\hat{\bfR}= \left( \begin{array}{c c}
1 & 0 \\ 
1 & 0  \\
0 & 1 \\
0 & 1
\end{array} \right)\,,
\end{gather*}
with corresponding mean 
$	\E(\bfY)= \left(
 0.5046, 0.9650 \right)^{\prime} $
	and $ \rho=0.1156$.
$\rho_{\tau}$ for this fit can be estimated via simulation, giving $\hat{\rho}_{\tau}=0.1448$. A contour plot of the fit is available in Figure~\ref{examph:contour} and we see that the algorithm recovers the original structure of the data even better. Note that again the log--likelihood of the fitted MPH* distribution ($-12\,329.82$) outperforms the log--likelihood using the original MPH* distribution ($-12\,330.96$).

Finally we would like to remark that Algorithm \ref{alg:bivariate} already starts with a more specific structure on its parameters which resembles the one of the distribution from which the data come from.  On the other hand, Algorithm  \ref{alg:MEM} does not require any prior assumption on the initial structure of its parameters.  Thus, a better fit from Algorithm \ref{alg:bivariate} is expected, since Algorithm \ref{alg:MEM} needs to find a distribution in a larger set. In line with the non-identifiability issue, one sees that one can obtain a quite reasonable fit in that larger class that captures some main features of the original distribution, whereas the more specific Algorithm \ref{alg:bivariate} finds a fit that even exhibits nicely the original correlation pattern. Yet, the flexibility of Algorithm \ref{alg:MEM} is a considerable advantage when dealing with data sets without the additional knowledge about the underlying distribution. 
\end{example}

%
%

\begin{example}[Known distribution -- Marshall-Olkin exponential] \normalfont
We now would like to illustrate that Algorithms \ref{alg:MEM} and \ref{alg:bivariate} can be modified to fit a $\mbox{MPH}^{*}$ model to a theoretically given joint  distribution $H$. The idea is along the lines of  \cite{asmussen1996fitting} and consists of considering sequences of empirical distributions with increasing sample size. We exemplify this by
	considering a bivariate Marshall--Olkin exponential distribution, whose joint survival function is of the form
	\begin{align*}
	\ov{F}\left(y^{(1)},y^{(2)} \right)=\exp \left( -\lambda_1 y^{(1)}- \lambda_2 y^{(2)} - \lambda_{12} \max ( y^{(1)}, y^{(2)} ) \right) \, .
\end{align*}
	We take $\lambda_1=1$, $\lambda_2=3$ and $\lambda_{12}=1$, then the distribution has theoretical moments
$	\E(\bfY)= \left(
0.5, 0.25 \right) $ and 
	$\rho=0.2$.
 It is easy to see that the Marshall--Olkin bivariate exponential is upper-tail-independent, i.e., $\lambda_U=0$. Moreover, we approximate $\rho_\tau$ via simulation, obtaining $\hat{\rho}_\tau=0.2012$. Then we fit a MPH* distribution using the Algorithm~\ref{alg:MEM}. 
	With 3 phases and random initial values together with $2\,500$ steps in each EM algorithm (overall running time about $120$ seconds), we obtain the parameters
	\begin{gather*} 
\hat{\bfpi}=\left(
0.8233,\, 0.1633,\, 0.0134\right)\,, \\ 
 \hat{\bfT}=\left( \begin{array}{cccc}
2.5894 & 1.6637 & 0.7601 \\ 
0.0087 & -2.0699 & 0.2102  \\
0.1032 & 0.3465 & -4.2765  
\end{array} \right) \,, \\  
 \hat{\bfR}= \left( \begin{array}{c c}
0.4412 & 0.5588 \\ 
0.9514 & 0.0486 \\
0.2348 & 0.7652
\end{array} \right)\,,
\end{gather*}
which has corresponding moments
	$\E(\bfY)= \left(
0.4942, 0.2565 \right)^{\prime}$
	and $\rho=0.3260. $
$\lambda_U$ and $\rho_\tau$ for the resulting model can be approximated by simulation to be $\hat{\lambda}_U=0.0547$ and $\hat{\rho}_\tau=0.2745$. Together with the densities (Figure~\ref{ex:mph-mo-densities}) and QQ plots (Figure~\ref{ex:mph-mo-qq}), one sees that this algorithm recovers rather well the structure of the original joint distribution.

\begin{figure}[hbt]
	\centering
  	\begin{subfigure}{0.22\textwidth}
  		\includegraphics[width=\textwidth]{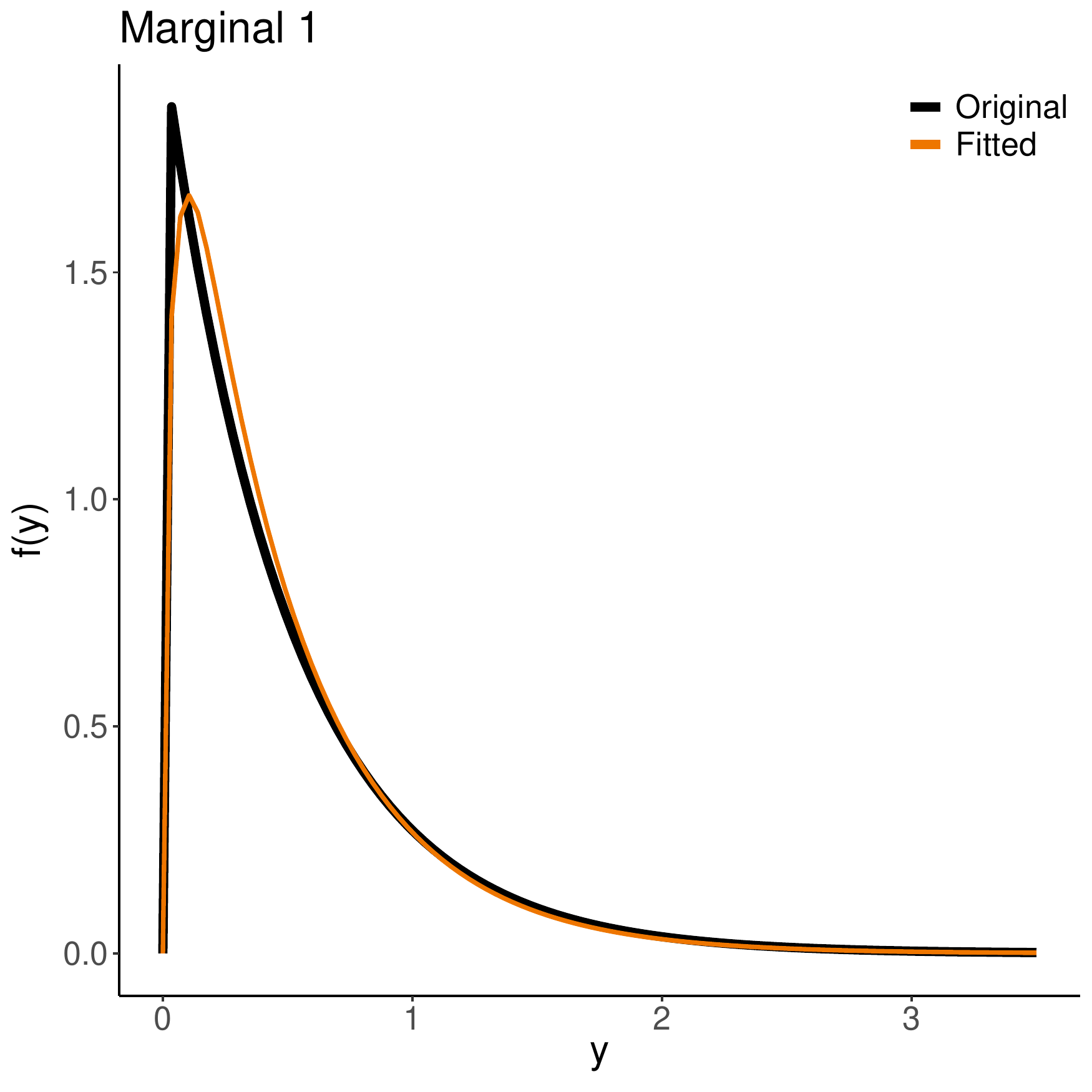}
  	\end{subfigure}
  	\begin{subfigure}{0.22\textwidth}
  		\includegraphics[width=\textwidth]{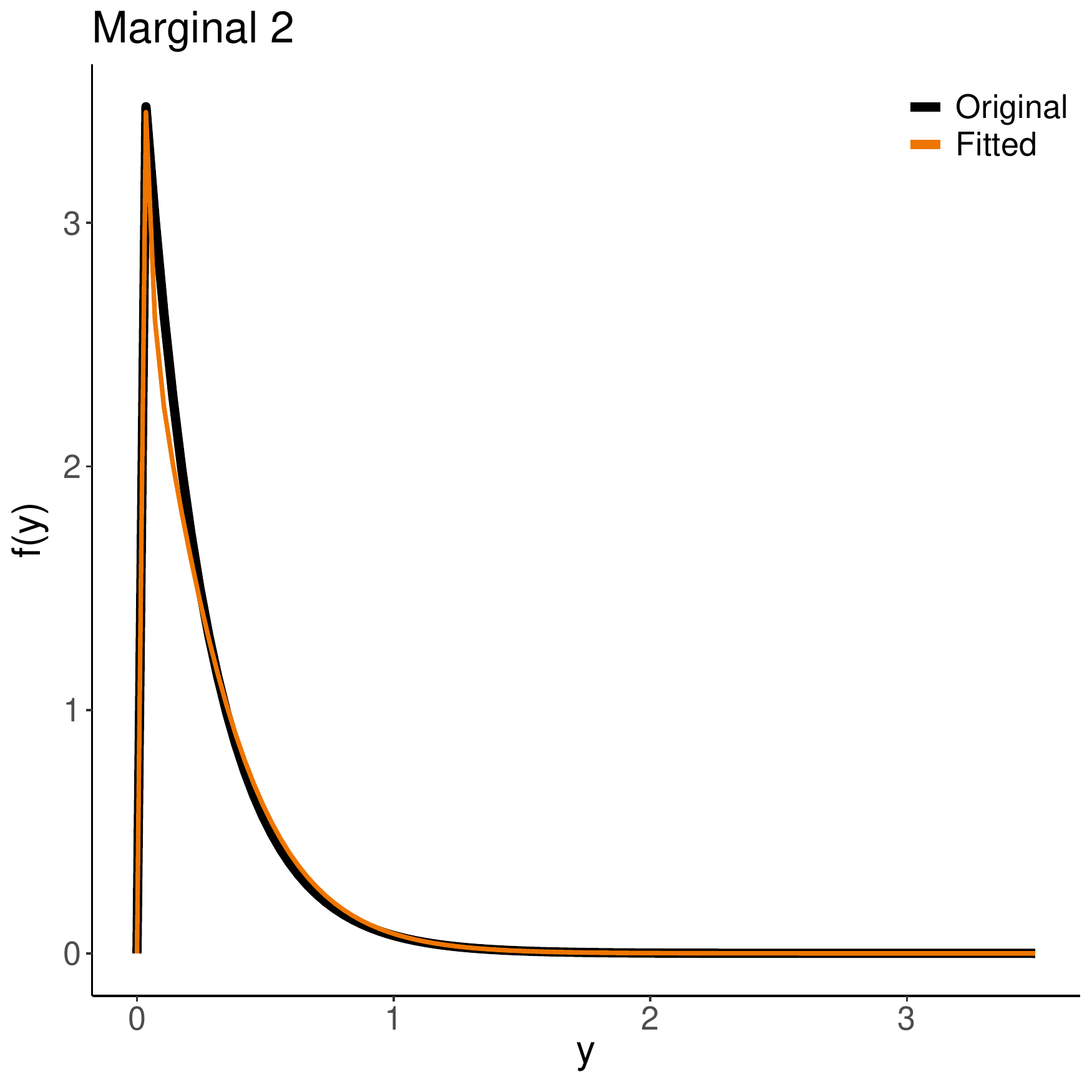}
  	\end{subfigure}
  	  	\begin{subfigure}{0.22\textwidth}
  		\includegraphics[width=\textwidth]{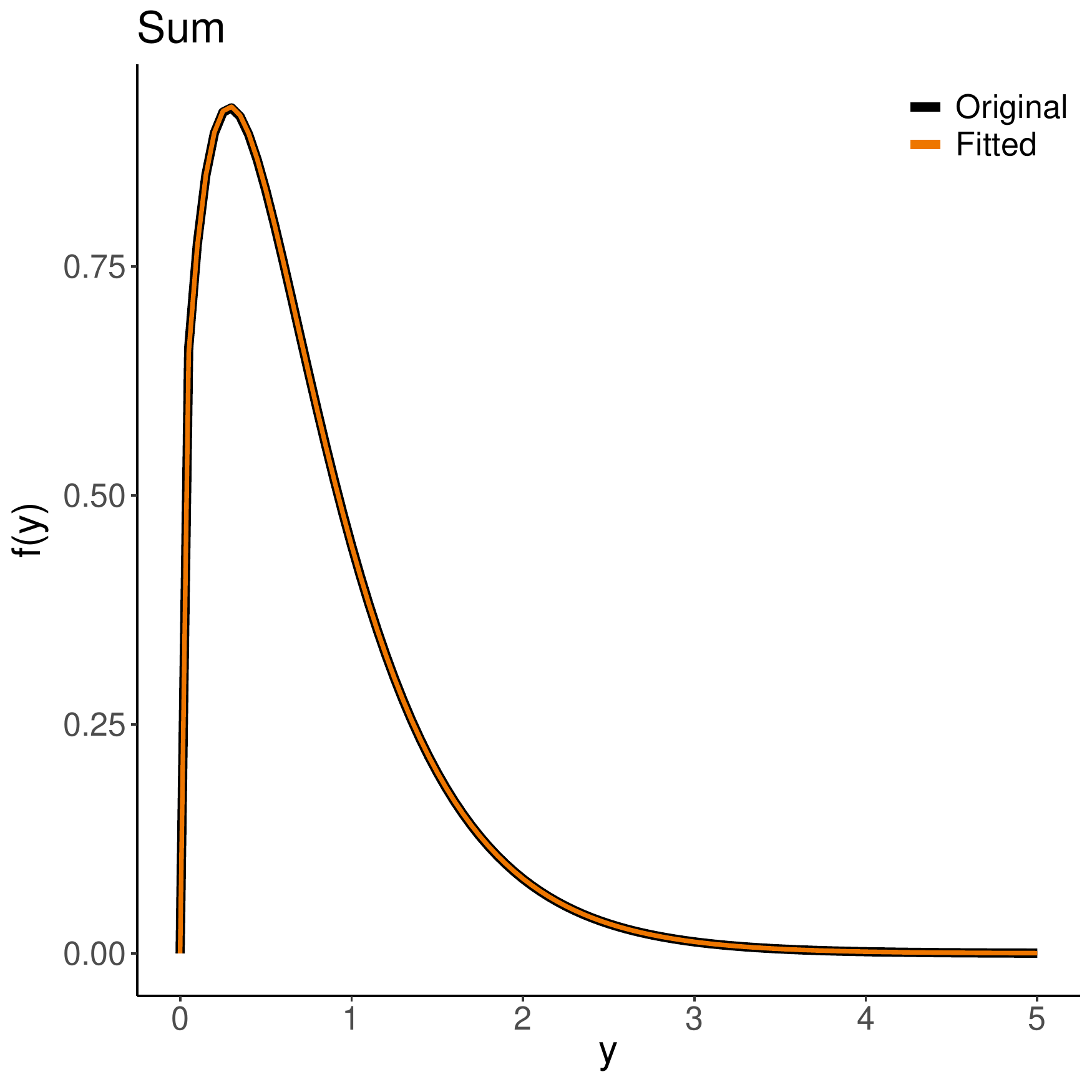}
  	\end{subfigure}
  \caption{ Densities of the original Marshall--Olkin distribution versus densities of the MPH* distribution fitted using Algorithm \ref{alg:MEM}.}\label{ex:mph-mo-densities}
\end{figure}

\begin{figure}[hbt]
	\centering
  	\begin{subfigure}{0.22\textwidth}
  		\includegraphics[width=\textwidth]{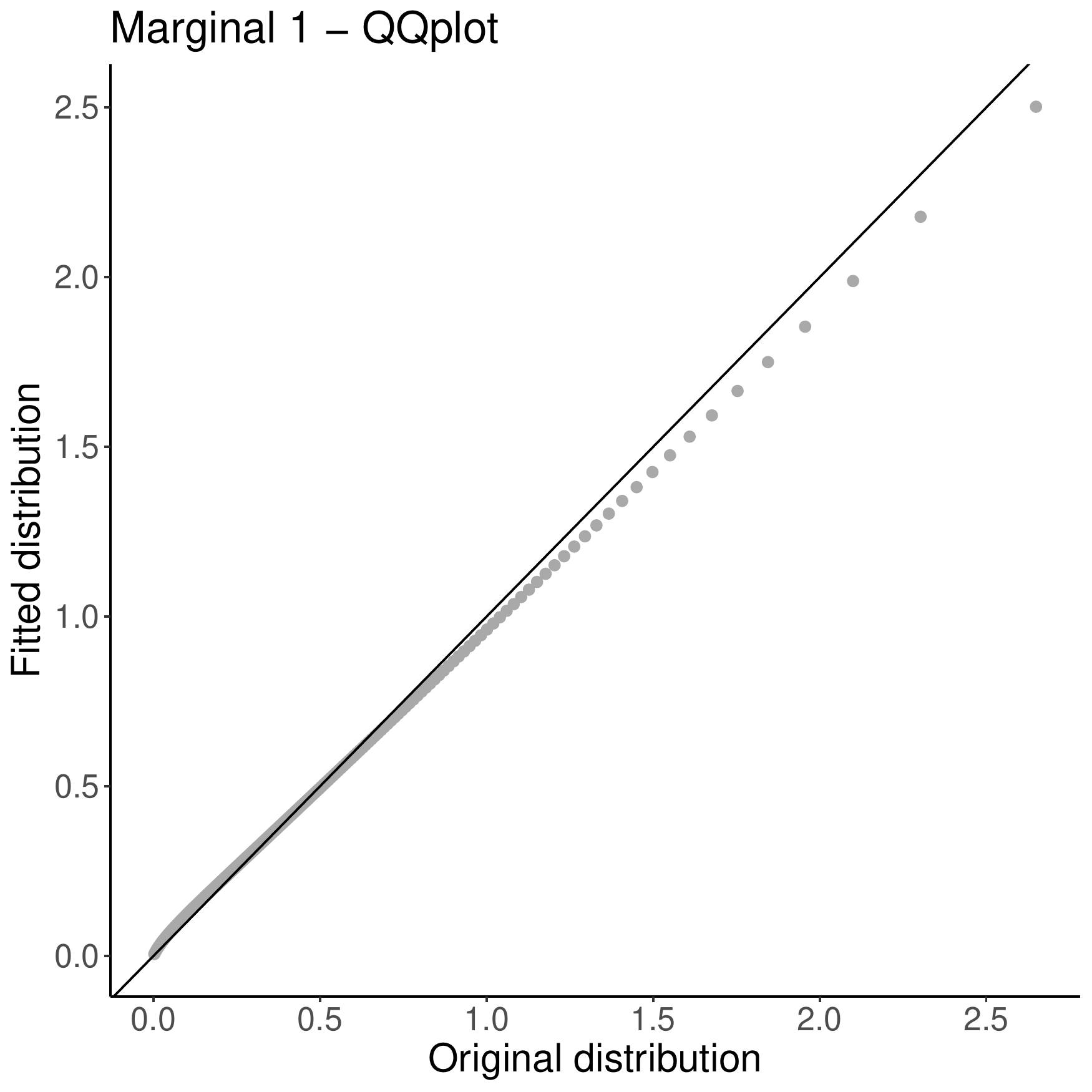}
  	\end{subfigure}
  	\begin{subfigure}{0.22\textwidth}
  		\includegraphics[width=\textwidth]{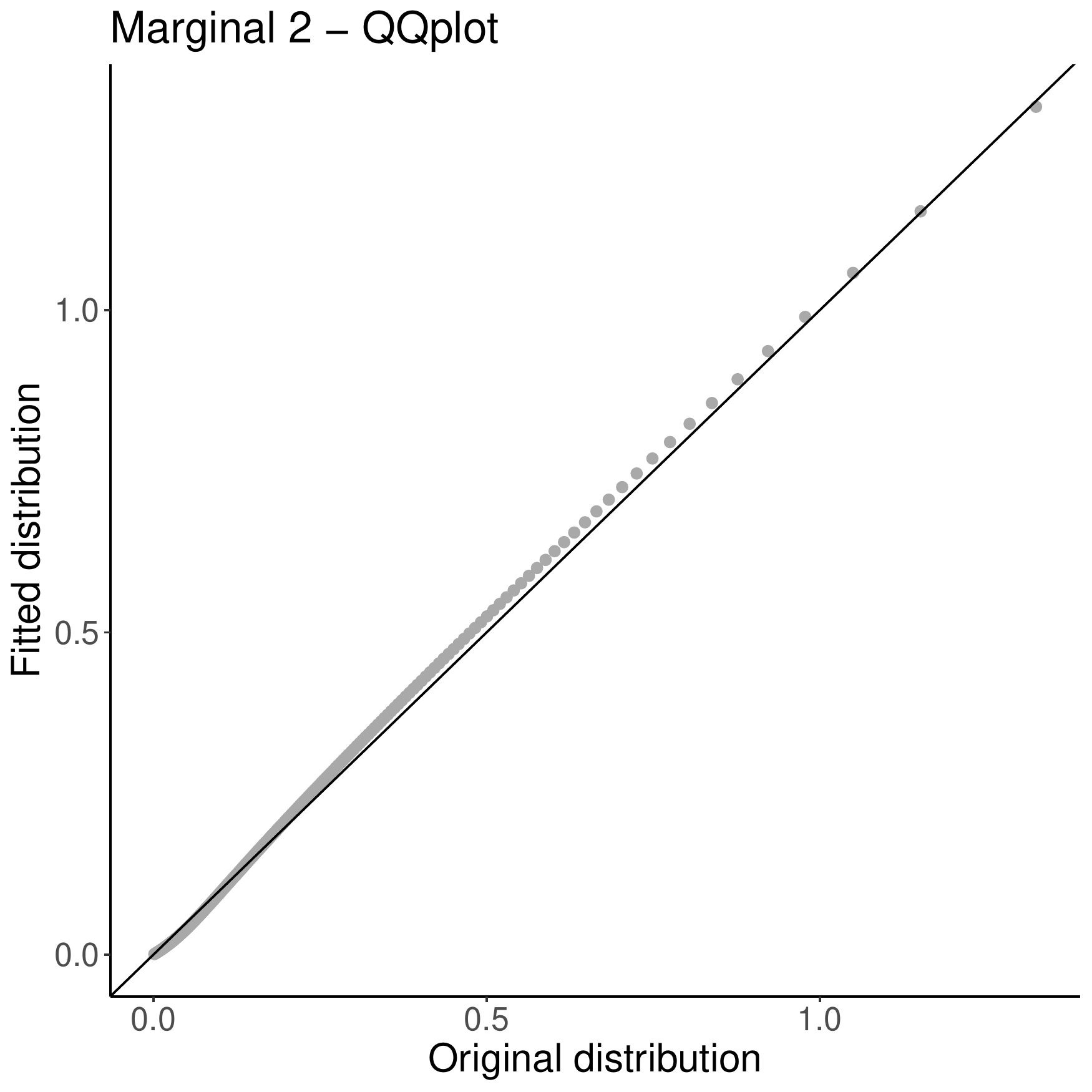}
  	\end{subfigure}
  	  	\begin{subfigure}{0.22\textwidth}
  		\includegraphics[width=\textwidth]{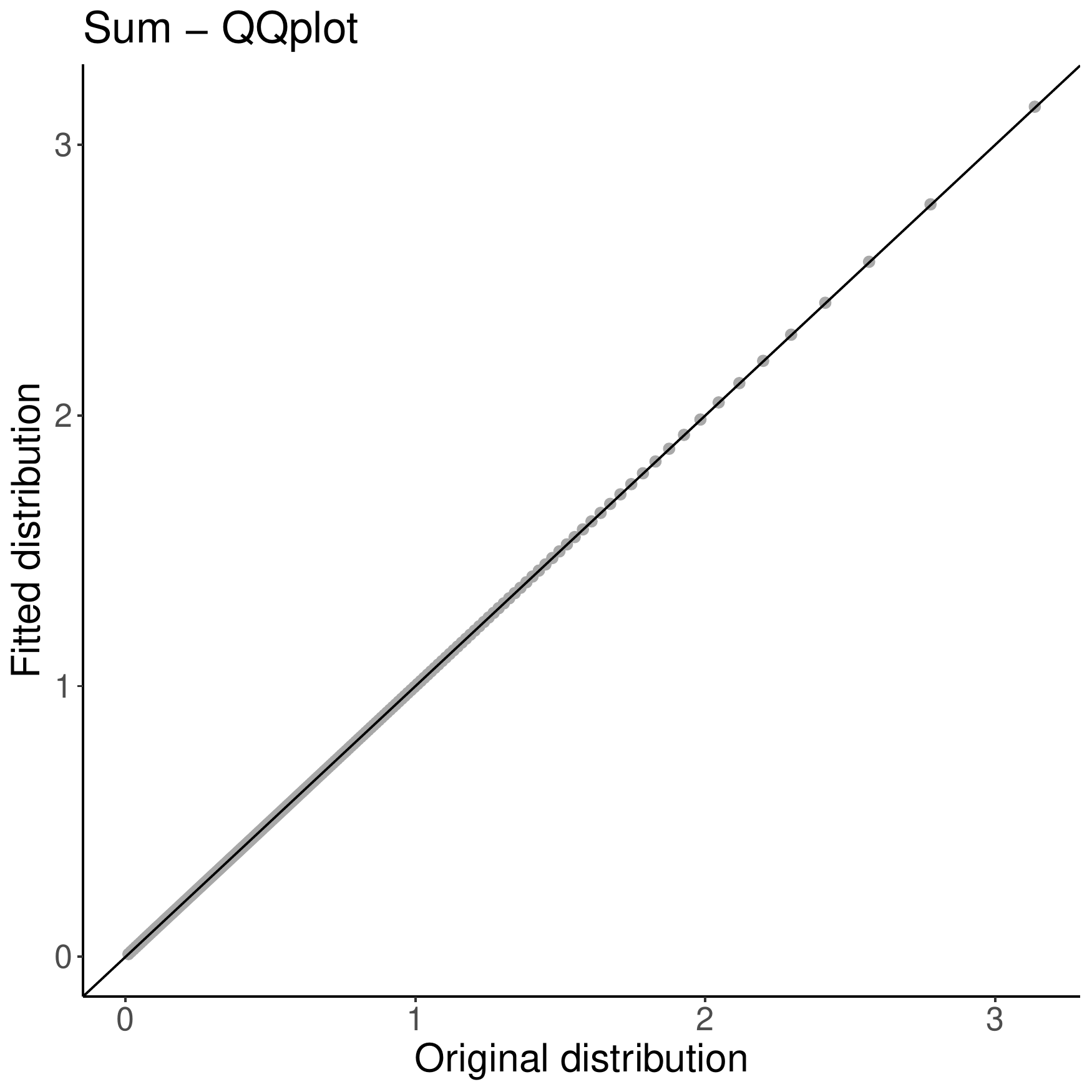}
  	\end{subfigure}
  \caption{ QQ plots of original Marshall--Olkin distribution versus fitted MPH* distribution using Algorithm~\ref{alg:MEM}. }\label{ex:mph-mo-qq}
\end{figure}

\begin{figure}[hbt]
	\centering
	\begin{subfigure}{0.22\textwidth}
  		\includegraphics[width=\textwidth]{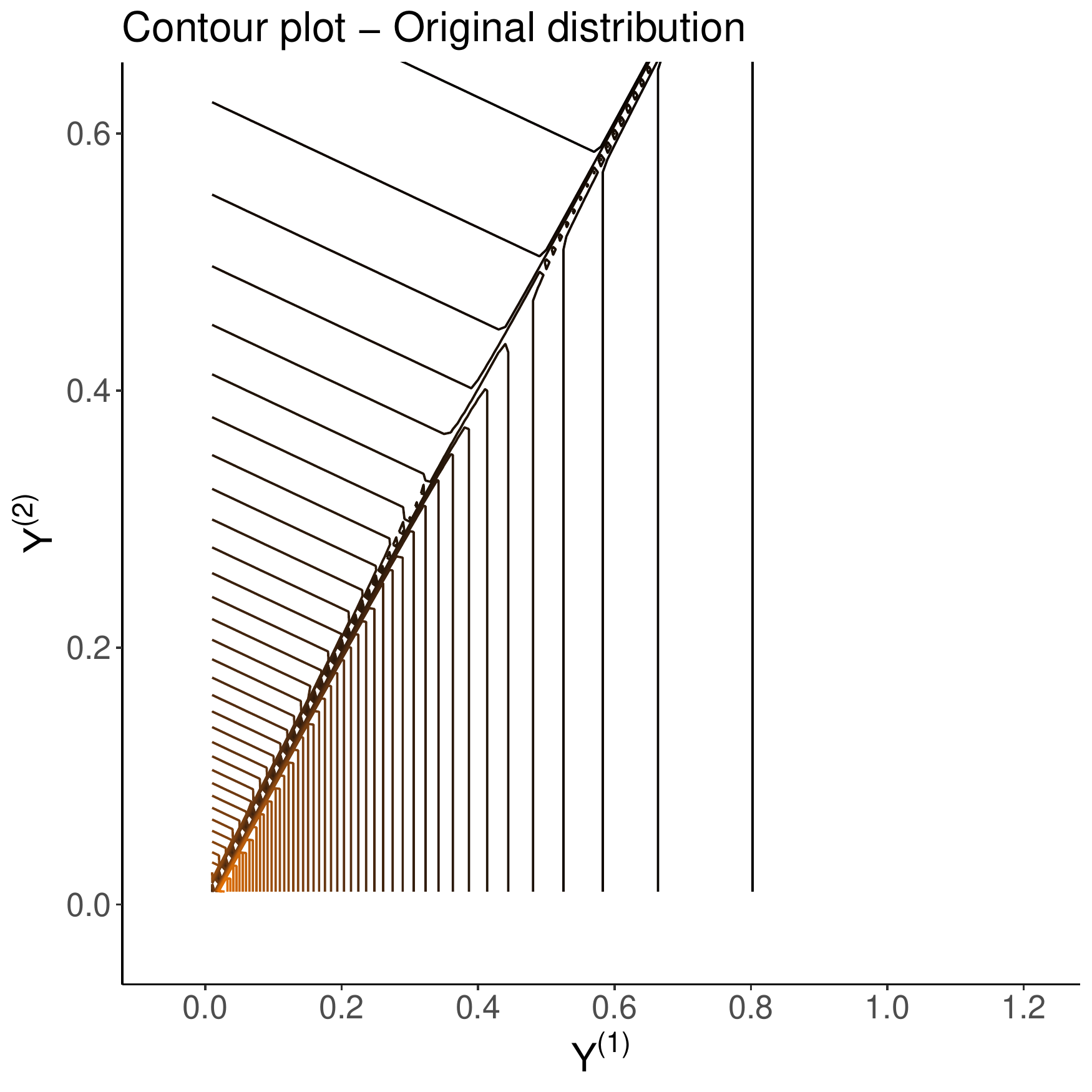}
  	\end{subfigure}
  	\begin{subfigure}{0.22\textwidth}
  		\includegraphics[width=\textwidth]{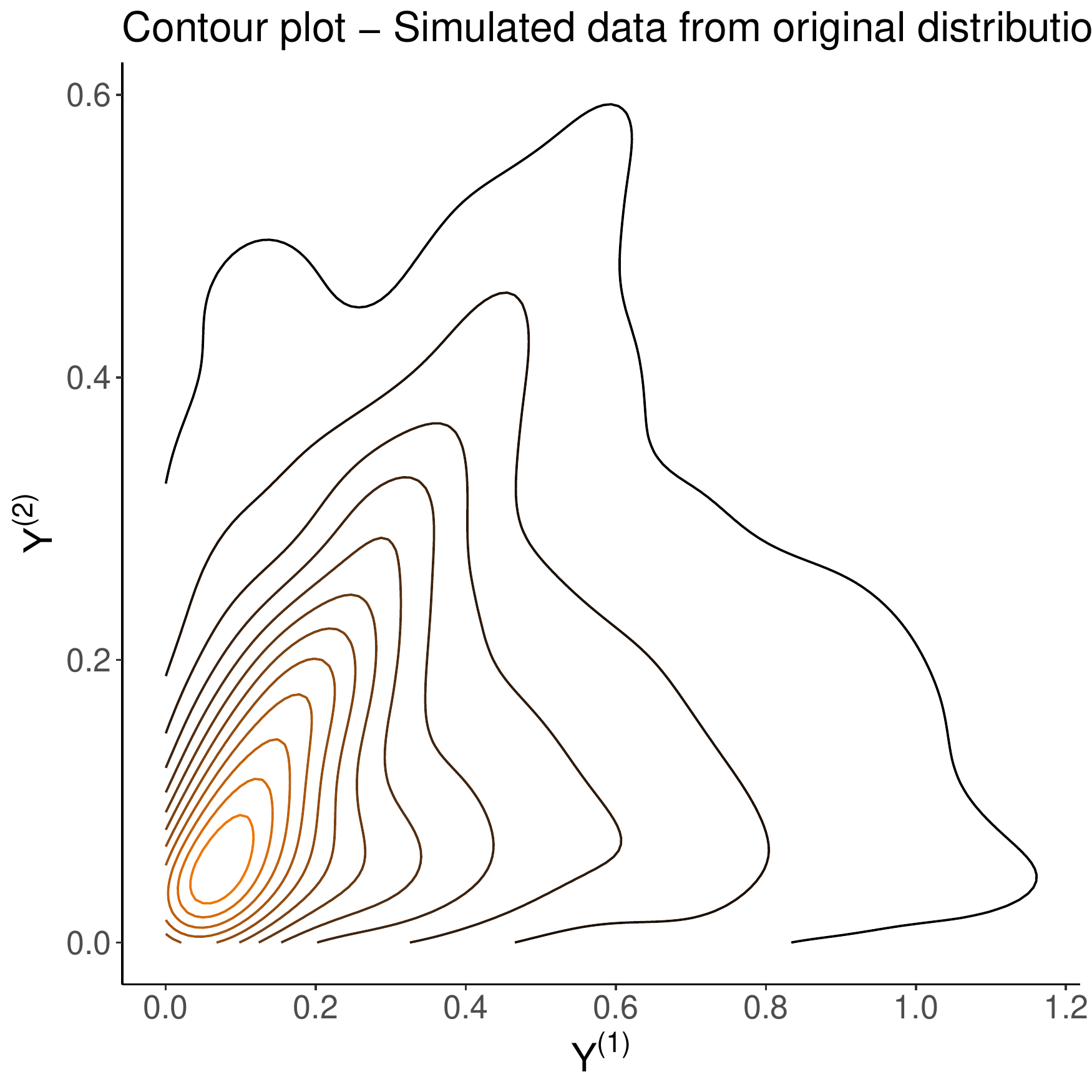}
  	\end{subfigure}
  	\begin{subfigure}{0.22\textwidth}
  		\includegraphics[width=\textwidth]{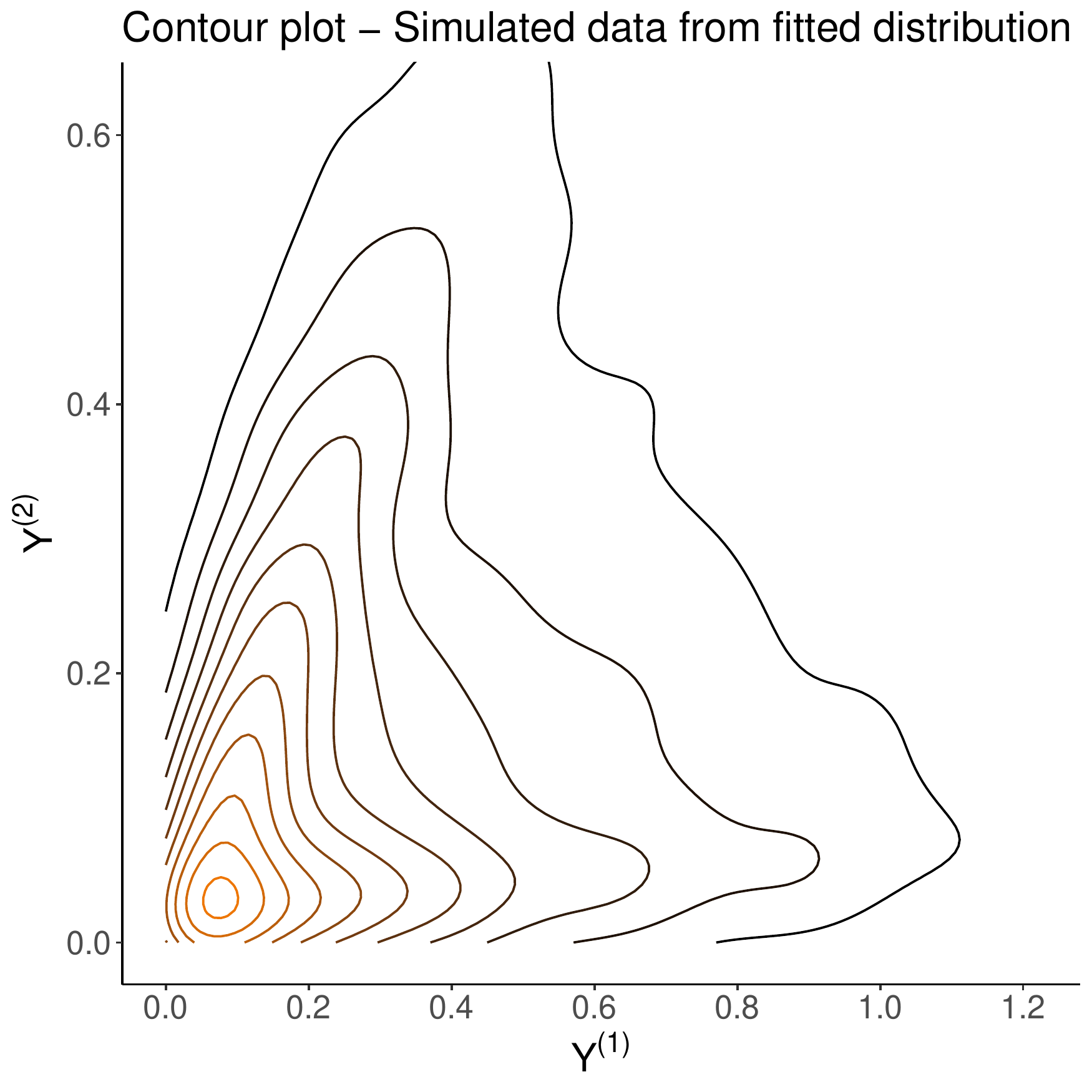}
  	\end{subfigure}
  \caption{Contour plot of the original Marshall--Olkin distribution (left) contour plot of simulated sample from Marshall--Olkin distribution (middle) and  contour plot of simulated sample from the distribution fitted with Algorithm~\ref{alg:MEM} (right).} \label{ex:mph-mo-count}
\end{figure}

\end{example}


\section{Multivariate inhomogeneous phase-type distributions}\label{sec:miph}
There are various possibilities for extending inhomogeneous phase-type distributions to more than one dimension. In the following we suggest one particular approach and provide an algorithm for the parameter estimation. 

\subsection{Definition and Properties}
Let $\bfY \sim 	\mbox{MPH}^{*}\left( \bfpi, \bfT, \bfR \right)$ and define $\bfX := (g_1(Y^{(1)}), \dots, g_d(Y^{(d)}))^{\prime} $, where $g_j : \R_{+} \to \R_{+}$ are increasing and differentiable functions {\markcol with range $\R_{+}$} for $j=1,\dots,d$, then we say that $\bfX$  has an \textit{inhomogeneous MPH*} distribution.

Several of its properties follow directly from the definition:  
\begin{enumerate}
	\item The marginals $X^{(j)} = g_j(Y^{(j)})$ are IPH distributed, since each $Y^{(j)}$ is phase--type distributed, $j=1,\dots,d$.
	\item {\markcol If $g_{j}$ is strictly} increasing for all $j=1,\dots,d$, the copula of $\bfX$ is the same as the copula of $\bfY$ (see e.g.\ \cite[Prop.7.7]{mcneil2015quantitative}).
	\item For fixed $g_j (\cdot)$, $j=1,\dots,d$, this new class is dense in $\R_{+}^{d}$ (by the denseness of the MPH* class). 
\end{enumerate}
In the sequel we will provide an algorithm for parameter estimation for which an explicit expressions of the bivariate density is needed. We therefore restrict it to the bivariate case.

\subsection{Parameter estimation in the bivariate case}
In the bivariate case, for any  $\bfY \sim \mbox{MPH}^{*}( \bfpi, \bfT, \bfR) $ with parameters \eqref{par:biv}, it is easy to see that the density of $\bfX = (g_1(Y^{(1)}), g_2(Y^{(2)}))^{\prime}$ is given by
\begin{align} \label{eq:denbiviph}
f_{\bfX} \left( x^{(1)}, x^{(2)} \right)=\bfalp \ex^{\bfT_{11} g_{1}^{-1}(x^{(1)}) } \bfT_{12}\ex^{\bfT_{22} g_{2}^{-1}(x^{(2)}) }(-\bfT_{22})\bfe \, \frac{1}{ g_{1}^{\prime}(g_{1}^{-1}(x^{(1)}))g_{2}^{\prime}(g_{2}^{-1}(x^{(2)})) } \,.
\end{align}

If we assume that $ g_{j} ( \,\cdot\, ; \vect{\beta}_{j} )$ is a parametric non--negative function depending on the vector $\vect{\beta}_{j}$, $j=1,2$, and let $\vect{\beta}  = (\vect{\beta}_{1} , \vect{\beta}_{2} )$. Then, 
we can formulate an algorithm analogous to  Algorithm~\ref{alg:transPH}:

\begin{algorithm} [EM algorithm for bivariate inhomogeneous MPH* distributions]\label{alg:BivIPH} \

	0. Initialize with some ``arbitrary'' $( \bfalp,\bfT , \vect{\beta} )$.

1. Transform the data into $y_{i}^{(j)}:= g_{j}^{-1}( x_{i}^{(j)}; \vect{\beta}_{j})$, $i=1,\dots,N$, $j=1,2$, and apply the E-- and M--steps of Algorithm~\ref{alg:bivariate} by which we obtain the estimators $( \hat{\bfalp},\hat{\bfT})$.

2.  Compute  
\begin{align*}
	\hat{\vect{\beta}}  & = \argmax_{\vect{\beta}} \sum_{i=1}^{N} \log (f_{\bfX}(x_{i}^{(1)}, x_{i}^{(2)}; \hat{\bfalp}, \hat{\bfT},\vect{\beta} ))
\end{align*}

3. Assign $(\bfalp, \bfT ,\vect{\beta}) = (\hat{\bfalp}, \hat{\bfT},\hat{\vect{\beta}})$ and GOTO 1. 
\end{algorithm}

We now consider particular multivariate distributions obtained through such a transformation of an MPH* random vector.

\subsection{Multivariate matrix--Pareto models}

Let $\bfX = (g_1(Y^{(1)}), \dots, g_d(Y^{(d)}))^{\prime} $,  where\\ $\bfY \sim 	\mbox{MPH}^{*}\left( \bfpi, \bfT, \bfR \right)$ and 
\begin{equation}\label{gj}
	g_j(y)=\beta_{j}(e^{y}-1),\quad \beta_j>0,\;j=1,\dots,d.
\end{equation} 
Then we say that $\bfX$  follows a \textit{multivariate matrix--Pareto distribution}. Some special properties of this class of distributions are: 
\begin{enumerate}
	\item Marginal distributions are matrix--Pareto distributed.
	\item Moments and cross--moments of $\bfX$ can be obtained from the moment generating function of $\bfY$ (see \cite[Theorem 8.1.2]{Bladt2017}), provided that they exist.
	\item Products of the type $ \prod_{i=1}^{M}\left(\frac{X^{(j_{i})}}{\beta_{j_{i}}} +1 \right)^{a_i}$ are matrix--Pareto distributed, $a_i>0$, $j_{i} \in \{1, \dots,d \}$, $i=1,\dots,M$, $M\leq N$, since linear combinations of $\bfY$ are PH distributed. 
\end{enumerate} 
In the bivariate case, \eqref{gj} and \eqref{eq:denbiviph}, lead to 
\begin{align*}
	f_{\bfX}( x^{(1)},x^{(2)})=\bfalp \left(\frac{x^{(1)}}{\beta_1} +1\right)^{\bfT_{11} -\mat{I} } \bfT_{12} \left(\frac{x^{(2)}}{\beta_2} +1\right)^{\bfT_{22}-\mat{I}}(-\bfT_{22} ) \, \bfe \,\frac{1}{\beta_1 \beta_2} \,. 
\end{align*}
and
\begin{align*} 
	\bar{F}_{\bfX}( x^{(1)}, x^{(2)})=\bfalp \left( -\bfT_{11} \right)^{-1} \left(\frac{x^{(1)}}{\beta_1} +1 \right)^{\bfT_{11}} \bfT_{12} \left(\frac{x^{(2)}}{\beta_2} +1\right)^{\bfT_{22}}\bfe \,
\end{align*}
 Moreover, in this bivariate case, linear combinations of $X^{(i)}$ are regularly varying, and the respective index is the real part of the eigenvalue with largest real part of the sub--intensity matrices of the marginals, which follows from \cite[Lem. 2.1]{davis1996limit} and asymptotic independence of $\bfY$.
	The general case is not clear since the condition of asymptotic independence does not hold.


\begin{remark} \rm 
	The Marshall--Olkin Pareto distribution (see \cite{hanagal1996multivariatePareto}) is a particular case of this class of distributions.  
\end{remark}

\subsubsection{Parameter estimation} \label{subsub:MPEst}
As in the univariate case, if we assume the simpler transformation $g_j(y)=e^{y}-1$, $j=1,\dots,d$, then
we can use fitting methods of the MPH* class by taking the logarithm of the marginal observations. I.e., we can apply Algorithms~\ref{alg:MEM} and \ref{alg:bivariate} to the transformed data $y_{i}^{(j)}:= \log( x_{i}^{(j)}+1)$, $i=1,\dots,N$, $j=1,\dots,d$, to estimate the parameters $\left( \bfpi, \bfT, \bfR \right)$. We exemplify the use of this method in two examples. 

\begin{example}\normalfont (Mardia type I)
	We generated an i.i.d.\ sample of size $10\,000$ from a (translated) Mardia type I Pareto distribution (see \cite{mardia1962multivariate}) with parameters $\sigma_1=\sigma_2=1$ and $\alpha=2$. This distribution has theoretical numerical values $\E(\bfX)=(1, 1)^{\prime}$, $\lambda_U=0$ and $\rho_{\tau}=0.2$.	
The simulated sample has numerical values $\hat{\E}(\bfX)= \left( 0.9812  , 0.9712 \right)^{\prime} $ and $\hat{\rho}_\tau=0.2049$. 

We fitted a bivariate matrix--Pareto distribution using Algorithm~\ref{alg:bivariate} with $p_1=p_2=2$ (i.e., $p=4$) and $2\,000$ steps on the transformed data (with a running time of $530$ seconds), getting the following parameters: 
	
\begin{gather*} 
 \hat{\bfalp}=\left(
 0.1666 ,\, 0.8334 ,\, 0 ,\, 0 \right)\,, \\ 
\hat{\bfT}=\left( \begin{array}{cccc}
-2.0022 & 0.0238 & 1.9784 & 0 \\ 
9.2506 & -11.2967 & 0.1204  & 1.9256 \\
0 & 0 & -2.0175  & 1.8247 \\
0 & 0 & 0.0834 & -12.8829
\end{array} \right)\,, \\ 
\hat{\bfR}= \left( \begin{array}{c c}
1 & 0 \\ 
1 & 0  \\
0 & 1 \\ 
0 & 1 
\end{array} \right)\,.
\end{gather*}
The tails of the marginals of the fitted distribution are determined by the real part of the eigenvalues with largest real part of the sub--intensity matrices of the marginal  distributions, which are $\lambda_{1}^{(\max)}=-1.9785$ and  $\lambda_{2}^{(\max)}=-2.0035$. These resemble well the ones of the original distribution. 
The fitted distribution has first moment 
$	\E(\bfX)= \left(
 1.0164, 0.9963 \right)^{\prime} $. Moreover, we estimated $\rho_{\tau}$ via simulation, obtaining $\hat{\rho}_{\tau}=0.1582$. The QQ plots are available in Figure~\ref{fig:MardiaQQ} and contour plots are depicted in Figure~\ref{fig:Marcontour}, from where it becomes clear that the algorithm recovers the structure of the data well. Again, the log--likelihood of the fitted bivariate matrix--Pareto ($-15\,550.52$) exceeds the log--likelihood using the original Mardia distribution ($-16\,146.65$).

	
\begin{figure}[hbt]
	\centering
  	\begin{subfigure}{0.22\textwidth}
  		\includegraphics[width=\textwidth]{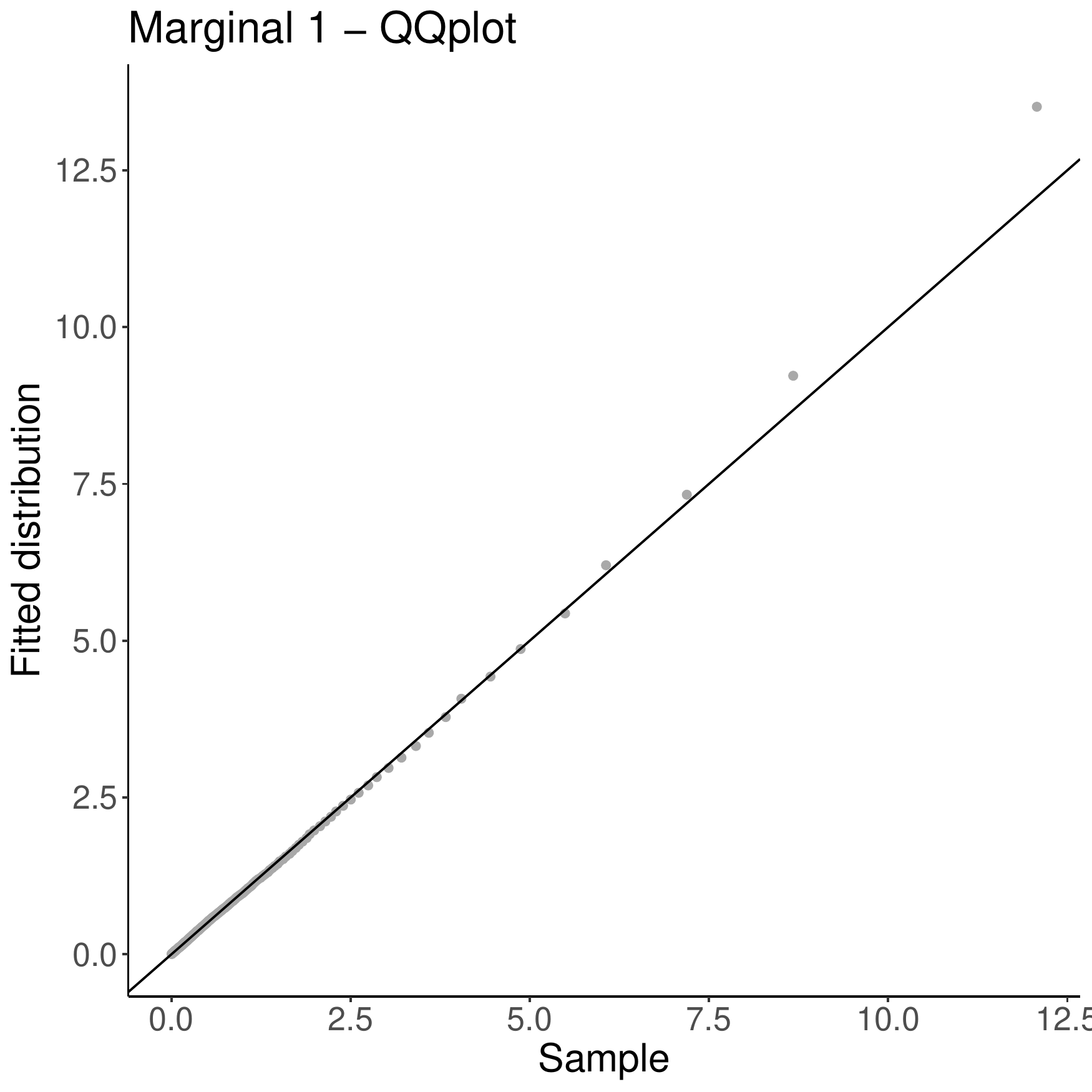}
  	\end{subfigure}
  	\begin{subfigure}{0.22\textwidth}
  		\includegraphics[width=\textwidth]{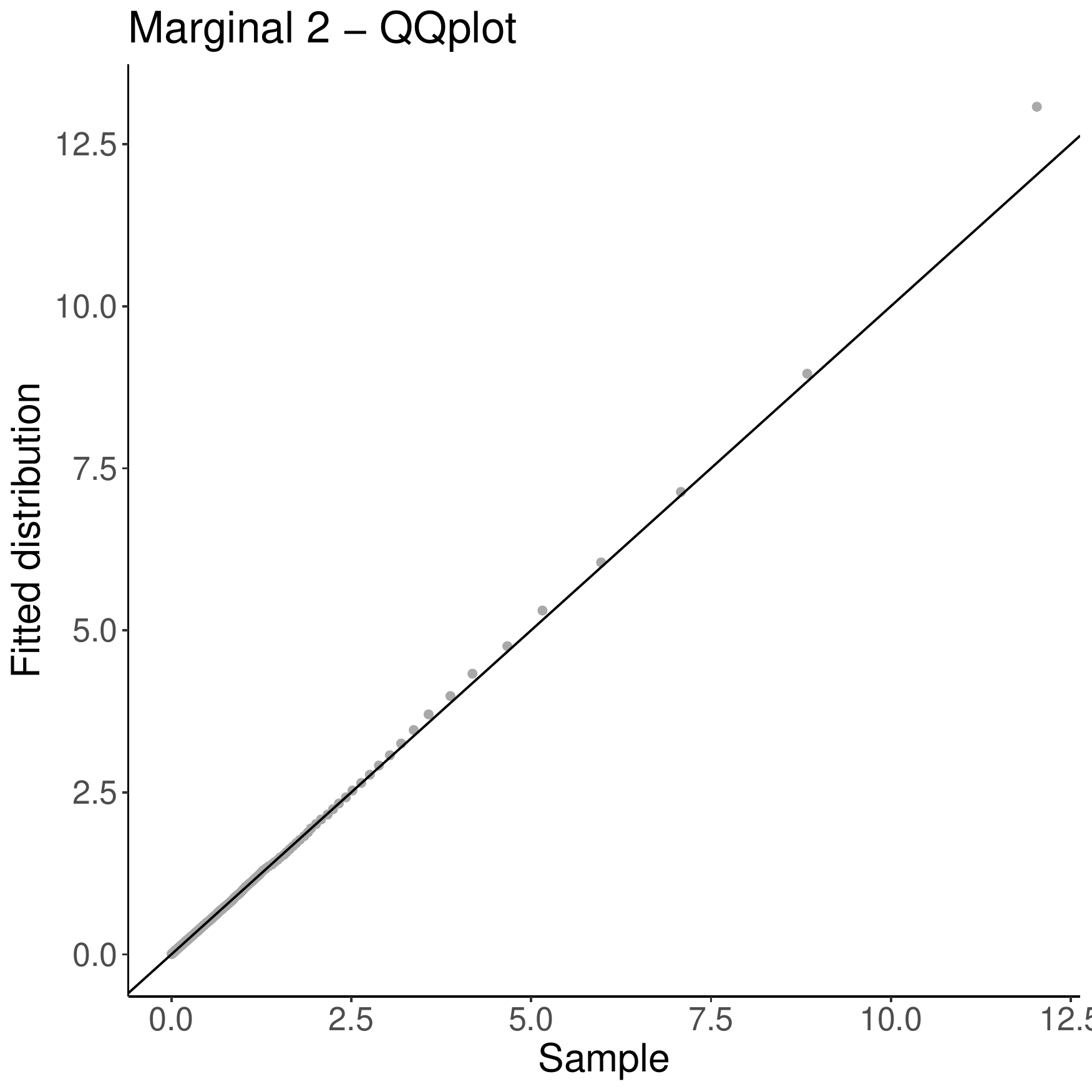}
  	\end{subfigure}
  	\begin{subfigure}{0.22\textwidth}
  		\includegraphics[width=\textwidth]{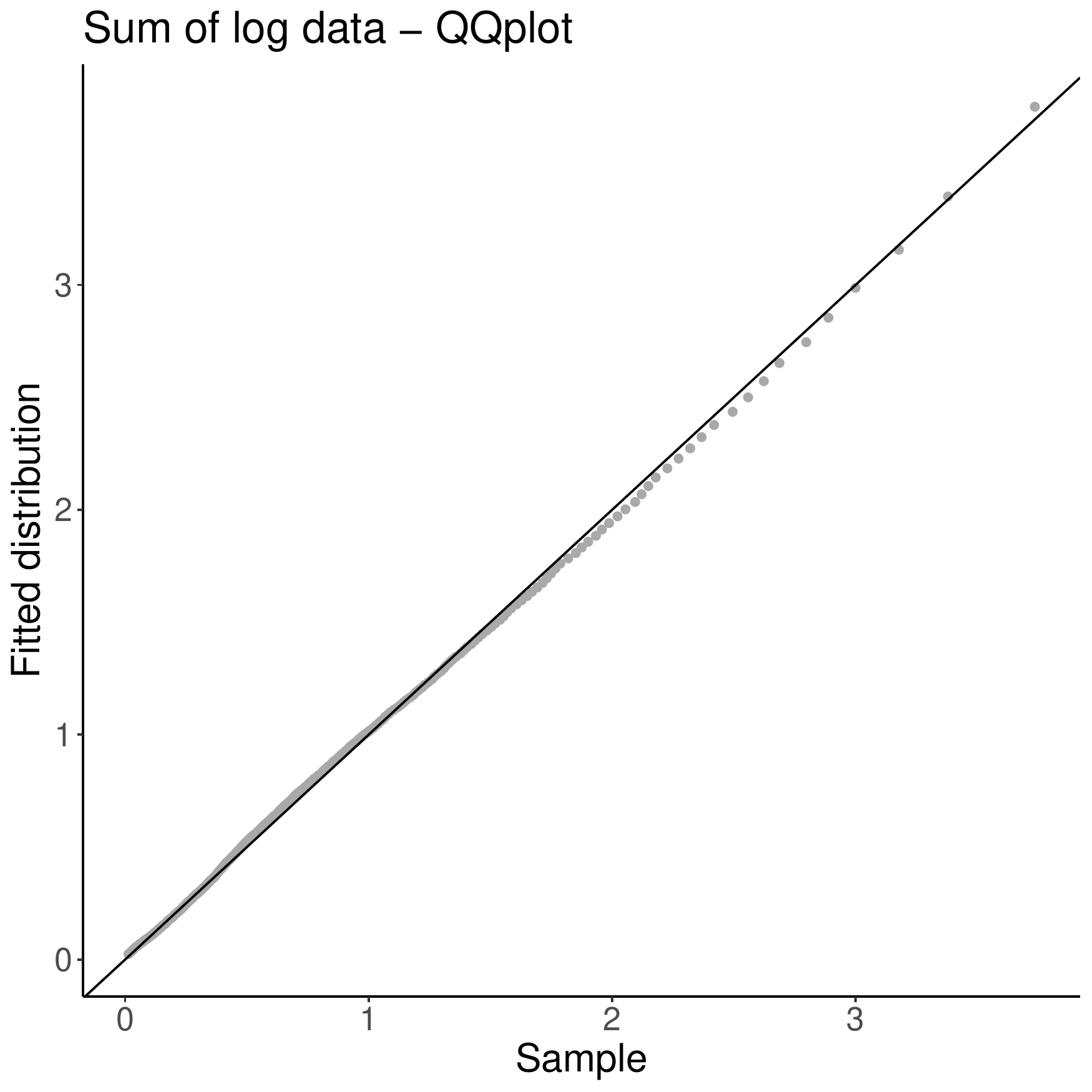}
  	\end{subfigure}
  \caption{QQ plots of simulated Mardia Pareto sample versus fitted multivariate matrix--Pareto distribution using Algorithm~\ref{alg:bivariate}.  } \label{fig:MardiaQQ}
\end{figure}

\begin{figure}[hbt]
	\centering
  	\begin{subfigure}{0.22\textwidth}
  		\includegraphics[width=\textwidth]{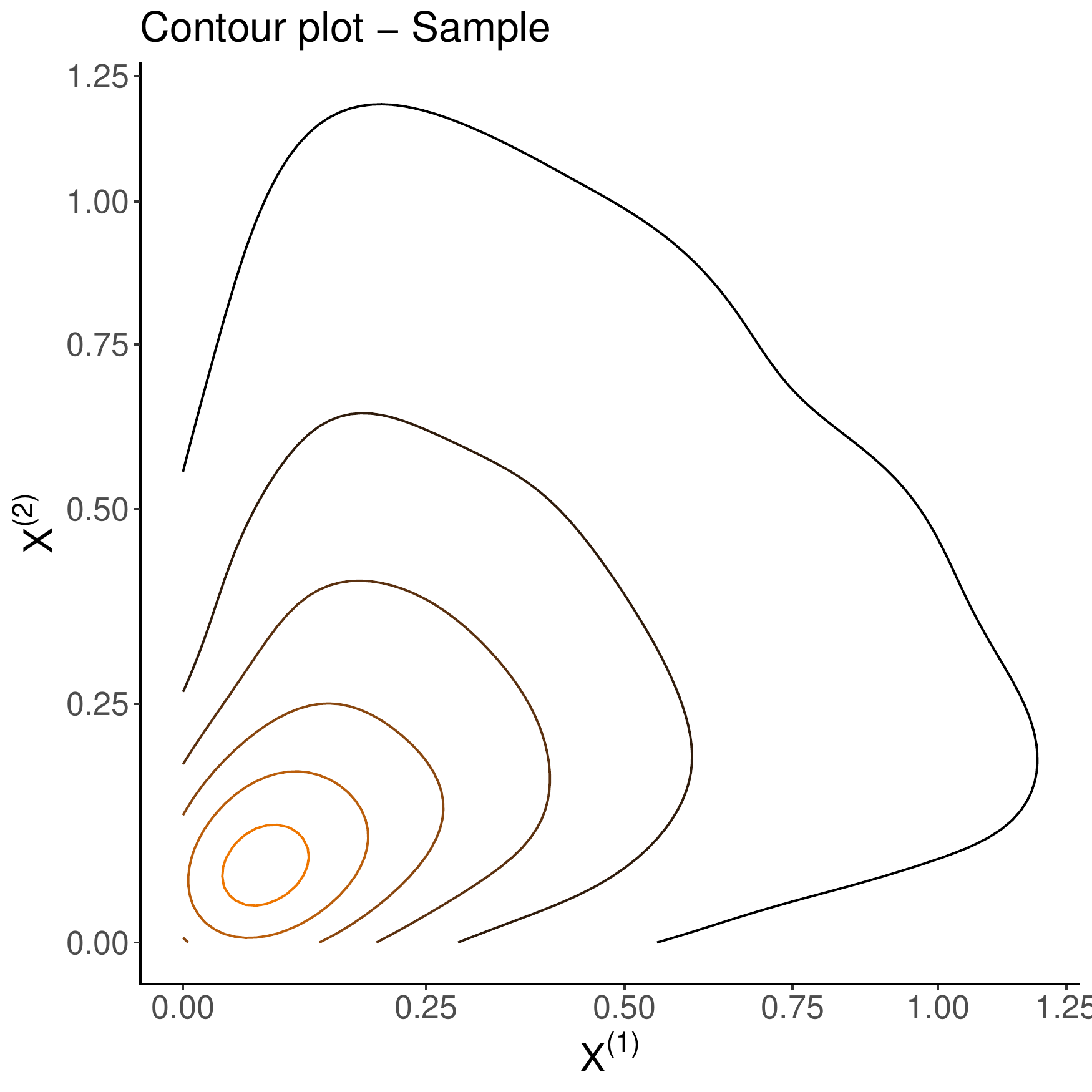}
  	\end{subfigure}
  	\begin{subfigure}{0.22\textwidth}
  		\includegraphics[width=\textwidth]{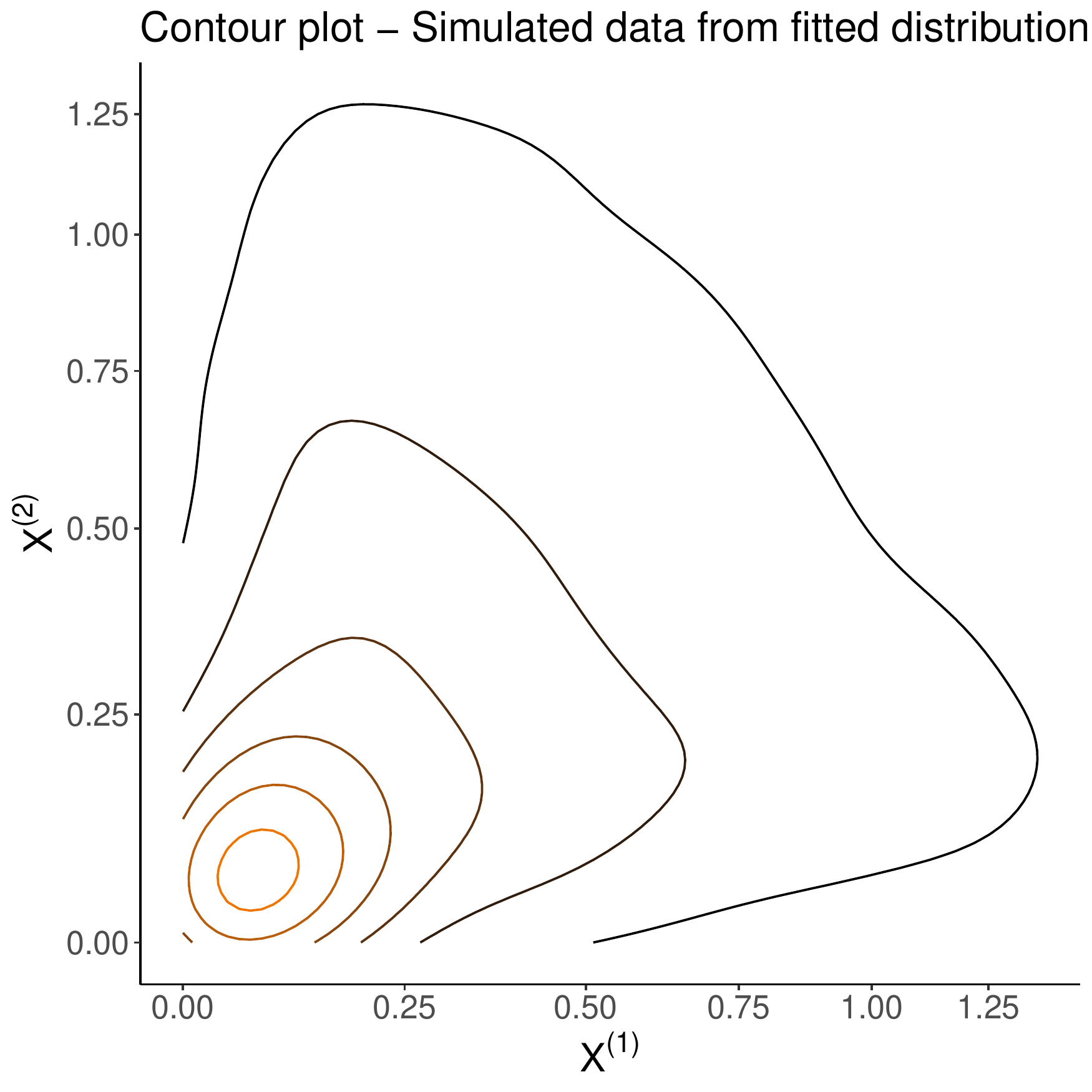}
  	\end{subfigure}
  \caption{Contour plot of simulated Mardia Pareto sample (left) and  contour plot of a simulated sample from the  Matrix--Pareto distribution fitted using Algorithm~\ref{alg:bivariate} (right). }\label{fig:Marcontour}
\end{figure}

\begin{figure}[hbt]
	\centering
  	\begin{subfigure}{0.22\textwidth}
  		\includegraphics[width=\textwidth]{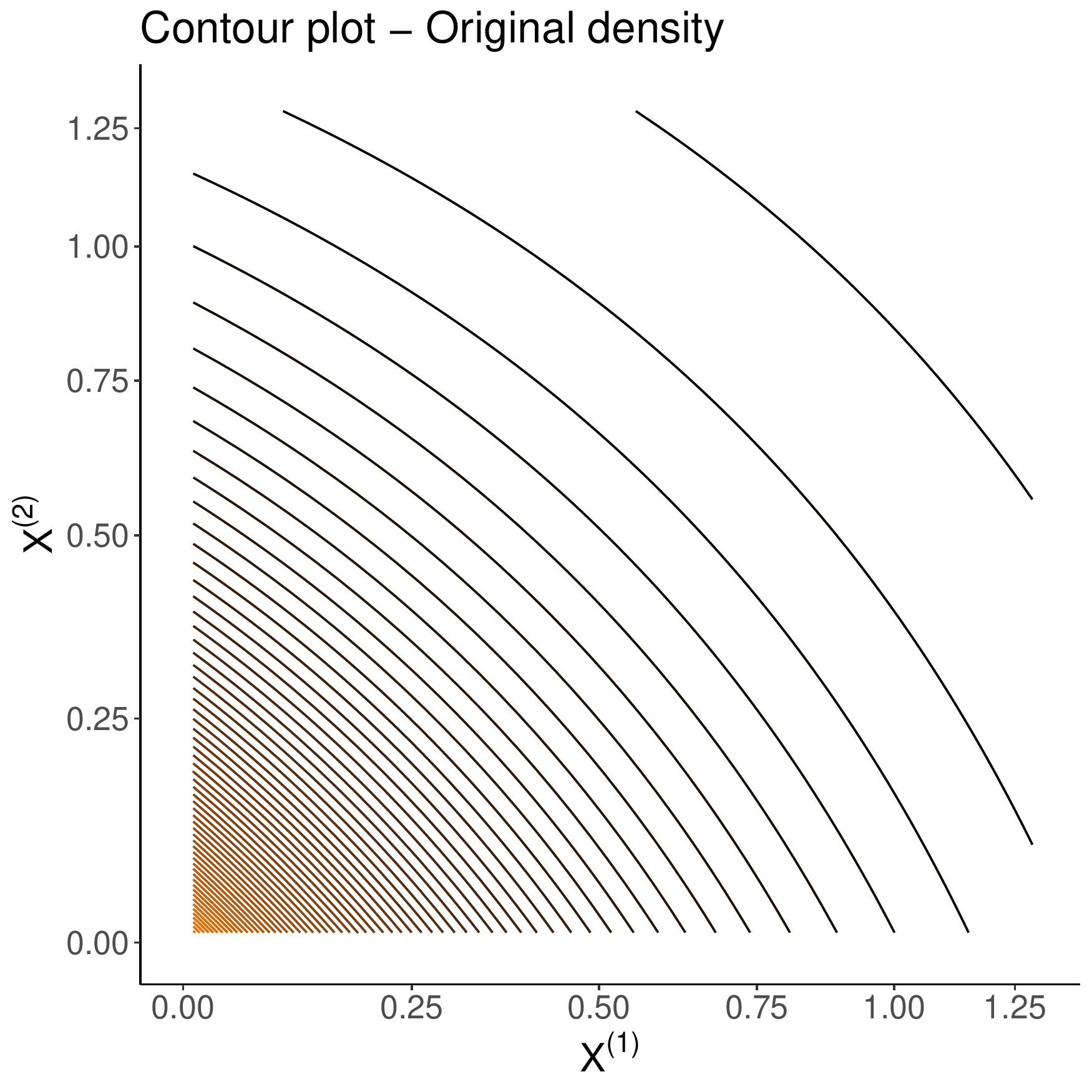}
  	\end{subfigure}
  	\begin{subfigure}{0.22\textwidth}
  		\includegraphics[width=\textwidth]{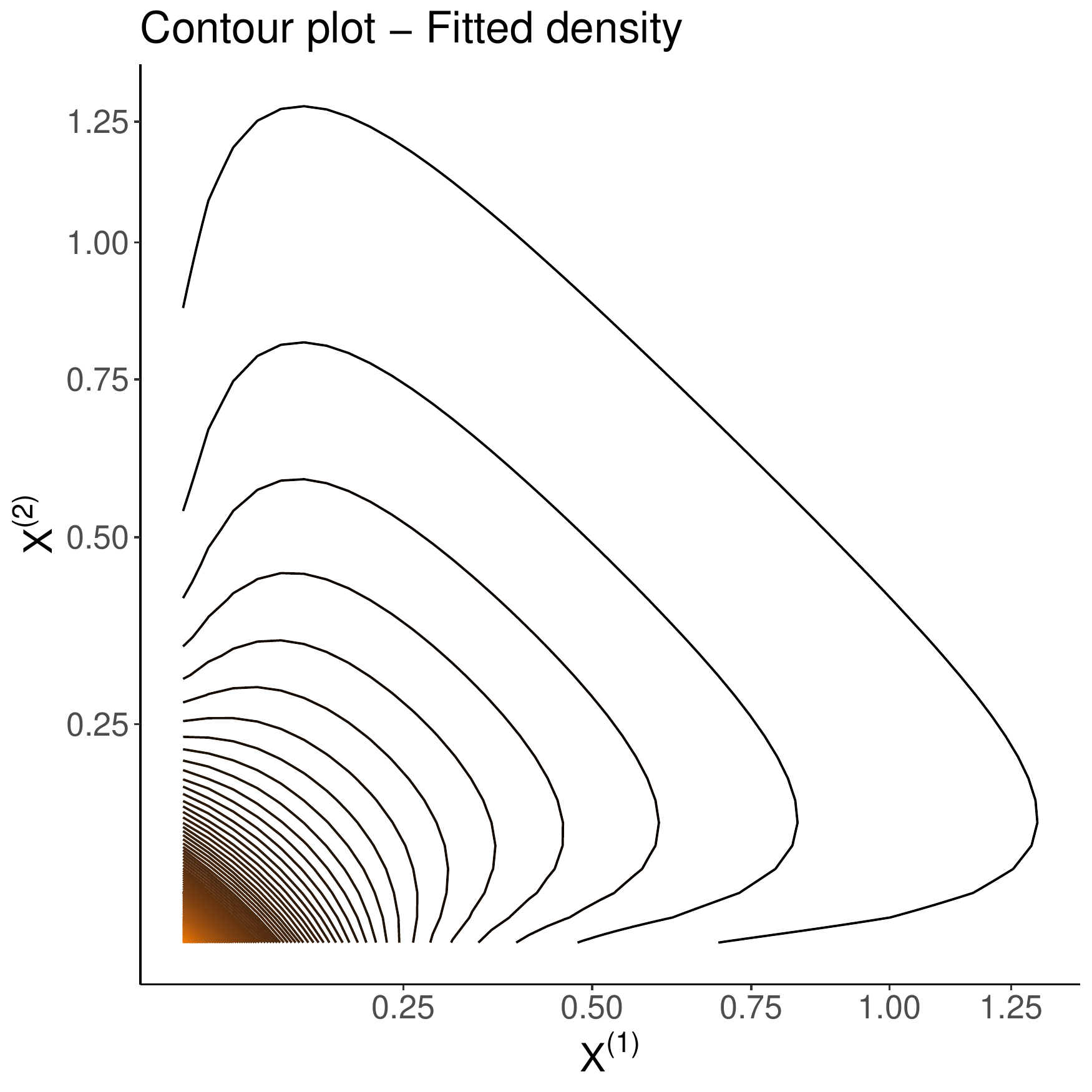}
  	\end{subfigure}
  \caption{Contour plot of original Mardia Pareto (left) and  contour plot of Matrix--Pareto distribution fitted using Algorithm~\ref{alg:bivariate} (right). }\label{fig:Marcontour2}
\end{figure}
\end{example}
The concrete structure of the intensity and reward matrix underlying Algorithm~\ref{alg:bivariate} restricts its application to tail-independent models. While in the previous example this was justified, in situations with tail-dependent data one should rather look for fits in the general MPH$^*$ class by using Algorithm~\ref{alg:MEM} on the transformed data. The following example illustrates such an approach.

\begin{example}\normalfont (Gumbel copula with matrix--Pareto marginals)
We generated an i.i.d.\ sample of size $5\,000$ from a three-dimensional random vector with first marginal being a matrix--Pareto with parameters 
\begin{gather*} 
 {\bfpi_1}=\left(
0.2,\, 0.8\right)\,, \\ 
{\bfT_1}=\left( \begin{array}{ccc}
-2 & 0  \\ 
 0 & -3 \\
\end{array} \right) \,, \\
\beta_{1} = 1 \,,
\end{gather*}
second marginal being a matrix--Pareto with parameters 
\begin{gather*} 
 {\bfpi_2}=\left(
0.5, 0.5\right)\,, \\ 
{\bfT_2}=\left( \begin{array}{ccc}
-2 & 0 \\ 
 0 & -1.5 
\end{array} \right)  \,,\\
\beta_{2} = 1 \,,
\end{gather*}
and the third marginal being a conventional Pareto with shape parameter $2.5$, and a Gumbel copula with parameter $\theta=4$. The choice of a Gumbel copula instead of a multivariate matrix--Pareto model based on an MPH$^*$ construction is to show that the algorithm can be used to model any type of dependence structure. The Gumbel copula  is known to have positive tail dependence. This distribution has theoretical numerical values $\E(\bfX)=(0.6, 1.5, 0.6667)^{\prime}$, $\lambda_U=2-2^{1/4}\approx 0.8108$ and $\rho_{\tau}=0.75$, and the real part of the eigenvalues that determine the heaviness of the tails are $\lambda_{1}^{(\max)}=-2$,  $\lambda_{2}^{(\max)}=-1.5$ and $\lambda_{3}^{(\max)}=-2.5$. 
The simulated sample has numerical values $\hat{\E}(\bfX)= \left( 0.6035  , 1.4721, 0.6732 \right)^{\prime} $, $\hat{\lambda}_U (X_1, X_2)=0.8142$, $\hat{\lambda}_U (X_1, X_3)=0.7429$, $\hat{\lambda}_U (X_2, X_3)=0.7857$,  $\hat{\rho}_\tau(X_1,X_2)=0.7513$,  $\hat{\rho}_\tau(X_1,X_3)=0.7526$ and $\hat{\rho}_\tau(X_2,X_3)=0.7544$. 
Then, we fitted a multivariate matrix--Pareto distribution using Algorithm~\ref{alg:MEM} with $p=6$ and $5\,000$ steps on the transformed data (running time $6\,501$ seconds), obtaining the following parameters: 
\begin{gather*} 
 \hat{\bfalp}=\left(
0.0062\,, 0.0606\,, 0.4394\,, 0.0600\,, 0.0608\,, 0.3731 \right)\,, \\ 
\hat{\bfT}=\left( \begin{array}{cccccc}
-16.9867 & 1.7532 & 1.3902 & 3.2835 & 1.0736 & 2.1224 \\ 
 0.2710 & -1.3282 & 0.2835 & 0.0757 & 0.0562 & 0.2987 \\ 
0.3714 & 0.0161 & -2.8735 & 0.6658 & 0.8798 & 0.9208 \\ 
1.1742 & 0.5387 & 0.2421 &-4.0067 & 0.5156 & 0.5571 \\ 
0.0786 & 0.5130 & 0.2916 & 1.1454 & -4.2894 & 1.5148 \\ 
2.0244 & 0.2912 & 0.7007 & 0.7777 & 2.0237 & -6.2507
\end{array} \right) \,, \\ 
\hat{\bfR}= \left( \begin{array}{ccc}
0.9414 &0.0558 &0.0028 \\ 
0.2856 &0.4594 &0.2549 \\ 
0.2934 &0.5233 &0.1833 \\ 
0.0906 &0.5548 &0.3546 \\ 
0.4230 &0.0750 &0.5020 \\ 
0.0748 &0.5992 &0.3260
\end{array} \right)\,.
\end{gather*}

The real part of the eigenvalues with largest real part of the sub--intensity matrices of the marginal distributions are $\lambda_{1}^{(\max)}=-2.5370$,  $\lambda_{2}^{(\max)}=-1.5941$ and $\lambda_{3}^{(\max)}=-2.4167$, which are close to the ones of the original distribution. 
The fitted distribution has first moment 
$	\E(\bfX)= \left(
 0.6006, 1.5344, 0.6850 \right)^{\prime} $. Moreover, we estimated $\lambda_U$ and $\rho_{\tau}$ via simulation, obtaining $\hat{\lambda}_U(X_1,X_2)=0.70403 $, $\hat{\lambda}_U(X_1,X_3)=0.736438 $, $\hat{\lambda}_U(X_2,X_3)=0.672417 $,  $\hat{\rho}_{\tau}(X_1,X_2)=0.6924$, $\hat{\rho}_{\tau}(X_1,X_3)=0.7470$ and $\hat{\rho}_{\tau}(X_2,X_3)=0.7728$. Comparing all numerical properties, together with the QQ plots (see Figure~\ref{fig:MP3qq}) and contour plots (see Figure~\ref{fig:MP3contour}), we see that this algorithm also recovers relatively well the structure of the data. 

\begin{figure}[H]
	\centering
  	\begin{subfigure}{0.22\textwidth}
  		\includegraphics[width=\textwidth]{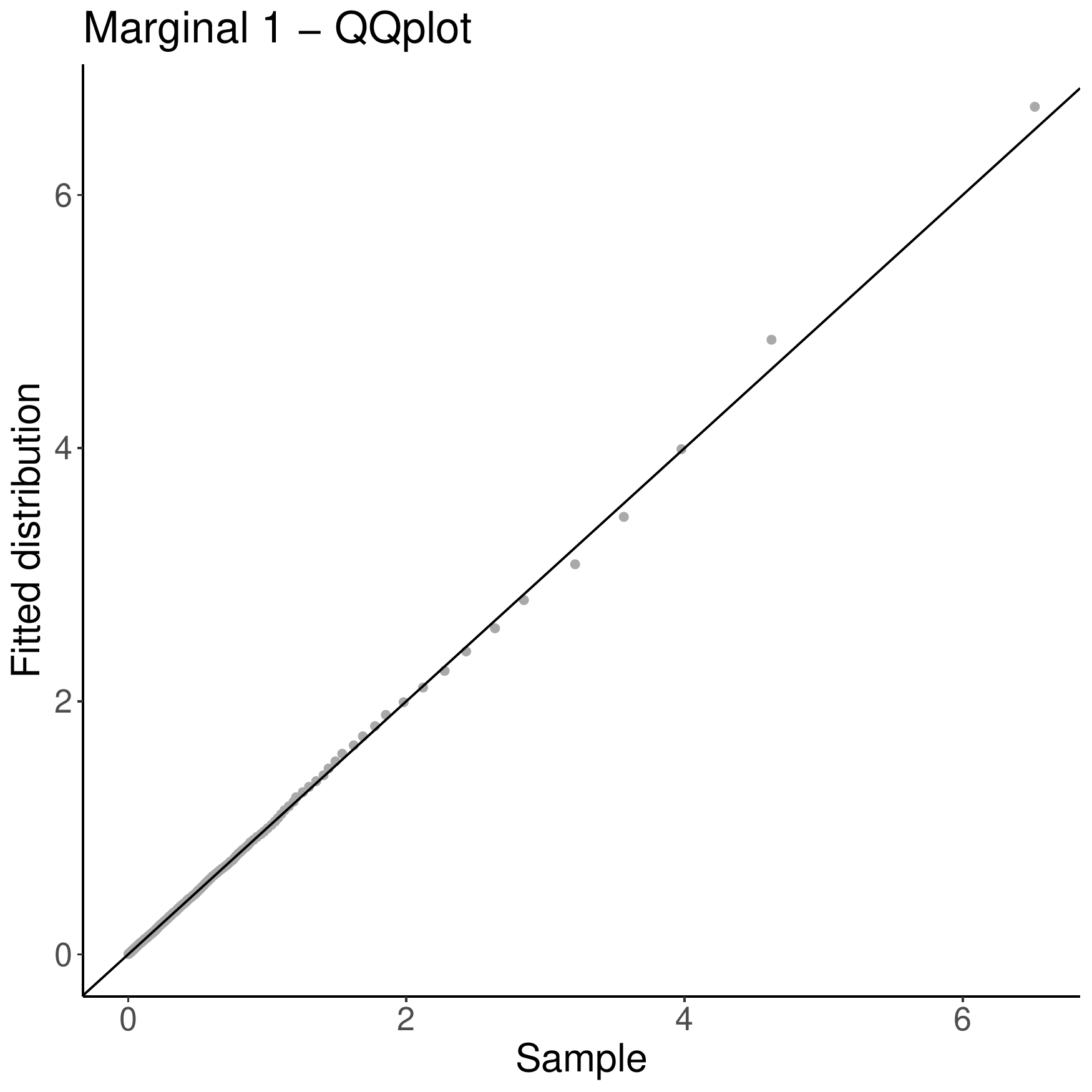}
  	\end{subfigure}
  	\begin{subfigure}{0.22\textwidth}
  		\includegraphics[width=\textwidth]{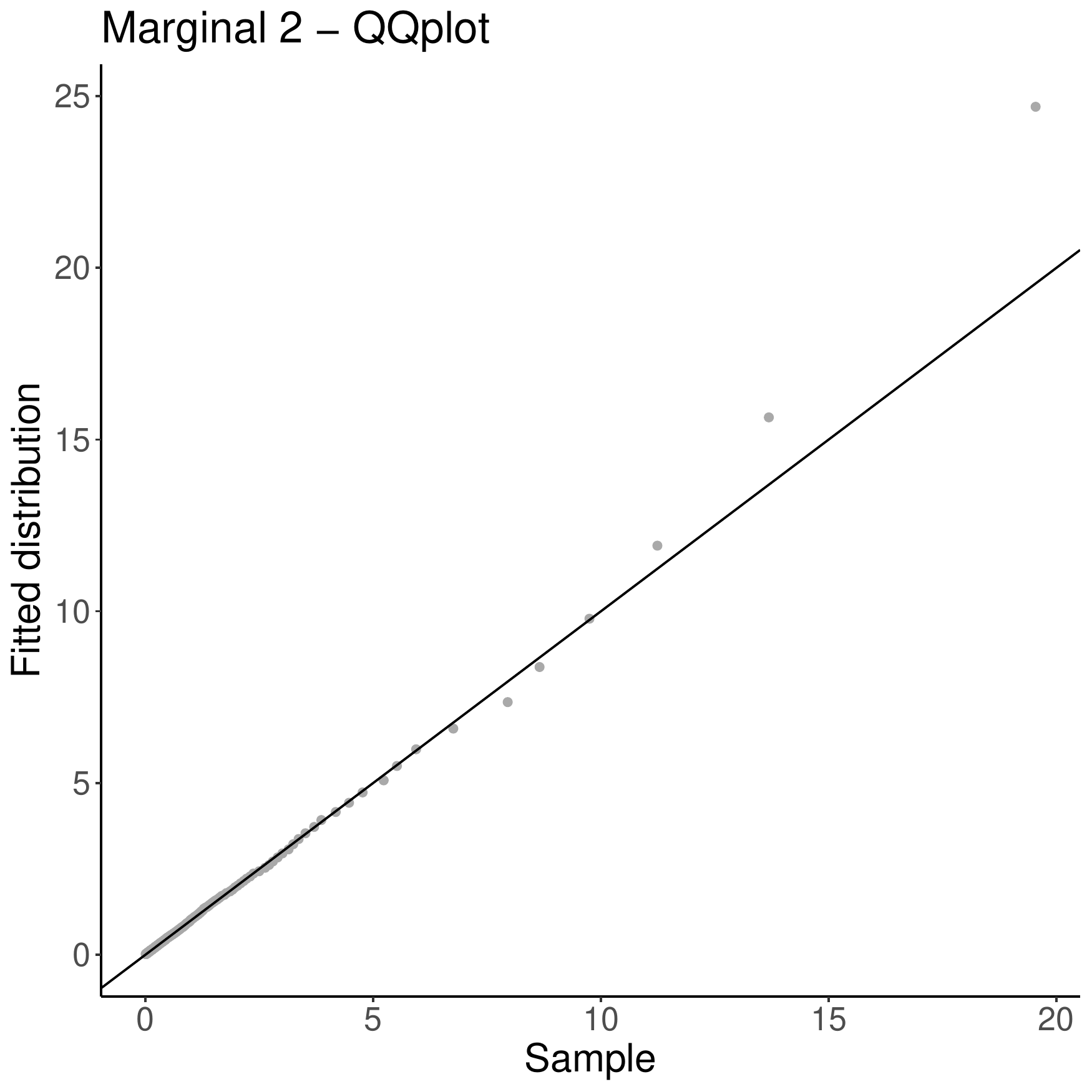}
  	\end{subfigure}
  	\begin{subfigure}{0.22\textwidth}
  		\includegraphics[width=\textwidth]{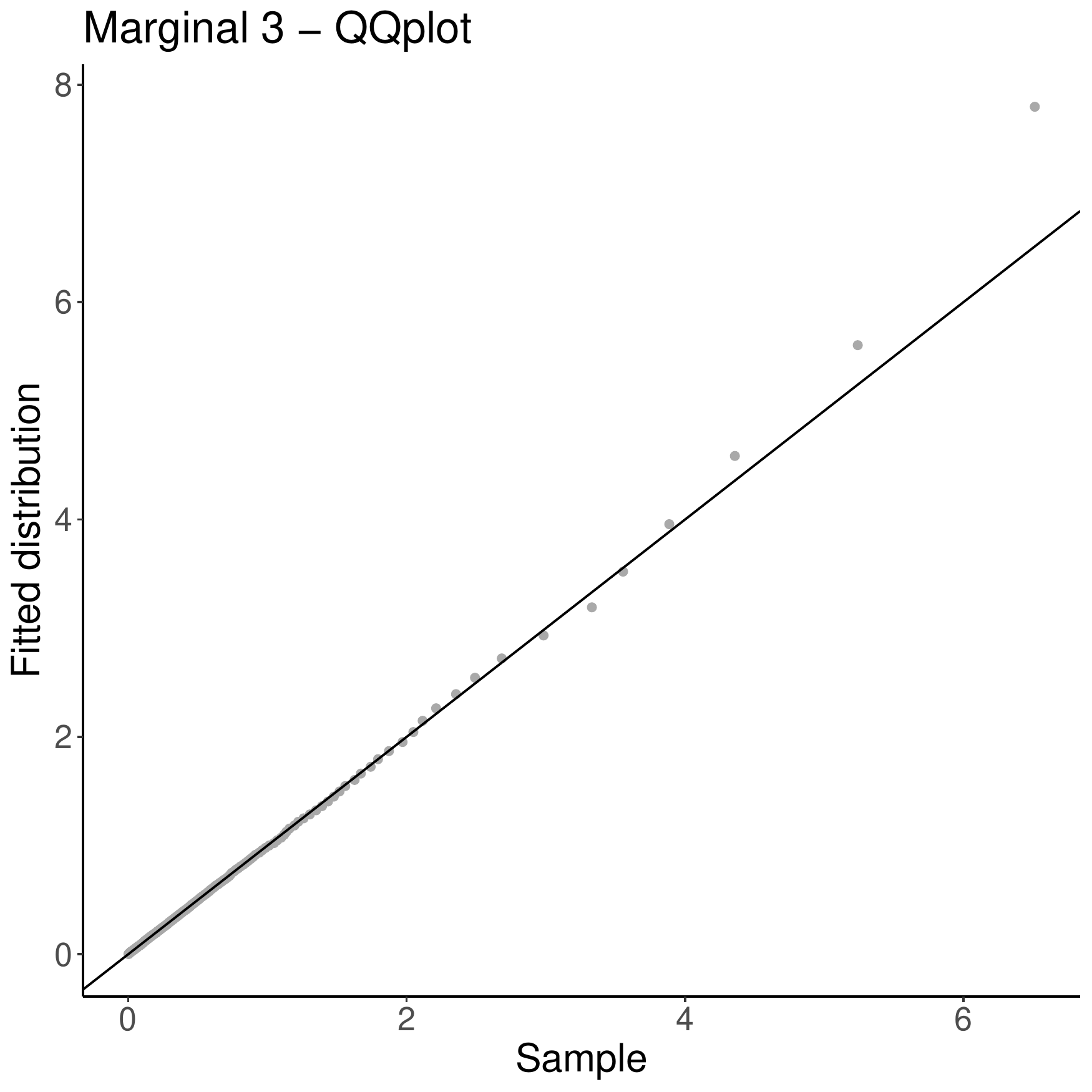}
  	\end{subfigure}
  	\begin{subfigure}{0.22\textwidth}
  		\includegraphics[width=\textwidth]{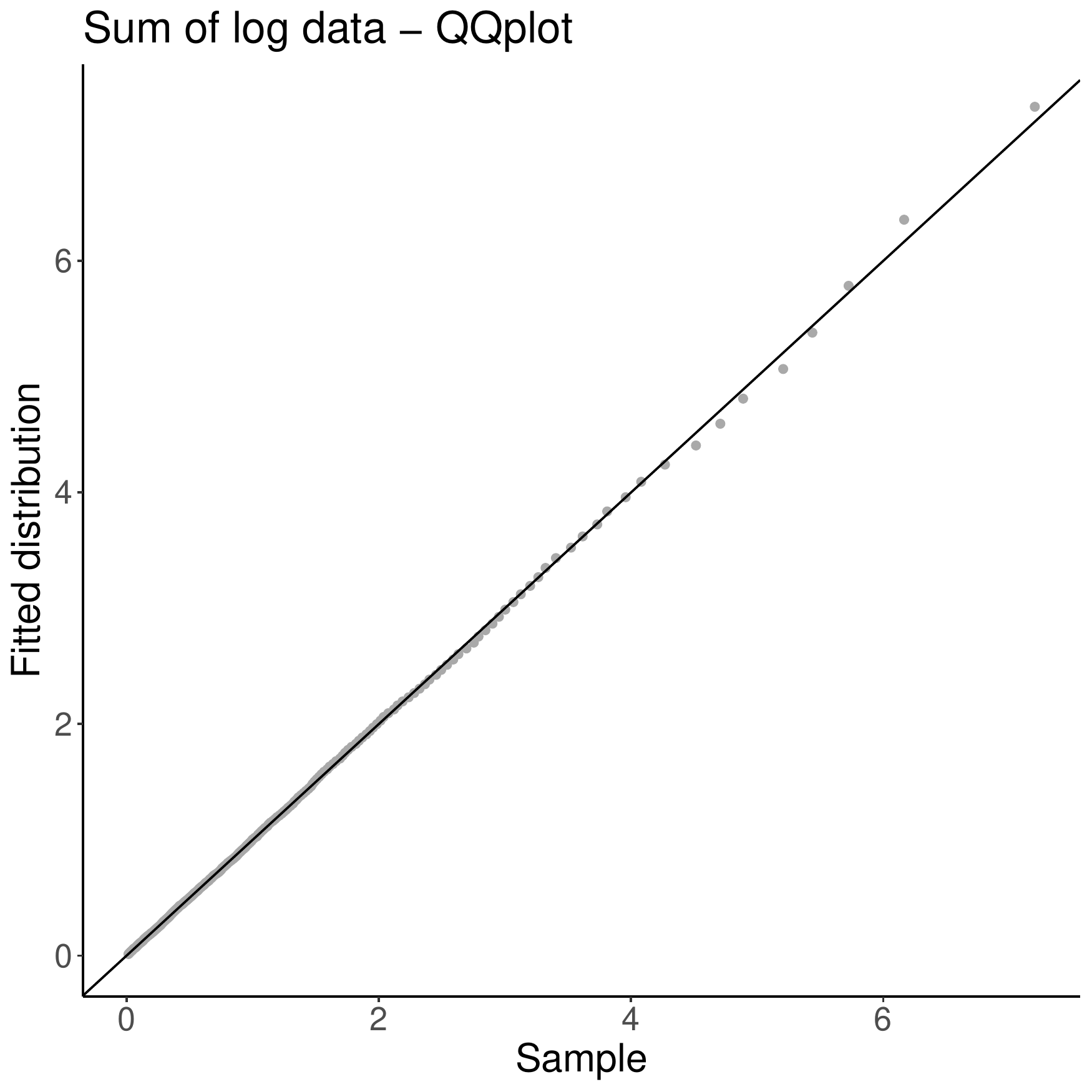}
  	\end{subfigure}
  \caption{QQ plots of simulated sample versus fitted multivariate matrix--Pareto distribution using Algorithm~\ref{alg:MEM}. } \label{fig:MP3qq}
\end{figure}

\begin{figure}[hbt]
	\centering
  	\begin{subfigure}{0.22\textwidth}
  		\includegraphics[width=\textwidth]{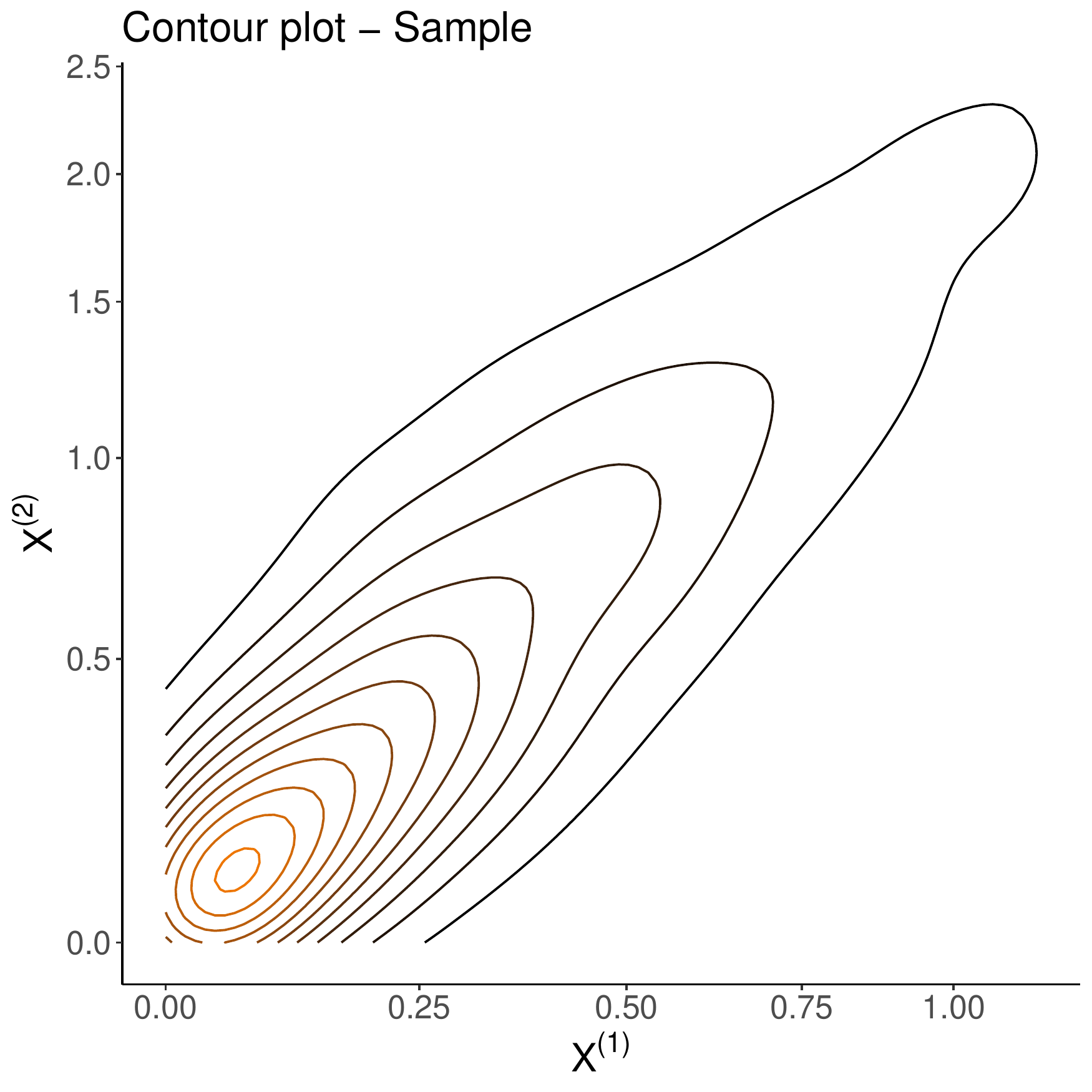}
  	\end{subfigure}
  	\begin{subfigure}{0.22\textwidth}
  		\includegraphics[width=\textwidth]{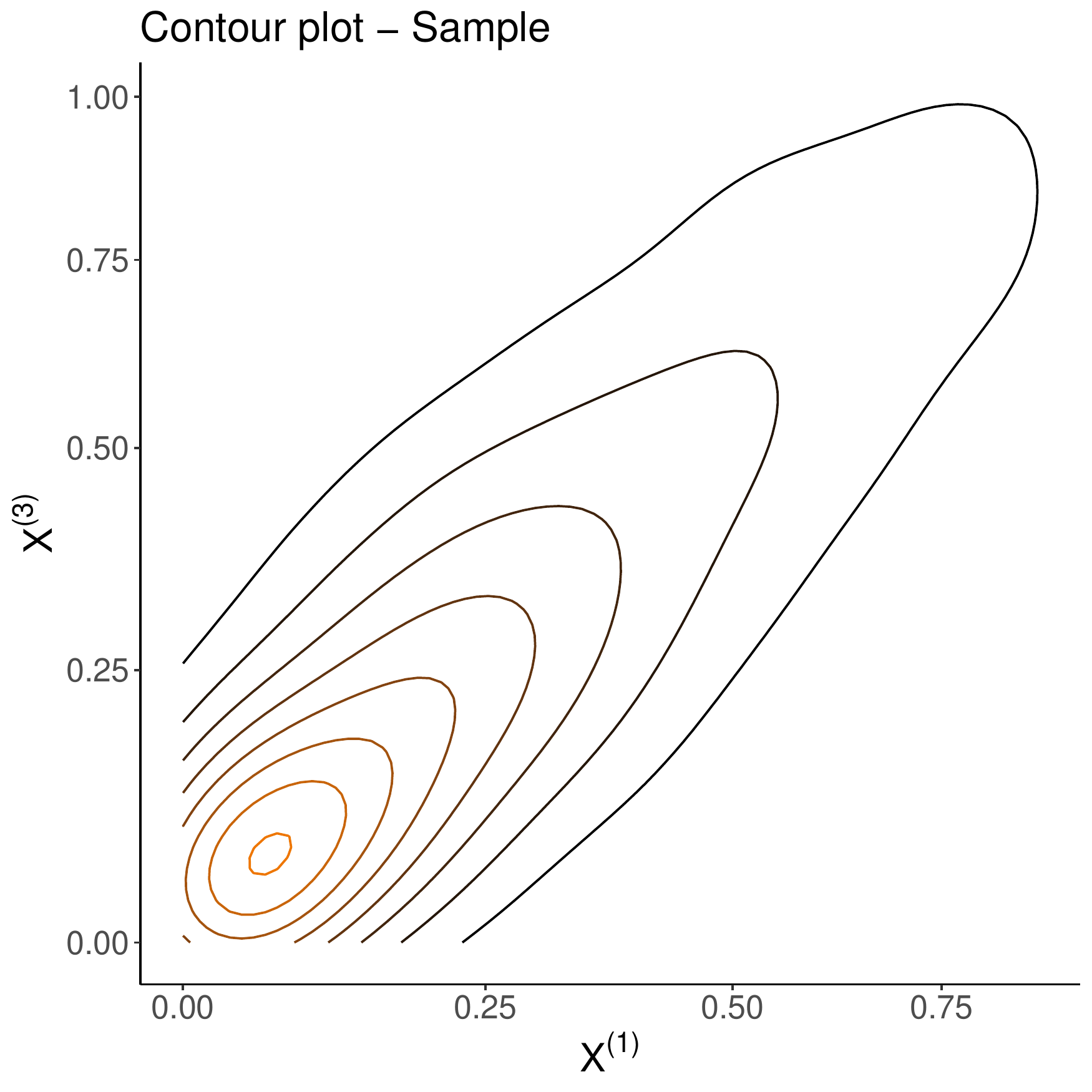}
  	\end{subfigure}
  	\begin{subfigure}{0.22\textwidth}
  		\includegraphics[width=\textwidth]{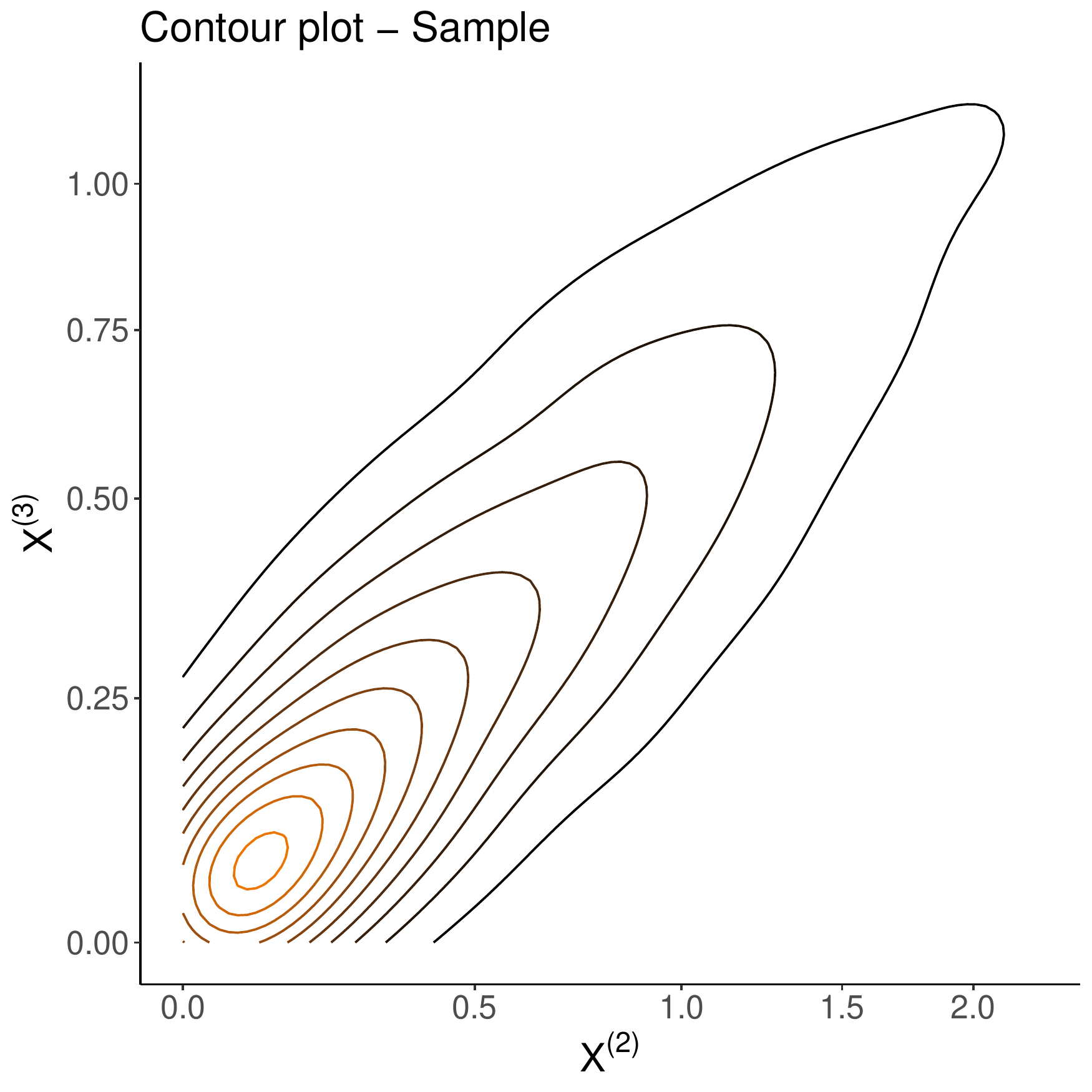}
  	\end{subfigure} \\
  	\begin{subfigure}{0.22\textwidth}
  		\includegraphics[width=\textwidth]{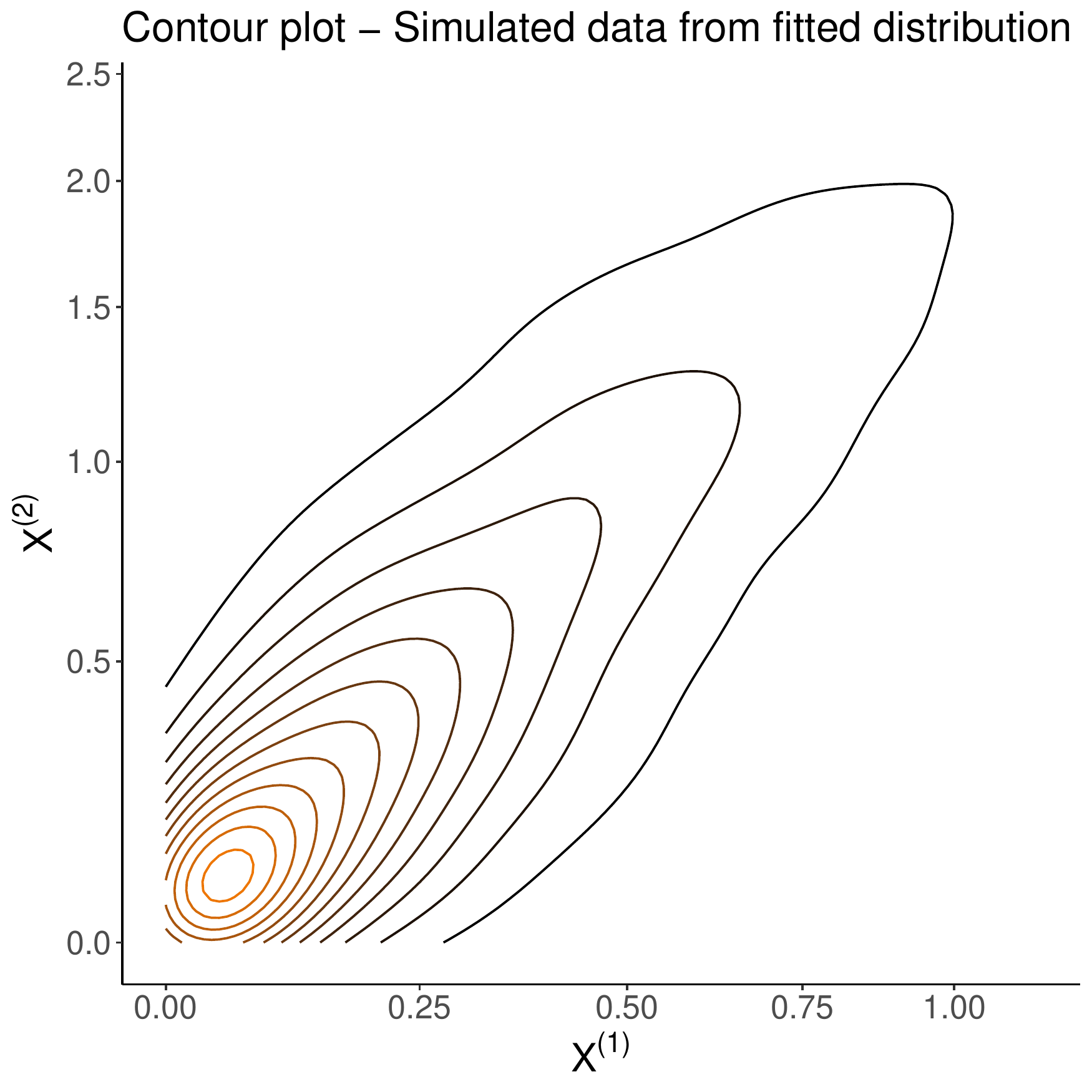}
  	\end{subfigure}
  	\begin{subfigure}{0.22\textwidth}
  		\includegraphics[width=\textwidth]{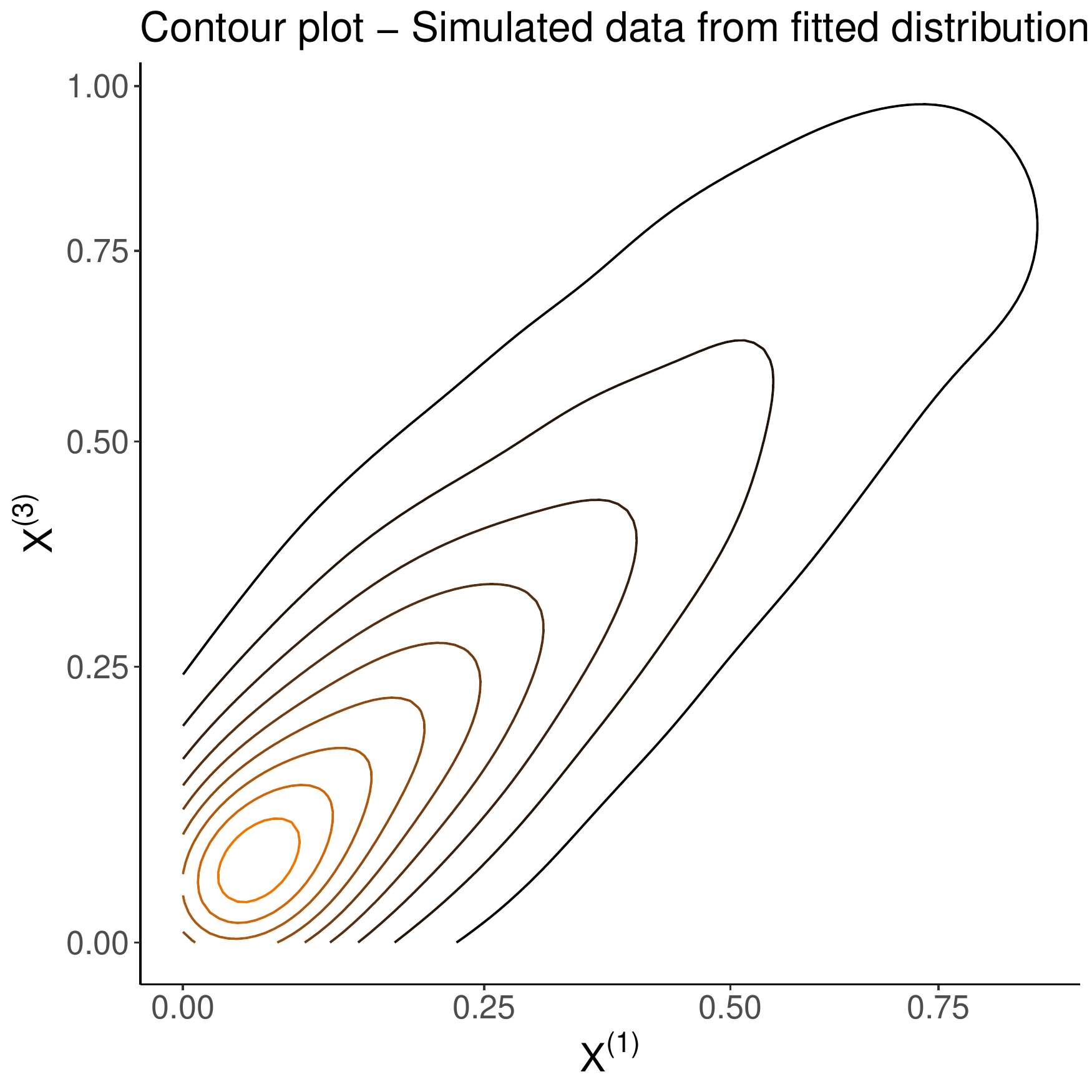}
  	\end{subfigure}
  	\begin{subfigure}{0.22\textwidth}
  		\includegraphics[width=\textwidth]{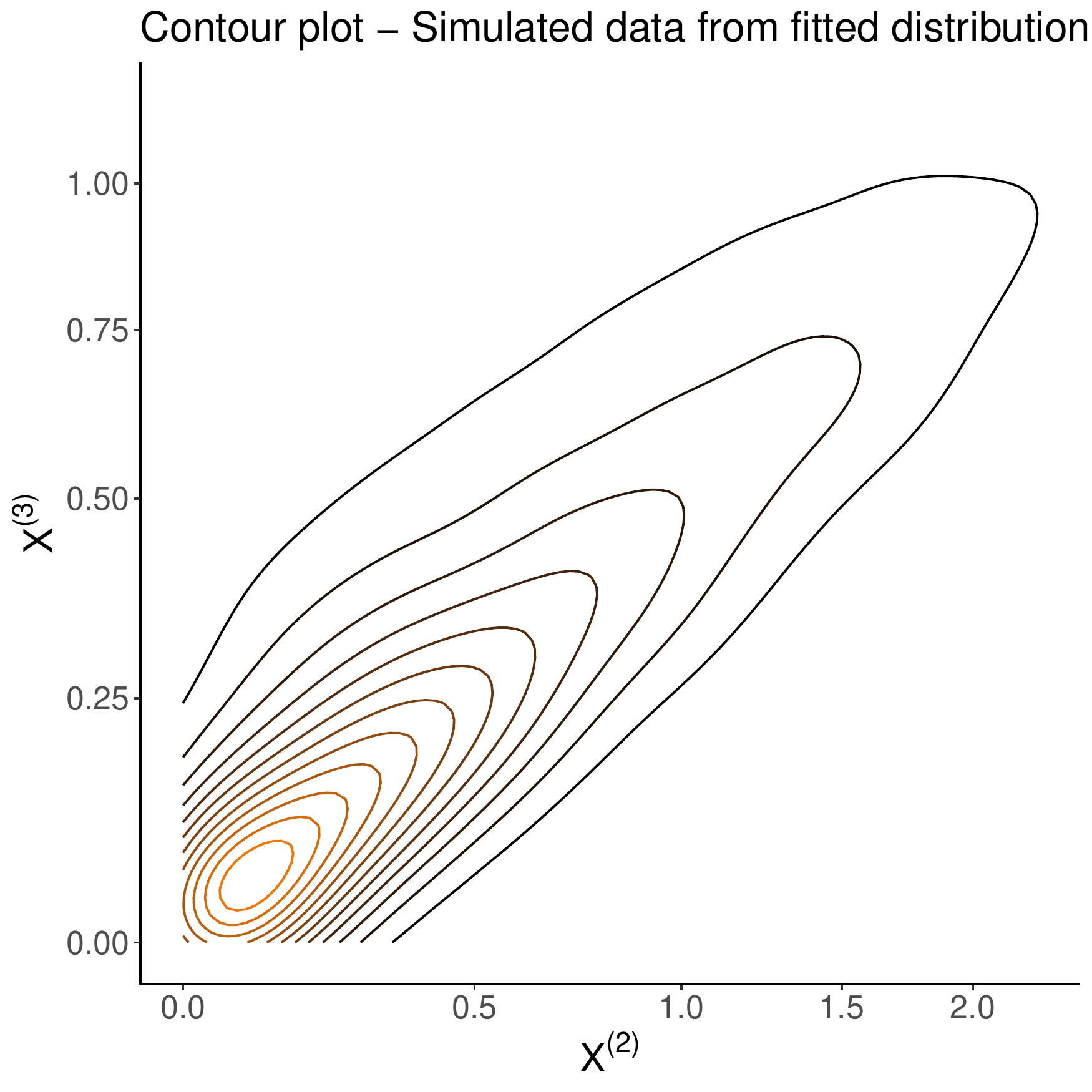}
  	\end{subfigure}
  \caption{Contour plots of simulated sample from Gumbel copula with matrix-Pareto marginals (top) and  contour plots of a simulated sample from the  Matrix--Pareto distribution fitted using Algorithm~\ref{alg:MEM} (bottom). }\label{fig:MP3contour}
\end{figure}

 We would like to comment on the relative inaccuracy of the obtained estimate for $\lambda_{1}^{(\max)}$, especially when compared to the convincing estimates for the other marginals. The first marginal here is in fact a mixture of Pareto distributions with shape parameters $2$ and $3$, and mixing probabilities $0.2$ and $0.8$, respectively. That is, a major proportion of the data for the first marginal are to be expected to stem from the Pareto distribution with shape parameter $3$. Since Algorithm~\ref{alg:MEM} fits body and tail at the same time, one cannot expect the resulting estimate for the tail to be as good as techniques of classical extreme value theory, which focusses on the tail fit only. To illustrate this point, we also applied the algorithm to data simulated from changed initial probabilities in the first marginal according to $\bfpi_{1} =(0.8, \, 0.2)$ (and otherwise identical parameters). Now most of the data points will stem from the heavier distribution in the mixture. Indeed, the estimate for $\lambda_{1}^{(\max)}$ in that case turns out to be $-2.0152$, while the tail estimate for the other marginals remains almost at the same value ($-1.5966$ and $-2.4154$). 

\end{example}

Next we present a parameter--dependent example with real data, employing Algorithm~\ref{alg:BivIPH}.
\begin{example}\normalfont {(Danish fire insurance data)}
Consider the famous Danish fire insurance claim data set (see e.g.\ \cite{grun2019extending}). We propose here a bivariate matrix--Pareto distribution as a model for the components building and content with observations in the set $(1,\infty) \times (1,\infty)$.  
To that end, we first translate the sample to the origin, thus
the sample has numerical values $\hat{\E}(\bfX)= \left( 3.2476  , 4.6856 \right)^{\prime} $ and  $\hat{\rho}_\tau=0.2143$. 
Then, we fit a bivariate matrix--Pareto distribution with $p_1=p_2=2$ using Algorithm~\ref{alg:BivIPH} with $2\,000$ steps (with a step--length of $0.05$ and gradient ascent until the norm of the derivative is less than $0.1$, the running time is $438$ seconds), obtaining the following parameters: 
\begin{gather*} 
 \hat{\bfalp}=\left(
0,\, 1,\, 0,\, 0\right)\,, \\ 
\hat{\bfT}=\left( \begin{array}{cccccc}
-2.6333 & 0 &  0.0005 & 2.6328   \\ 
1.2788 & -3.8336 &  2.5548 & 0  \\
0 & 0 & -10.9822 & 0   \\
 0 & 0 & 2.4131 & -2.4732
\end{array} \right) \,,\\ 
\hat{\bfR}= \left( \begin{array}{c c}
1 & 0 \\ 
1 & 0  \\
0 & 1 \\ 
0 & 1 
\end{array} \right) \,,\\
\hat{\beta}_1= 5.0377\,, \quad
\hat{\beta}_2= 13.4538 \,.
\end{gather*}
The real part of the eigenvalues that determines the heaviness of the tails are $\lambda_{1}^{(\max)}=-2.6333$ and $\lambda_{2}^{(\max)}=-2.4732$, respectively. 
The fitted distribution has mean
$	\E(\bfX)= \left( 3.1698, 4.6803 \right)^{\prime} $. Moreover, estimating $\rho_{\tau}$ via simulation gives  $\hat{\rho}_{\tau}=0.2190$. From the QQ plots (Figure~\ref{fig:qqDanish}) and contour plots (Figure~\ref{fig:DanishContour}), we see that the fitted distribution is a reasonable model for the data.


\begin{figure}[H]
	\centering
  	\begin{subfigure}{0.22\textwidth}
  		\includegraphics[width=\textwidth]{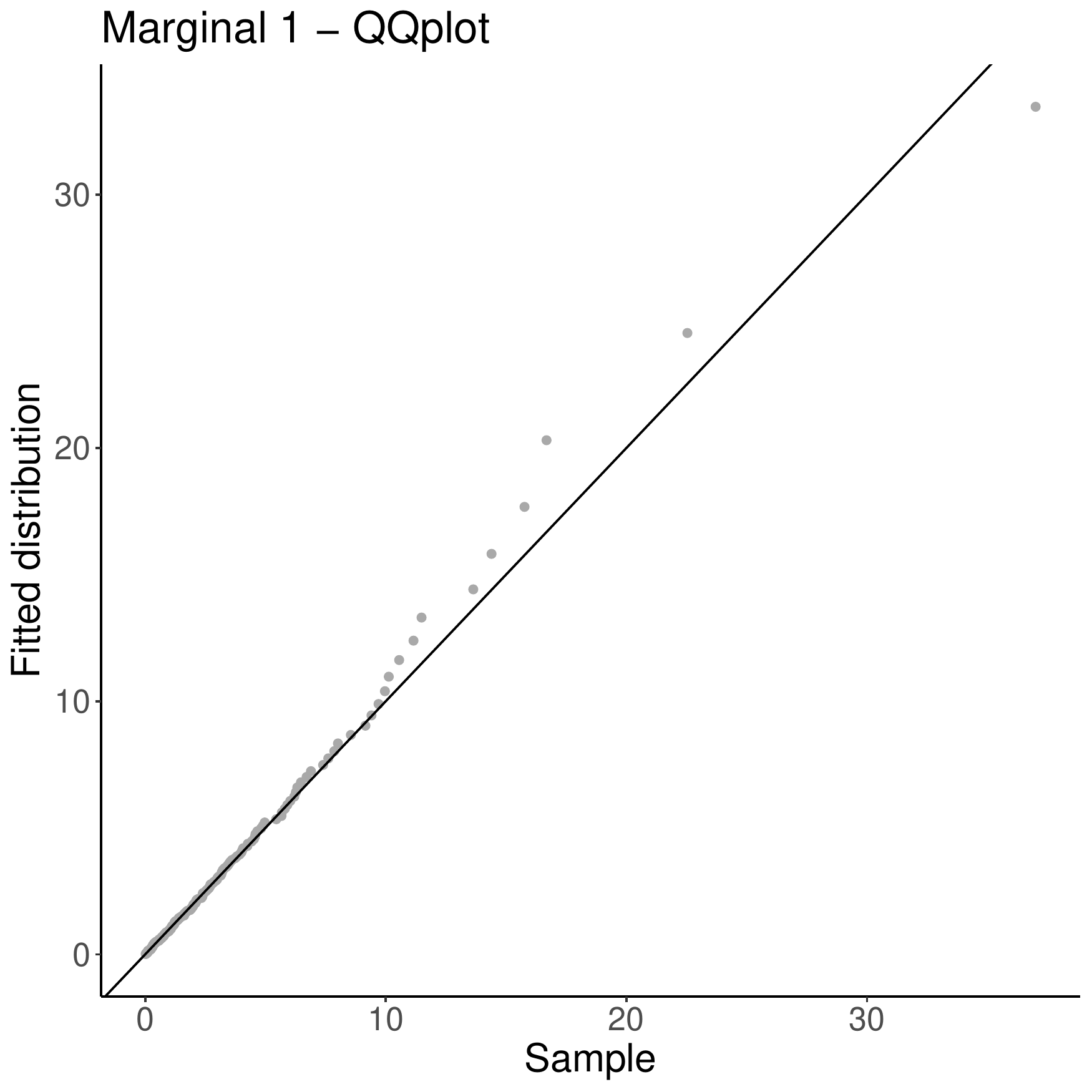}
  	\end{subfigure}
  	\begin{subfigure}{0.22\textwidth}
  		\includegraphics[width=\textwidth]{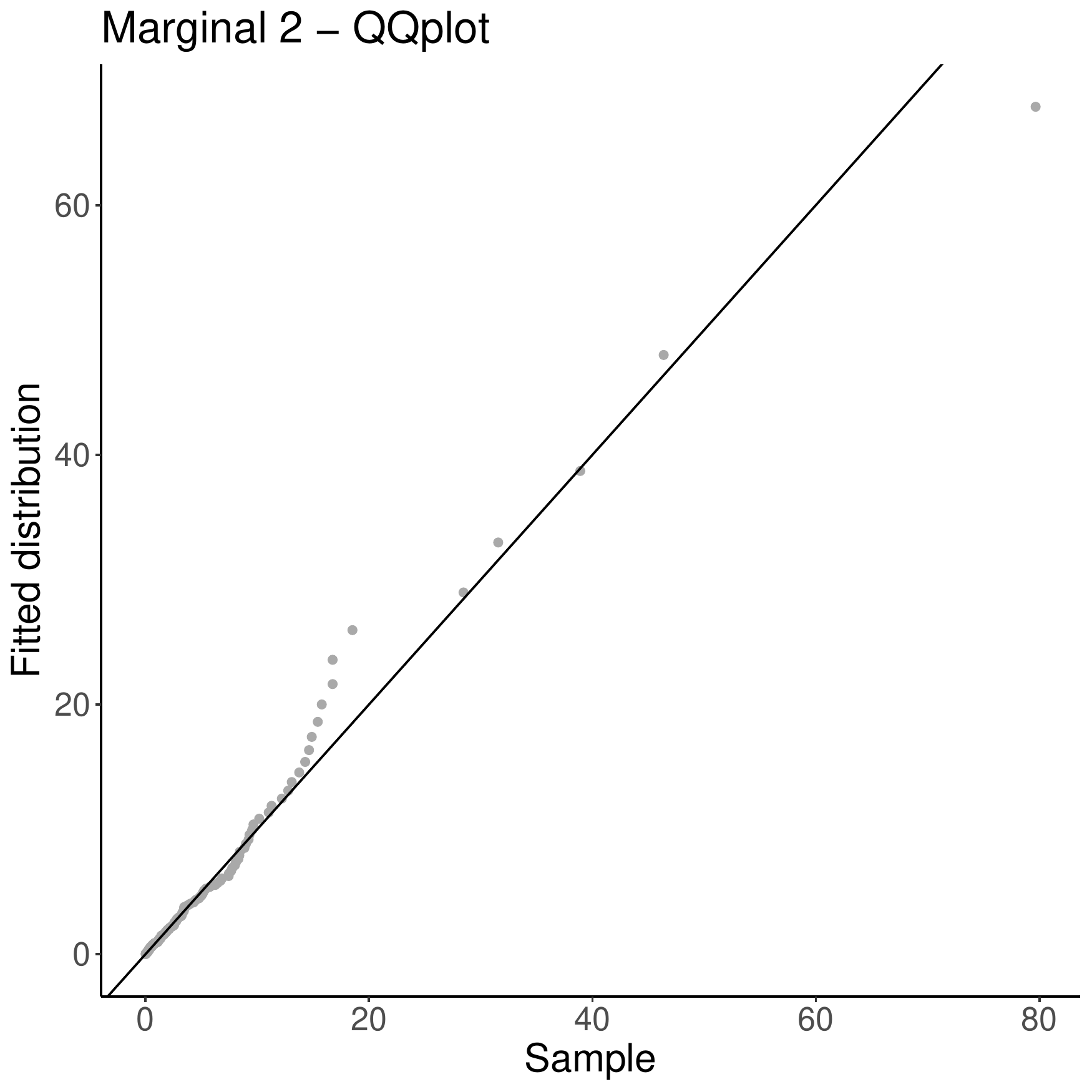}
  	\end{subfigure}
  \caption{QQ plots of Danish fire insurance claim size sample versus fitted bivariate matrix--Pareto distribution using Algorithm~\ref{alg:BivIPH}. } \label{fig:qqDanish}
\end{figure}

\begin{figure}[hbt]
	\centering
  	\begin{subfigure}{0.22\textwidth}
  		\includegraphics[width=\textwidth]{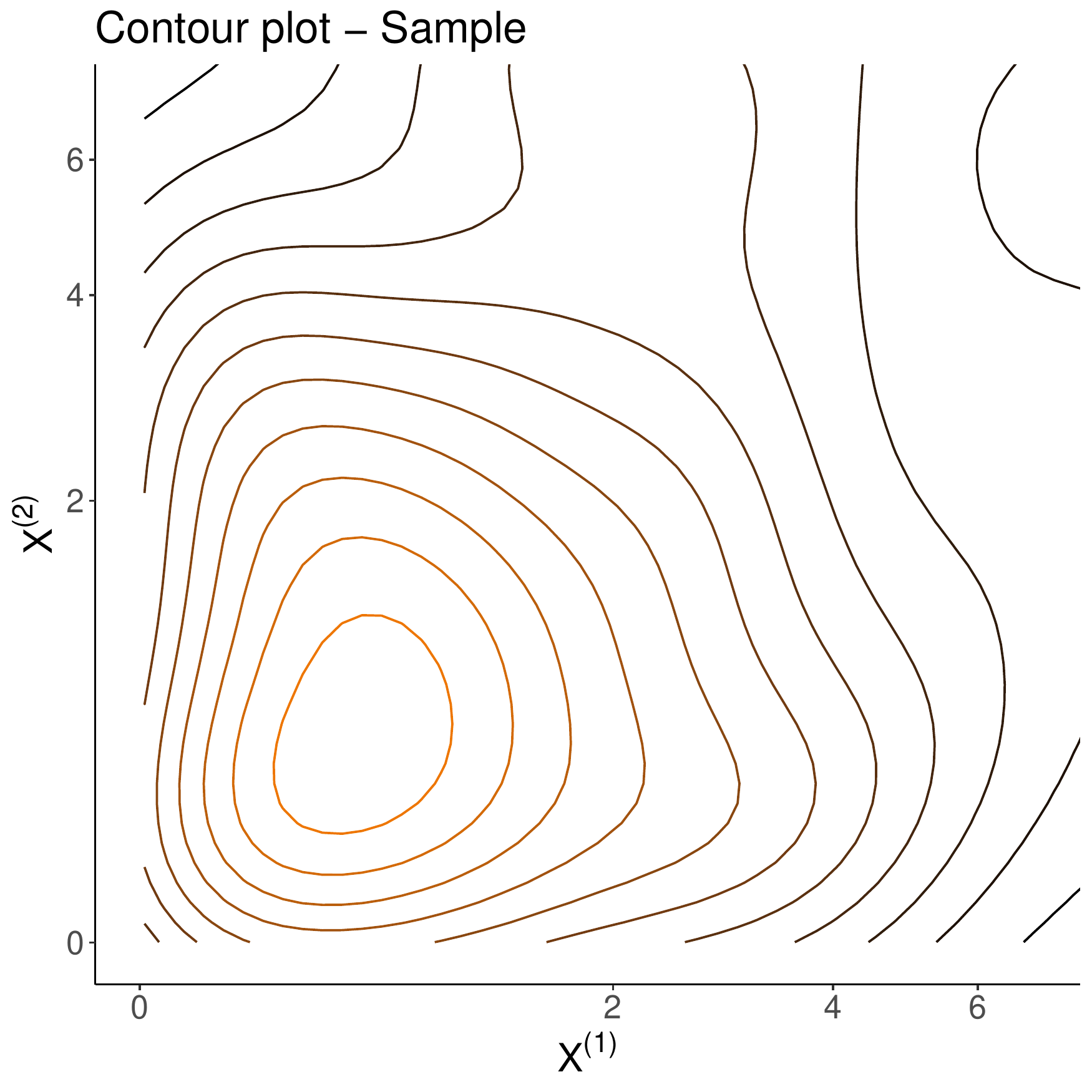}
  	\end{subfigure}
  	\begin{subfigure}{0.22\textwidth}
  		\includegraphics[width=\textwidth]{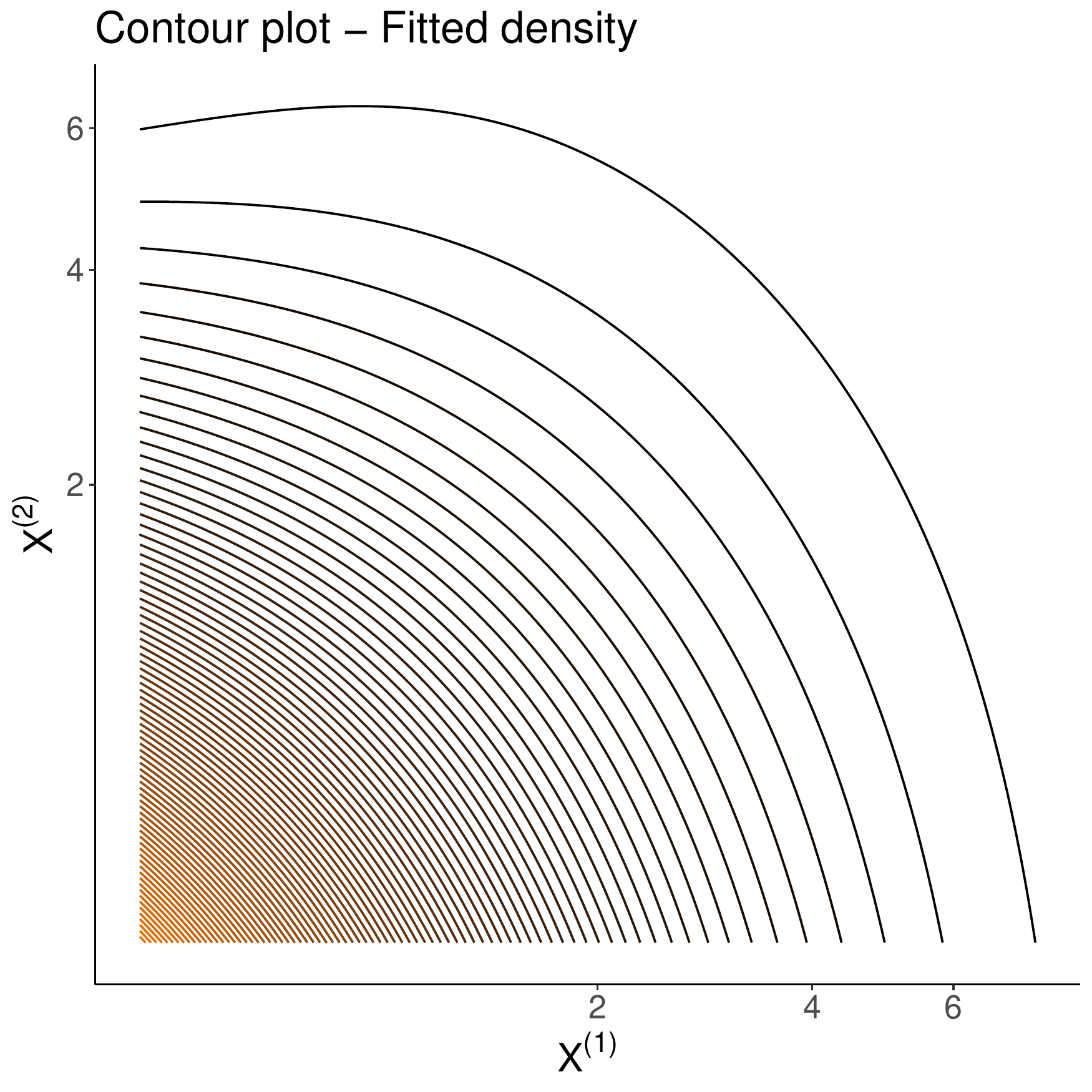}
  	\end{subfigure}
  	\begin{subfigure}{0.22\textwidth}
  		\includegraphics[width=\textwidth]{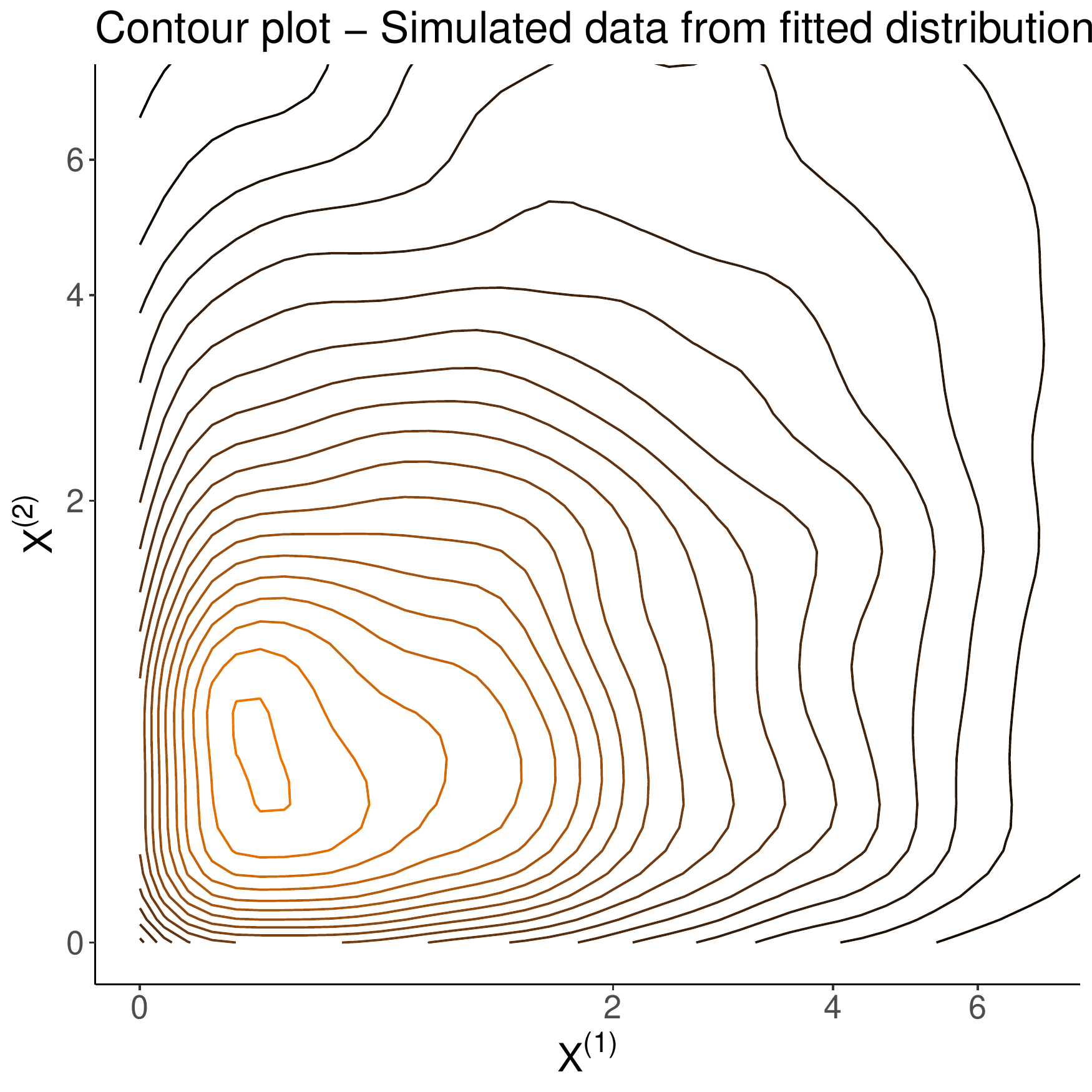}
  	\end{subfigure}
  \caption{Contour plot of the Danish fire insurance claim size sample (left),  contour plot of the Matrix--Pareto distribution fitted using Algorithm~\ref{alg:BivIPH} (middle) and contour plot of the simulated sample from the fitted distribution (right). }\label{fig:DanishContour}
\end{figure}

This bivariate data set was recently also studied in \cite[Sec. 4.5.2]{albrecher2017reinsurance}, where a splicing model with a bivariate mixed Erlang for the body and a bivariate Generalized Pareto distribution (GPD) for the tail was proposed. That approach required a threshold selection for the fitting of the tails, and univariate extreme value analysis led there to  values of regular variation of around $2$ for the building component and  $1.67$ for the contents component. Even though our estimates are further away from these values than their bivariate model ($1.75$ and $1.54$, respectively), we would like to emphasize that the fitting of a matrix--Pareto distribution does not require any threshold selection. Furthermore, if we were to use more phases and a general form of the sub--intensity and reward matrices, the fit would quickly improve and for about 6 phases reach the accuracy of the bivariate GPD model, but then the overall number of parameters compared to the size of the present data set may not be considered commensurate, which is why we stick to the above choice.
Note that our proposed procedure is fully automatic and the respective implementation can easily be applied to any other data set as well. 
\end{example}

\subsection{Multivariate Matrix--Weibull models}
Let $\bfX = (g_1(Y^{(1)}), \dots, g_d(Y^{(d)}))^{\prime} $,  where $\bfY \sim 	\mbox{MPH}^{*}\left( \bfpi, \bfT, \bfR \right)$ and  $g_j(y)=y^{1/\beta_j}$, $\beta_j>0$, $j=1,\dots,d$, then we say that $\bfX$  has a \textit{multivariate matrix--Weibull distribution}. Some special properties of this type of distribution are: 

\begin{enumerate}
	\item Marginal distributions are matrix--Weibull distributed.
	\item For a vector $\boldsymbol{a} = (a^{(1)},\dots,a^{(d)})$, with $a^{(j)}>0$, $j=1,\dots,d$,   $\Delta(\boldsymbol{a}) \bfX $ is multivariate matrix--Weibull distributed.
\end{enumerate}
For the bivariate case we get 
\begin{align*}
	f_{\bfX} \left( x^{(1)}, x^{(2)} \right)=\bfalp \ex^{\bfT_{11} (x^{(1)})^{\beta_1}} \bfT_{12}\ex^{\bfT_{22} (x^{(2)})^{\beta_2}}(-\bfT_{22})\bfe \, \beta_1 (x^{(1)})^{\beta_1 - 1} \beta_2 (x^{(2)})^{\beta_2 - 1} \,.
\end{align*}
and
\begin{align*}
		\bar{F}_{\bfX}( x^{(1)},x^{(2)} )=\bfalp \left( -\bfT_{11} \right)^{-1} \ex^{\bfT_{11} (x^{(1)})^{\beta_{1}}} \bfT_{12}\ex^{\bfT_{22} (x^{(2)})^{\beta_2}}\bfe \,
\end{align*}

\begin{remark} \rm 
	The Marshall--Olkin Weibull distribution (see \cite{hanagal1996multivariate}) is a particular case of this distribution. 
\end{remark} 

\subsubsection{Parameter estimation}
In contrast to Section~\ref{subsub:MPEst}, all transformations are parameter--dependent, and the fitting procedures of the previous subsection are not applicable. However, we can apply Algorithm~\ref{alg:BivIPH} in the bivariate case. 

\begin{example}\rm {(Bivariate Matrix--Weibull)}
	We generate an i.i.d.\ sample of size $5\,000$ of a bivariate random vector with matrix--Weibull marginals with parameters
\begin{gather*} 
 {\bfpi_1}=\left(
0.8,\ 0.2\right)\,, \\ 
{\bfT_1}=\left( \begin{array}{ccc}
-1 & 0.5  \\ 
 0 & -0.5 \\
\end{array} \right) \,, \\
\beta_1= 0.4 
\end{gather*}
for the first marginal and 
\begin{gather*} 
 {\bfpi_2}=\left(
0.5,\ 0.25,\ 0.25\right)\,, \\ 
{\bfT_2}=\left( \begin{array}{ccc}
-1 & 1 & 0 \\ 
 0 & -0.5 & 0.5  \\
0 & 0 & -0.1  \\
\end{array} \right) \,, \\ 
\beta_2= 0.6 
\end{gather*}
 for the second marginal, and a Gaussian copula with parameter $\rho=0.5$. While any copula, or also simply a bivariate matrix--Weibull based on a MPH$^*$ construction could be used, we choose the Gaussian copula here to illustrate that the algorithm is able to work with any type of dependence structure. This distribution has theoretical mean $\E(\bfX)=(18.7997, 86.3711)^{\prime}$.  The sample has numerical values  $\hat{\E}(\bfX)=(18.7690, 88.1637)^{\prime}$
 and $\hat{\rho}_{\tau} = 0.3431$.	

%
	
	We fit a bivariate matrix--Weibull distribution with $p_1=p_2=3$ using Algorithm~\ref{alg:BivIPH} with $1\,500$ steps (with a running time of $3\,930$ seconds for a step--length of $10^{-5}$), getting the following parameters: 
\begin{gather*} 
 \hat{\bfalp}=\left(
0.2101,\, 0,\, 0.7899,\, 0,\, 0,\, 0\right)\,, \\ 
\hat{\bfT}=\left( \begin{array}{cccccc}
-4.2507 & 0.5527 & 1.2916 & 0 & 0 & 2.4064 \\ 
0 & -0.3069 & 0 & 0 & 0.3069 & 0 \\
0.0089 & 0.2575 & -0.6903 & 0.4238 & 0 & 0  \\
 0 & 0 & 0  & -0.0946 & 0 & 0.0946 \\
0 & 0 & 0 & 0.0360 & -0.0360 & 0 \\
0 & 0 & 0 & 0.0026 & 0 & -0.2542
\end{array} \right) \,, \\ 
\hat{\bfR}= \left( \begin{array}{c c}
1 & 0 \\ 
1 & 0  \\
1 & 0  \\
0 & 1 \\ 
0 & 1 \\ 
0 & 1 
\end{array} \right)\,, \\
\beta_1= 0.4689 \,,\quad
\beta_2= 0.7340 \,.
\end{gather*}
One sees that the algorithm estimates the shape parameters of the matrix--Weibull marginals reasonably well. The fitted distribution has mean $\E(\bfX)=(18.6988, 88.1166)^{\prime}$, and from simulated data we get
  $\hat{\rho}_{\tau} = 0.3236$. The QQ and contour plots are given in Figure~\ref{fig:qqBivMW} and  Figure~\ref{fig:MWcontour2}, respectively. The log--likelihood of the fitted bivariate matrix--Weibull is $-39\,748$, which is to be compared with the log--likelihood $-39\,687.19$ using the original distribution.


\begin{figure}[H]
	\centering
  	\begin{subfigure}{0.22\textwidth}
  		\includegraphics[width=\textwidth]{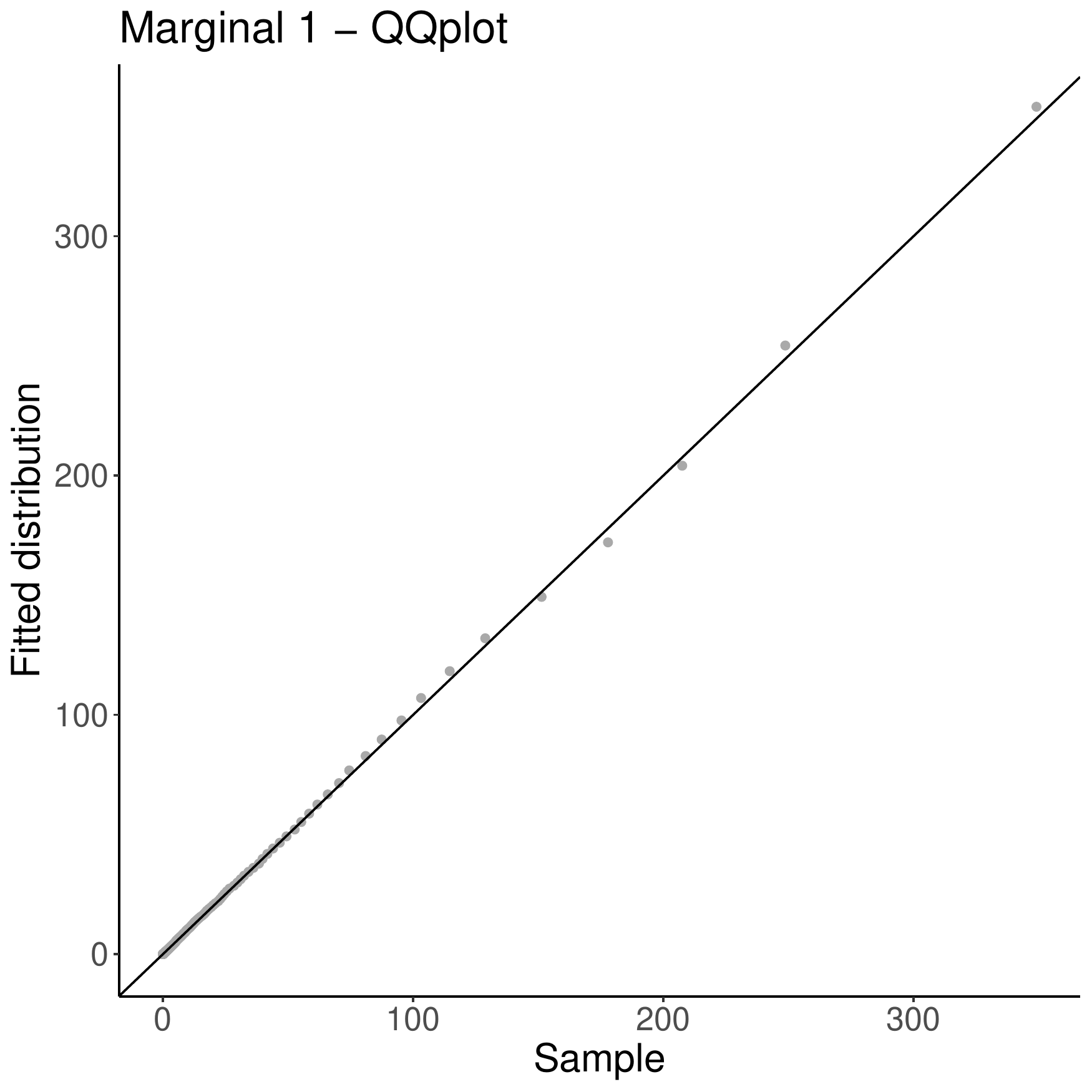}
  	\end{subfigure}
  	\begin{subfigure}{0.22\textwidth}
  		\includegraphics[width=\textwidth]{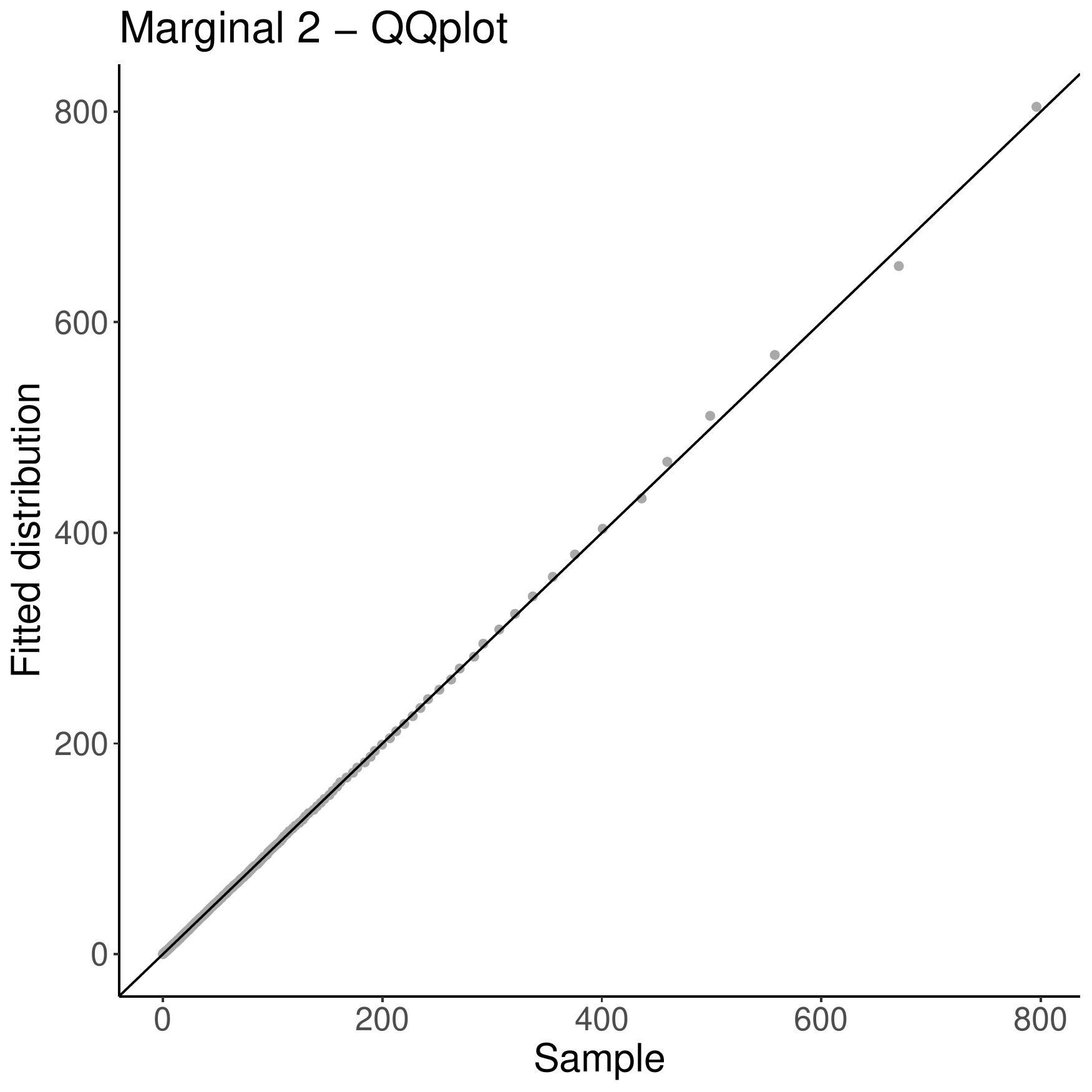}
  	\end{subfigure}
  \caption{QQ plots of sample versus fitted bivariate matrix--Weibull distribution using Algorithm~\ref{alg:BivIPH}. } \label{fig:qqBivMW}
\end{figure}

\begin{figure}[hbt]
	\centering
  	\begin{subfigure}{0.22\textwidth}
  		\includegraphics[width=\textwidth]{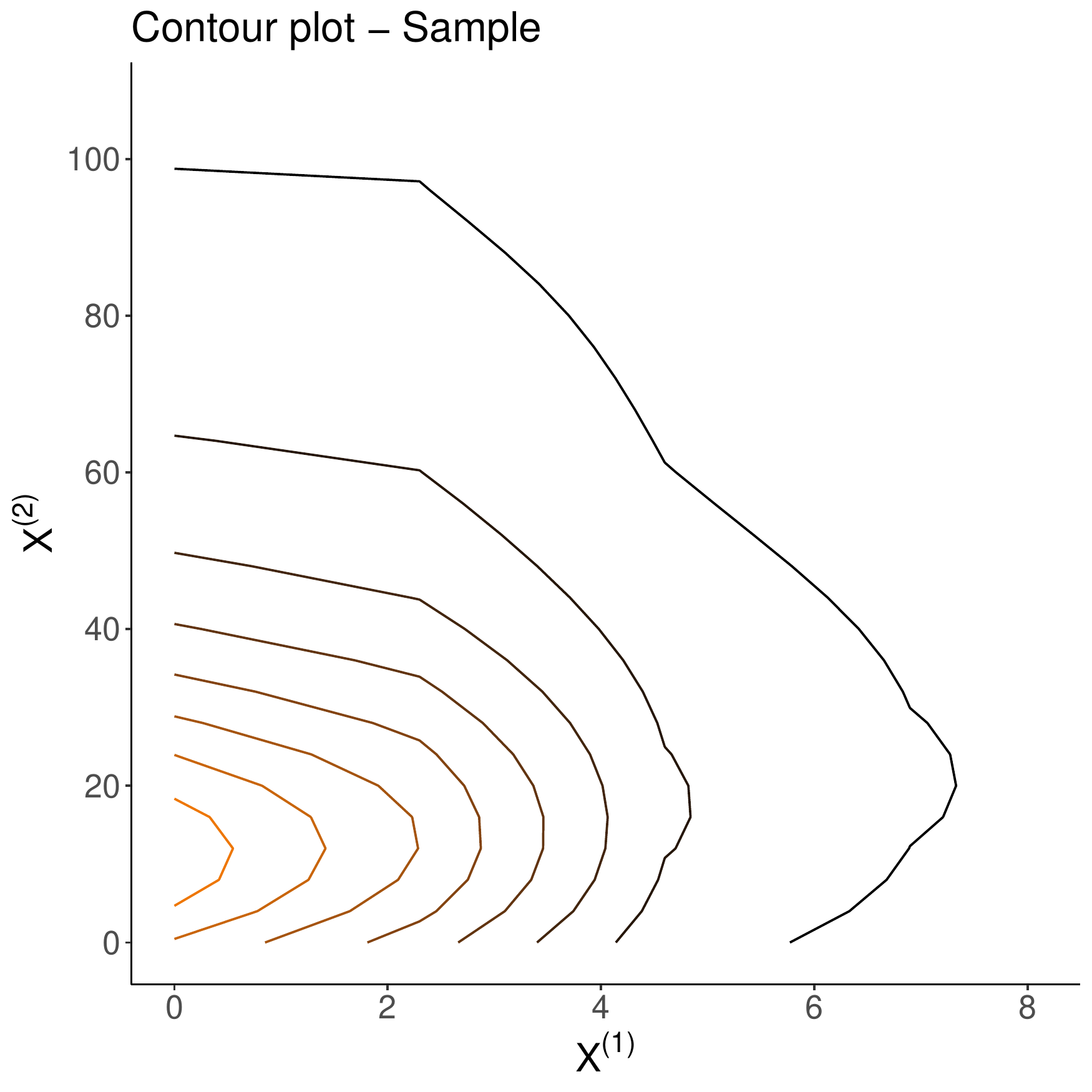}
  	\end{subfigure}
  	\begin{subfigure}{0.22\textwidth}
  		\includegraphics[width=\textwidth]{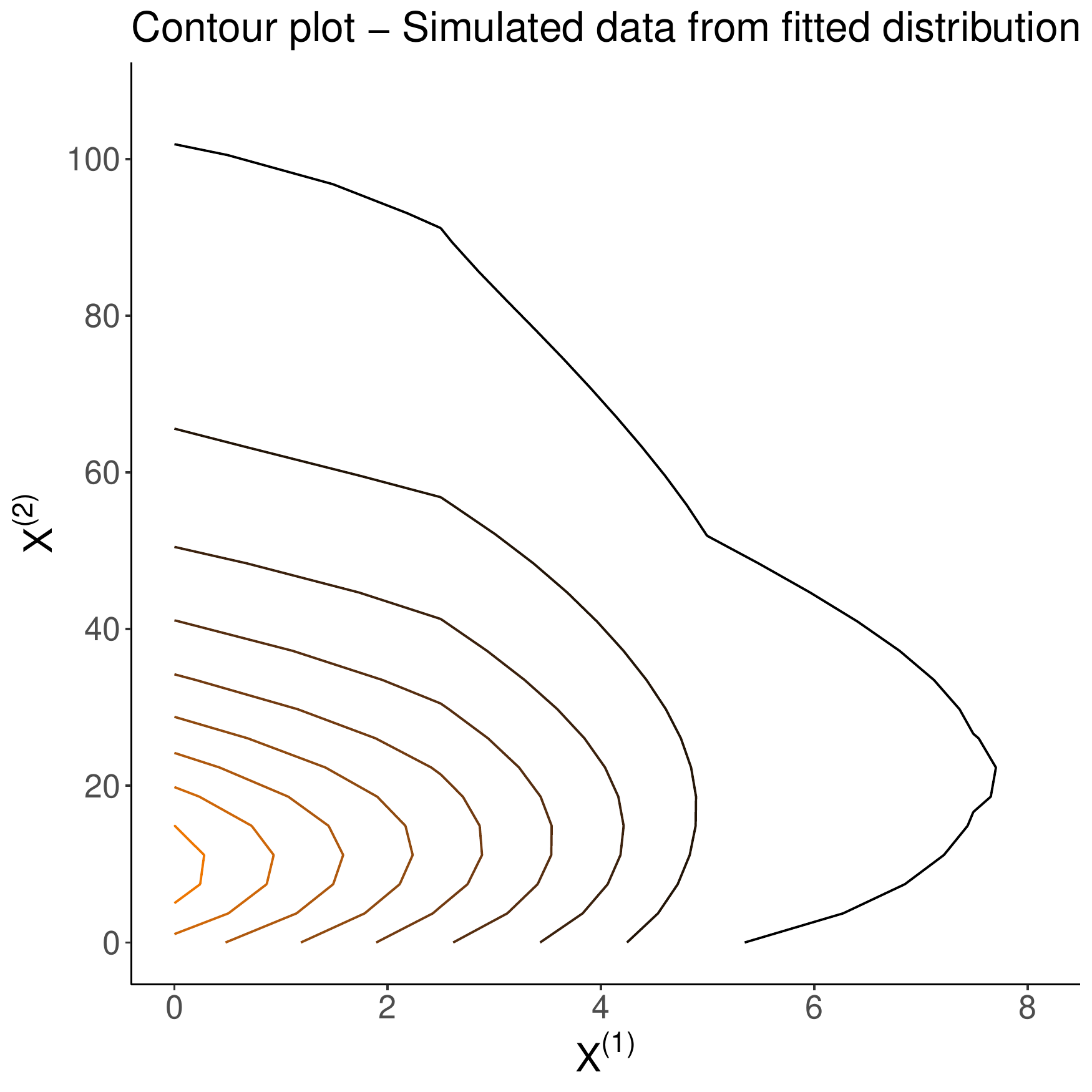}
  	\end{subfigure}
  \caption{Contour plot of simulated sample (left) and  contour plot of a simulated sample from the bivariate Matrix--Weibull distribution fitted using Algorithm~\ref{alg:BivIPH} (right). }\label{fig:MWcontour}
\end{figure}

\begin{figure}[hbt]
	\centering
  	\begin{subfigure}{0.22\textwidth}
  		\includegraphics[width=\textwidth]{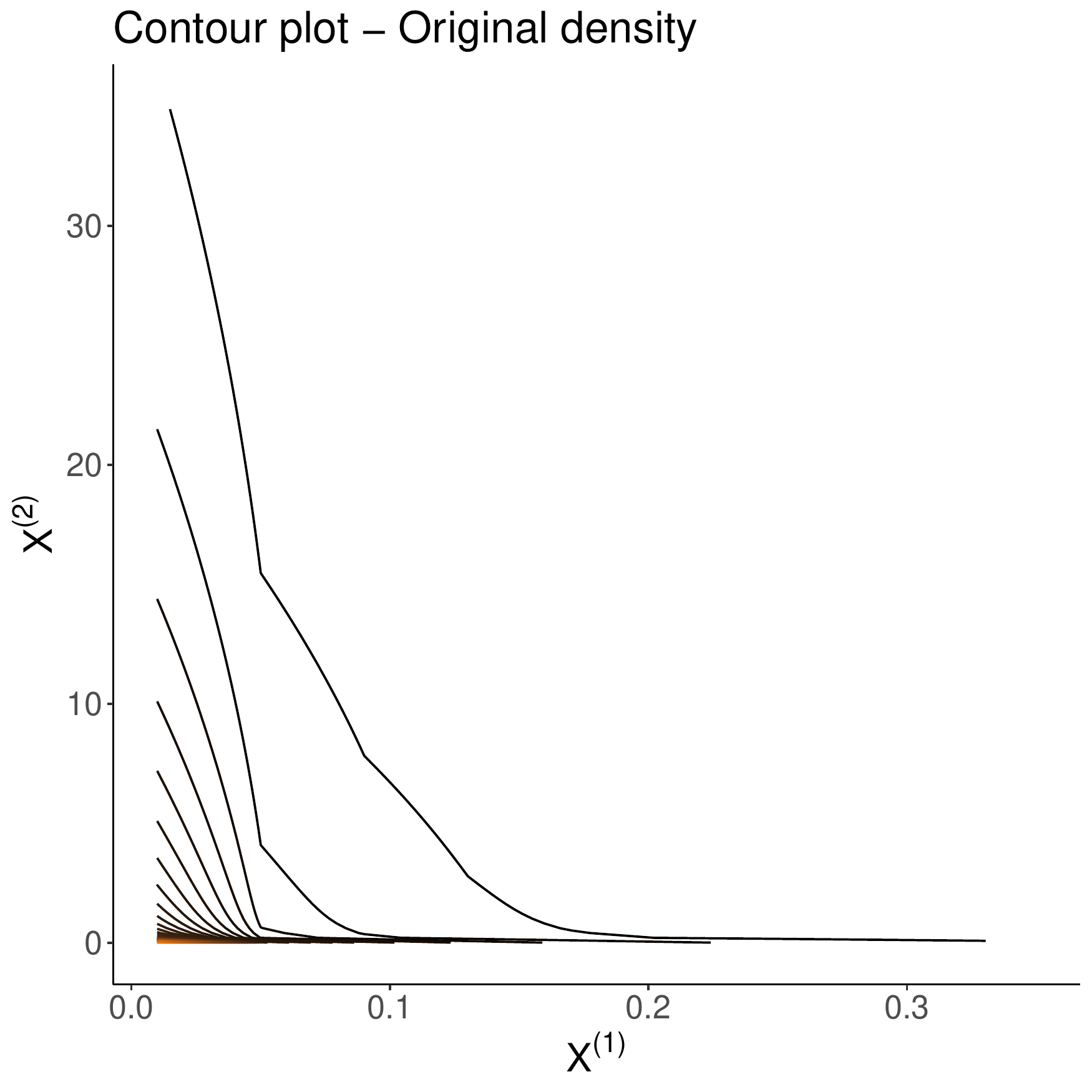}
  	\end{subfigure}
  	\begin{subfigure}{0.22\textwidth}
  		\includegraphics[width=\textwidth]{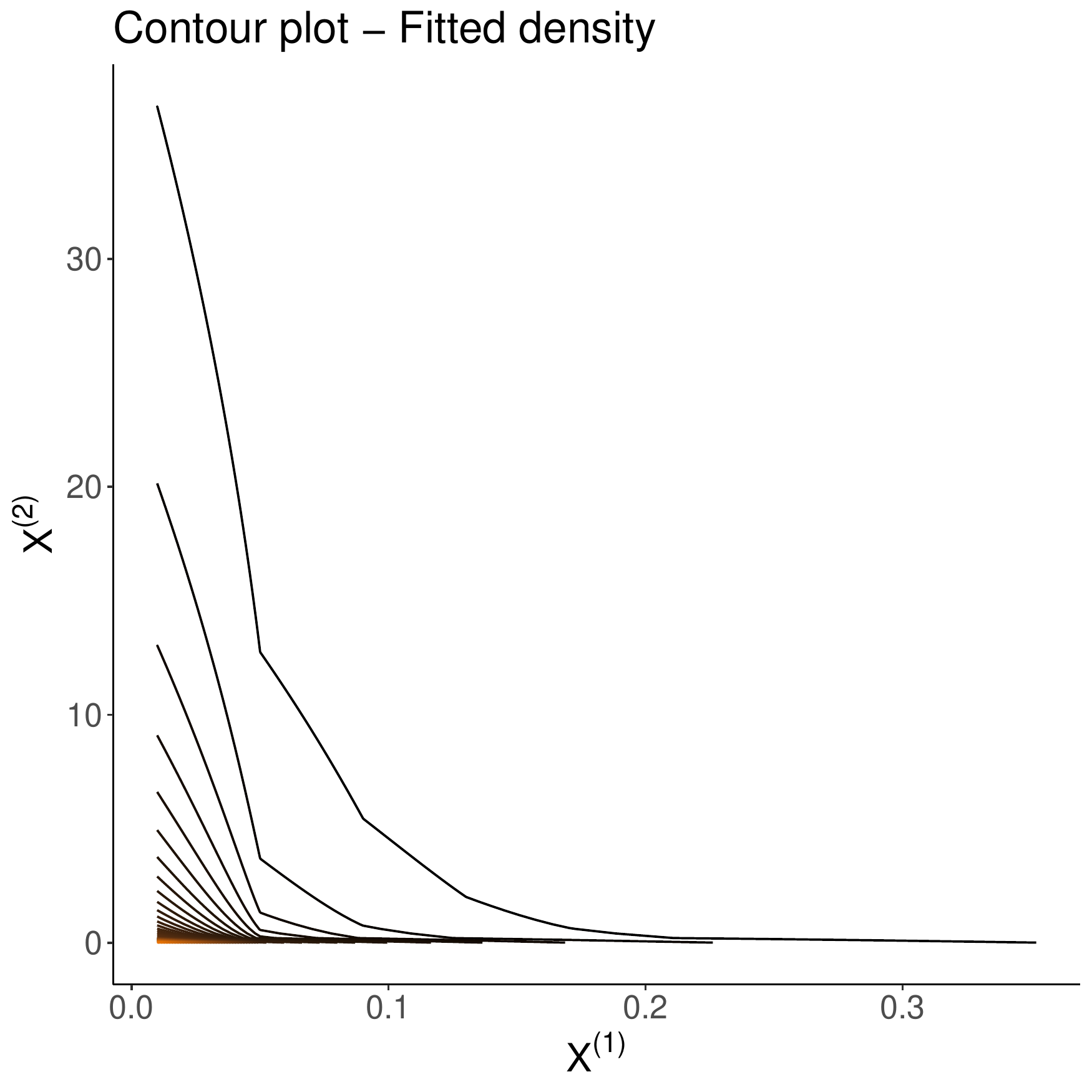}
  	\end{subfigure}
  \caption{Contour plot of original bivarite matrix-Weibull (left) and  contour plot of bivariate Matrix--Weibull distribution fitted using Algorithm~\ref{alg:BivIPH} (right). }\label{fig:MWcontour2}
\end{figure}
\end{example}

\begin{remark} \rm 
	In all examples of this section the marginals were assumed to be of the same type (both matrix--Pareto or both matrix--Weibull). We would like to mention that the generality of Algorithm~\ref{alg:BivIPH} also allows to fit models with marginals of different types (e.g.\ one marginal matrix--Pareto and the other matrix--Weibull). 
\end{remark}

\section{Conclusion}\label{sec:conclusions}
In this paper we provided a guide for the statistical fitting of homogeneous and inhomogeneous phase--type distributions to data, both for the univariate and multivariate case. For that purpose, we derived a new EM algorithm for IPH distributions that are obtained through parameter--dependent transformations. In addition, we  introduced new classes of multivariate distributions with IPH marginals and some attractive properties. As a by-product, we amended  the estimation method proposed by \cite{breuer2016semi} for the homogeneous MPH$^*$ case and illustrated its usefulness and flexibility. We furthermore discussed extensions for censored data and the fitting of the phase--type classes to given continuous joint distribution functions. The performance of the proposed algorithms was exemplified in various numerical examples, both on simulated and real data. In order to facilitate the implementation of the proposed algorithms for fitting this general class of distributions to given data, a respective R package is in preparation and will be made available on Cran. 


\section*{Acknowledgement}
We are grateful to Steffen L. Lauritzen for some important clarifications concerning the EM algorithm.  {\markcol We would like to thank two anonymous reviewers and the editor for the careful reading and constructive remarks.}

\bibliographystyle{apalike}
\bibliography{references.bib}

\end{document}